\input amstexl


\catcode`\@=11
\ifx\amstexloaded@\relax\else
 \errmessage{AmS-TeX must be loaded before LamS-TeX}\fi
\ifx\laxread@\undefined\else\catcode`\@=\active \fi
\def\err@#1{\errmessage{LamS-TeX error: #1}}
\def^^L{\par}
\let\+\tabalign
\def\newcount{\alloc@0\count\countdef\insc@unt}
\def\newdimen{\alloc@1\dimen\dimendef\insc@unt}
\def\newskip{\alloc@2\skip\skipdef\insc@unt}
\def\newmuskip{\alloc@3\muskip\muskipdef\@cclvi}
\def\newbox{\alloc@4\box\chardef\insc@unt}
\let\newtoks\relax
\def\newhelp#1#2{\newtoks#1#1\expandafter{\csname#2\endcsname}}
\def\newtoks{\alloc@5\toks\toksdef\@cclvi}
\def\newread{\alloc@6\read\chardef\sixt@@n}
\def\newwrite{\alloc@7\write\chardef\sixt@@n}
\def\newfam{\alloc@8\fam\chardef\sixt@@n}
\def\newlanguage{\alloc@9\language\chardef\@cclvi}
\def\newinsert#1{\global\advance\insc@unt by\m@ne
  \ch@ck0\insc@unt\count
  \ch@ck1\insc@unt\dimen
  \ch@ck2\insc@unt\skip
  \ch@ck4\insc@unt\box
  \allocationnumber=\insc@unt
  \global\chardef#1=\allocationnumber
  \wlog{\string#1=\string\insert\the\allocationnumber}}
\def\newif#1{\count@\escapechar \escapechar\m@ne
  \expandafter\expandafter\expandafter
   \edef\@if#1{true}{\let\noexpand#1=\noexpand\iftrue}%
  \expandafter\expandafter\expandafter
   \edef\@if#1{false}{\let\noexpand#1=\noexpand\iffalse}%
  \@if#1{false}\escapechar\count@}

\def\Err@#1{\errhelp\defaulthelp@\err@{#1}}
{\catcode`\@=\active
 \edef\next{\gdef\noexpand@{\futurelet\noexpand\next
  \csname at\string@\endcsname}}
 \next
}
\def\at@{\ifcat\noexpand\next a\let\next@\at@@\else
 \ifcat\noexpand\next0\let\next@\at@@\else
 \ifcat\noexpand\next\relax\let\next@\at@@\else
 \let\next@\at@@@\fi\fi\fi\next@}
\def\at@@@{\errhelp\athelp@\err@{Invalid use of @}}
\def\at@@#1{\expandafter
 \ifx\csname\string#1@at\endcsname\relax\let\next@\at@@@\else
 \DN@{\csname\string#1@at\endcsname}\fi\next@}
\def\atdef@#1{\expandafter\def\csname\string#1@at\endcsname}
\newif\iftest@
\def\tagin@#1{\tagin@false
 \DN@##1\tag##2##3\next@{\test@true\ifx\tagin@##2\test@false\fi}%
 \next@#1\tag\tagin@\next@\tagin@false\iftest@\tagin@true\fi}
\let\lkerns@\relax
\def\nolinebreak{\RIfM@\mathmodeerr@\nolinebreak\else
 \ifhmode\saveskip@\lastskip\unskip
 \nobreak\ifdim\saveskip@>\z@\hskip\saveskip@\fi\lkerns@
 \else\vmodeerr@\nolinebreak\fi\fi}
\def\allowlinebreak{\RIfM@\mathmodeerr@\allowlinebreak\else
 \ifhmode\saveskip@\lastskip\unskip
 \allowbreak\ifdim\saveskip@>\z@\hskip\saveskip@\fi\lkerns@
 \else\vmodeerr@\allowlinebreak\fi\fi}
\def\linebreak{\RIfM@\mathmodeerr@\linebreak\else
 \ifhmode\unskip\unkern\break\lkerns@
 \else\vmodeerr@\linebreak\fi\fi}
\let\nkerns@\relax
\def\newline{\RIfM@\mathmodeerr@\newline\else
 \ifhmode\unskip\unkern\null\hfill\break\nkerns@
 \else\vmodeerr@\newline\fi\fi}%
\def\newbox@{\alloc@@4\box\chardef\insc@unt}
\def\newcount@{\alloc@@0\count\countdef\insc@unt}
\def\accentedsymbol#1#2{\expandafter\newbox@\csname\exstring@#1@box\endcsname
 \setbox\csname\exstring@#1@box\endcsname\hbox{$\m@th#2$}%
 \define#1{\copy\csname\exstring@#1@box\endcsname{}}}
\def\rightadd@#1\to#2{\toks@{\\#1}\toks@@\expandafter{#2}\xdef#2{\the\toks@@
 \the\toks@}\toks@{}\toks@@{}}
\def\fontlist@{\\\tenrm\\\sevenrm\\\fiverm\\\teni\\\seveni\\\fivei
 \\\tensy\\\sevensy\\\fivesy\\\tenex\\\tenbf\\\sevenbf\\\fivebf
 \\\tensl\\\tenit}
\def\font@#1=#2 {\rightadd@#1\to\fontlist@\font#1=#2 }
\def\ismember@#1#2{\global\let\Next@ F\let\next@= #2%
 {\def\\##1{\let\nextii@##1\ifx\nextii@\next@\global\let\Next@ T\fi}#1}%
 \test@false\ifx\Next@ T\test@true\fi\let\next@\relax}
\def\FNSS@#1{\let\FNSS@@#1\FN@\FNSS@@@}
\def\FNSS@@@{\ifx\next\space@\def\FNSS@@@@. {\FN@\FNSS@@@}\else
 \def\FNSS@@@@.{\FNSS@@}\fi\FNSS@@@@.}
\atdef@"{\unskip
 \DN@{\ifx\next`\DN@`{\FN@\nextii@}%
  \else\ifx\next\lq\DN@\lq{\FN@\nextii@}%
  \else\DN@####1{\FN@\nextiii@}\fi\fi
  \next@}%
 \DNii@{\ifx\next`\DN@`{\sldl@``}%
  \else\ifx\next\lq\DN@\lq{\sldl@``}%
  \else\DN@{\dlsl@`}\fi\fi\next@}%
 \def\nextiii@{\ifx\next'\DN@'{\srdr@''}%
  \else\ifx\next\rq\DN@\rq{\srdr@''}%
  \else\DN@{\drsr@'}\fi\fi\next@}%
 \FNSS@\next@}
\def\root{%
  \DN@{\ifx\next\uproot\let\next@\nextii@\else
   \ifx\next\leftroot\let\next@\nextiii@\else
   \let\next@\plainroot@\fi\fi\next@}%
  \DNii@\uproot##1{\uproot@##1\relax\FNSS@\nextiv@}%
  \def\nextiv@{\ifx\next\leftroot\let\next@\nextv@\else
   \let\next@\plainroot@\fi\next@}%
  \def\nextv@\leftroot##1{\leftroot@##1\relax\plainroot@}%
  \def\nextiii@\leftroot##1{\leftroot@##1\relax\FNSS@\nextvi@}%
  \def\nextvi@{\ifx\next\uproot\let\next@\nextvii@\else
   \let\next@\plainroot@\fi\next@}%
  \def\nextvii@\uproot##1{\uproot@##1\relax\plainroot@}%
  \bgroup\uproot@\z@\leftroot@\z@
 \FNSS@\next@}
\def\loop#1\repeat{\def\iterate{#1\relax\expandafter\iterate\fi}%
 \iterate\let\iterate\relax}
\def\gloop@#1\repeat{\gdef\iterate@{#1\relax\expandafter\iterate@\fi}%
 \iterate@\global\let\iterate@\relax}
\def\printoptions{\W@{Do you want S(yntax check),
  G(alleys) or P(ages)?^^JType S, G or P, follow by <return>: }\loop
 \read\m@ne to\ans@
 \edef\next@{\def\noexpand\Ans@{\ans@}}\uppercase\expandafter{\next@}%
 \ifx\Ans@\S@\test@true\syntax\else
 \ifx\Ans@\G@\test@true\galleys\else
 \ifx\Ans@\P@\test@true\else
 \test@false\fi\fi\fi
 \iftest@\else\W@{Type S, G or P, follow by <return>: }%
 \repeat}
\expandafter\let\csname A@;\endcsname;
\expandafter\let\csname A@:\endcsname:
\expandafter\let\csname A@?\endcsname?
\expandafter\let\csname A@!\endcsname!
\def\APdef#1{\def\next@{\expandafter\let\csname A@\string#1\endcsname#1}%
 \afterassignment\next@\def#1}
\let\fextra@\,
\def\tdots@{\unskip
 \DN@{$\m@th\mathinner{\ldotp\ldotp\ldotp}\,
   \ifx\next,\,$\else\ifx\next.\,$\else
   \ifx\next;\,$\else
   \expandafter\ifx\csname A@\string;\endcsname\next\fextra@$\else
   \ifx\next:\,$\else
   \expandafter\ifx\csname A@\string:\endcsname\next\fextra@$\else
   \ifx\next?\,$\else
   \expandafter\ifx\csname A@\string?\endcsname\next\fextra@$\else
   \ifx\next!\,$\else
   \expandafter\ifx\csname A@\string!\endcsname\next\fextra@$\else
   $ \fi\fi\fi\fi\fi\fi\fi\fi\fi\fi}%
 \ \FN@\next@}
\def\extrap@#1{%
 \ifx\next,\DN@{#1\,}\else
 \ifx\next;\DN@{#1\,}\else
 \expandafter\ifx\csname A@\string;\endcsname\next\DN@{#1\fextra@}\else
 \ifx\next.\DN@{#1\,}\else\extra@
 \ifextra@\DN@{#1\,}\else
 \let\next@#1\fi\fi\fi\fi\fi\next@}
\def\dotsc{\DN@{\ifx\next;\plainldots@\,\else
 \expandafter\ifx\csname A@\string;\endcsname\next\plainldots@\fextra@\else
 \ifx\next.\plainldots@\,\else\extra@\plainldots@
 \ifextra@\,\fi\fi\fi\fi}%
 \FN@\next@}
\def\keybin@{\keybin@true
 \ifx\next+\else\ifx\next=\else\ifx\next<\else\ifx\next>\else\ifx\next-\else
 \ifx\next*\else\ifx\next:\else
 \expandafter\ifx\csname A@\string;\endcsname\next\else
 \keybin@false\fi\fi\fi\fi\fi\fi\fi\fi}
\def\boldkey#1{\ifcat\noexpand#1A%
  \ifcmmibloaded@{\fam\cmmibfam#1}\else
   \Err@{First bold symbol font not loaded}\fi
 \else
 \let\next=#1%
 \ifx#1!\mathchar"5\bffam@21 \else
 \expandafter\ifx\csname A@\string!\endcsname\next\mathchar"5\bffam@21 \else
 \ifx#1(\mathchar"4\bffam@28 \else\ifx#1)\mathchar"5\bffam@29 \else
 \ifx#1+\mathchar"2\bffam@2B \else\ifx#1:\mathchar"3\bffam@3A \else
 \expandafter\ifx\csname A@\string:\endcsname\next\mathchar"3\bffam@3A \else
 \ifx#1;\mathchar"6\bffam@3B \else
 \expandafter\ifx\csname A@\string;\endcsname\next\mathchar"6\bffam@3B \else
 \ifx#1=\mathchar"3\bffam@3D \else
 \ifx#1?\mathchar"5\bffam@3F \else
 \expandafter\ifx\csname A@\string?\endcsname\next\mathchar"5\bffam@3F \else
 \ifx#1[\mathchar"4\bffam@5B \else
 \ifx#1]\mathchar"5\bffam@5D \else
 \ifx#1,\mathchari@63B \else
 \ifx#1-\mathcharii@200 \else
 \ifx#1.\mathchari@03A \else
 \ifx#1/\mathchari@03D \else
 \ifx#1<\mathchari@33C \else
 \ifx#1>\mathchari@33E \else
 \ifx#1*\mathcharii@203 \else
 \ifx#1|\mathcharii@06A \else
 \ifx#10\bold0\else\ifx#11\bold1\else\ifx#12\bold2\else\ifx#13\bold3\else
 \ifx#14\bold4\else\ifx#15\bold5\else\ifx#16\bold6\else\ifx#17\bold7\else
 \ifx#18\bold8\else\ifx#19\bold9\else
  \Err@{\noexpand\boldkey can't be used with #1}%
 \fi\fi\fi\fi\fi\fi\fi\fi\fi\fi\fi\fi\fi\fi\fi
 \fi\fi\fi\fi\fi\fi\fi\fi\fi\fi\fi\fi\fi\fi\fi\fi\fi\fi}
\def\arabic#1{#1}
\def\alph#1{\count@#1\relax\advance\count@96 \ifnum\count@>122
 \Err@{\noexpand\alph invalid for numbers > 26}\else\char\count@\fi}
\def\Alph#1{\count@#1\relax\advance\count@64 \ifnum\count@>90
 \Err@{\noexpand\Alph invalid for numbers > 26}\else\char\count@\fi}

\def\Roman#1{\uppercase\expandafter{\romannumeral#1}}
\def\fnsymbol#1{\count@#1\relax
 \count@@\count@
 \advance\count@\m@ne\divide\count@7
 \count@@@\count@\advance\count@@@\@ne
 \multiply\count@7 \advance\count@@-\count@
 \count@\count@@@
 {\loop
  \ifcase\count@@\or*\or\dag\or\ddag\or\P\or\S\or\text{$\|$}\or\#\fi
  \advance\count@\m@ne\ifnum\count@>\z@\repeat}}
\def\cardnine@#1{\ifcase#1\or one\or two\or three\or four\or five\or
 six\or seven\or eight\or nine\fi}
\let\alloc@\alloc@@
\newcount\ten@
\ten@10
\def\cardinal#1{\count@#1\relax
 \ifnum\count@>99 \number\count@
 \else
  \ifnum\count@=\z@ zero%
  \else
   \ifnum\count@<\ten@\cardnine@\count@
   \else
    \ifnum\count@<20
     \advance\count@-\ten@
     \ifcase\count@ ten\or eleven\or twelve\or thirteen\or fourteen\or
      fifteen\or sixteen\or seventeen\or eighteen\or nineteen\fi
    \else
     \count@@\count@\count@@@\count@@
     \divide\count@\ten@\multiply\count@\ten@
     \advance\count@@@-\count@\divide\count@\ten@
     \ifcase\count@\or\or twenty\or thirty\or forty\or fifty\or sixty\or
      seventy\or eighty\or ninety\fi
     \ifnum\count@@@=\z@\else-\cardnine@\count@@@\fi
    \fi
   \fi
  \fi
 \fi}
\def\ordnine@#1{\ifcase#1\or first\or second\or third\or fourth\or fifth\or
 sixth\or seventh\or eighth\or ninth\fi}
\newcount\count@@@@
\def\ordsuffix@{\count@@@@\count@
 \divide\count@\ten@
 \count@@@\count@\count@@\count@
 \divide\count@@\ten@\multiply\count@@\ten@
 \advance\count@@@-\count@@
 \ifnum\count@@@=\@ne th%
 \else
  \count@@@\count@@@@
  \count@@\count@@@@
  \divide\count@@\ten@\multiply\count@@\ten@
  \advance\count@@@-\count@@
  \ifcase\count@@@ th\or st\or nd\or rd\else th\fi
 \fi}
\def\nordinal#1{\count@#1\relax\number\count@\ordsuffix@}
\def\spordinal#1{\count@#1\relax\number\count@$^{\text{\ordsuffix@}}$}
\def\ordinal#1{\count@#1\relax
 \ifnum\count@>99 \number\count@\ordsuffix@
 \else
   \ifnum\count@=\z@ zeroth%
  \else
    \ifnum\count@<\ten@\ordnine@\count@
    \else
     \ifnum\count@<20 \advance\count@-\ten@
      \ifcase\count@ tenth\or eleventh\or twelfth\or thirteenth\or
       fourteenth\or fifteenth\or sixteenth\or seventeenth\or eighteenth\or
       nineteenth\fi
     \else
      \count@@\count@
      \divide\count@\ten@\multiply\count@\ten@
      \count@@@\count@@\advance\count@@@-\count@
      \divide\count@\ten@
      \ifcase\count@\or\or twent\or thirt\or fort\or fift\or sixt\or sevent\or
       eight\or ninet\fi
      \ifnum\count@@@=\z@ ieth\else y-\ordnine@\count@@@\fi
     \fi
    \fi
  \fi
 \fi}
\font@\tensmc=cmcsc10
\textonlyfont@\smc\tensmc
\newtoks\noexpandtoks@
\noexpandtoks@{\let\arabic\relax\let\alph\relax\let\Alph\relax
 \let\Roman\relax\let\fnsymbol\relax\let\rm\relax
 \let\it\relax\let\bf\relax\let\sl\relax\let\smc\relax
 \let\/\relax\let\null\relax}
\def\noexpands@{\the\noexpandtoks@}
\def\Nonexpanding#1{\global\noexpandtoks@
 \expandafter{\the\noexpandtoks@\let#1\relax}}
\def\prevanish@{\saveskip@\z@\ifhmode\saveskip@\lastskip\unskip\fi}
\def\postvanish@{\ifdim\saveskip@>\z@\hskip\saveskip@\fi\FN@\postvanish@@}
\def\postvanish@@{\DN@.{}%
 \ifx\next\space@\ifdim\saveskip@>\z@\DN@. {}\fi\fi\next@.}
\def\invisible#1{\prevanish@\ignorespaces#1\unskip\postvanish@}
\def\vanishlist@{\\\invisible}
\let\noindent@\noindent
\def\noindent{\par\noindent@\FN@\pretendspace@}
\def\pretendspace@{\ismember@\vanishlist@\next
 \iftest@\nobreak\hskip-\p@\hskip\p@\fi}

\newtoks\everypartoks@
\def\noindent@@{\par\everypartoks@\expandafter{\the\everypar}\everypar{}%
 \noindent@\everypar\expandafter{\the\everypartoks@}}
\def\page{\Err@{\noexpand\page has no meaning by itself}}
\let\page@C\pageno
\let\page@P\empty
\let\page@Q\empty
\def\page@S#1{#1\/}
\def\page@F{\rm}
\def\page@N{\arabic}   
\newif\ifindexing@
\def\indexfile{\ifindexing@\else
 \alloc@@7\write\chardef\sixt@@n\ndx@
 \immediate\openout\ndx@=\jobname.ndx
 \global\indexing@true\fi}
\global\advance\insc@unt\m@ne
\ch@ck0\insc@unt\count
\ch@ck1\insc@unt\dimen
\ch@ck2\insc@unt\skip
\ch@ck4\insc@unt\box
\allocationnumber\insc@unt
\global\chardef\margin@\allocationnumber
\dimen\margin@\maxdimen
\count\margin@\z@
\skip\margin@\z@
\newif\ifindexproofing@
\def\indexproofing{\indexproofing@true}
\def\noindexproofing{\indexproofing@false}
\def\unmacro@#1:#2->#3\unmacro@{\def\macpar@{#2}\def\macdef@{#3}}
\def\starparts@#1{\def\stari@{#1}\def\starii@{#1}\let\stariii@\empty
 \test@false
 \DN@##1*##2##3\next@{\ifx\starparts@##2\test@false\else\test@true\fi}%
 \next@#1*\starparts@\next@
 \iftest@\DN@{\starparts@@#1\starparts@@}\else\let\next@\relax\fi\next@}
\def\starparts@@#1*#2\starparts@@{\def\starii@{#1}\def\stariii@{*#2}}
\def\windex@{\ifindexing@
 \expandafter\unmacro@\meaning\stari@\unmacro@
 \edef\macdef@{\string"\macdef@\string"}%
 \edef\next@{\write\ndx@{\macdef@}}\next@
 \write\ndx@{{\number\pageno}{\page@N}{\page@P}{\page@Q}}%
 \fi
 \ifindexproofing@
  \ifx\stariii@\empty\else
   \expandafter\unmacro@\meaning\stariii@\unmacro@\fi
  \insert\margin@{\hbox{\rm\vrule\height9\p@\depth2\p@\width\z@\starii@
  \ifx\stariii@\empty\else\tt\macdef@\fi}}\fi}
\catcode`\"=\active
\def"{\FN@\quote@}
\def\quote@{\ifx\next"\expandafter\quote@@\else\expandafter\quote@@@\fi}
\def\quote@@@#1"{\starparts@{#1}\starii@\windex@}
\def\quote@@"#1"{\prevanish@\starparts@{#1}\windex@\FN@\quote@@@@}
\def\quote@@@@{\ifx\next"\DN@"{\postvanish@}\else
 \let\next@\postvanish@\fi\next@}
\rightadd@"\to\vanishlist@
\def\idefine#1{\DN@{#1}\DNii@{\noexpand#1}%
 \afterassignment\idefine@\def\nextiii@}
\def\idefine@{\ifindexing@
 \expandafter\let\next@\nextiii@
 \expandafter\unmacro@\meaning\nextiii@\unmacro@
 \immediate\write\ndx@{\noexpand\define\nextii@\macpar@{\macdef@}}\fi}
\def\iabbrev*#1#2{\ifindexing@\toks@{#2}%
 \immediate\write\ndx@{\noexpand\abbrev*\noexpand#1{\the\toks@}}\fi}
\newread\laxread@
\newwrite\laxwrite@
\let\fnpages@\empty
\def\Finit@#1#2\Finit@{\let\nextii@#1\def\nextiii@{#2}}
\catcode`\~=11
\def\getparts@ @#1~#2~#3~#4~#5~#6{\def\nextiv@{#1}%
 \def\nextiii@{#2~#3~#4~#5~}\count@#6\relax}
\newif\ifdocument@
\def\document{\ifdocument@\else\global\document@true
 \let\fontlist@\empty
 \immediate\openin\laxread@=\jobname.lax\relax
 {\endlinechar\m@ne\noexpands@\catcode`\@=11 \catcode`\~=11
  \loop\ifeof\laxread@\else
   \read\laxread@ to\next@
   \ifx\next@\empty
   \else
    \expandafter\Finit@\next@\Finit@
    \if\nextii@ F%
     \expandafter\rightadd@\nextiii@\to\fnpages@
    \else
     \expandafter\getparts@\next@
     \edef\next@{\gdef\csname\nextiv@ @L\endcsname{\nextiii@\number\count@}}%
     \next@
    \fi
   \fi
  \repeat}%
 \immediate\closein\laxread@
 \immediate\openout\laxwrite@=\jobname.lax\relax\fi}
\let\thelabel@\relax
\def\thelabels@{\thelabel@ ~\thelabel@@ ~\thelabel@@@ ~\thelabel@@@@ ~}
\def\label#1{\prevanish@
 \ifx\thelabel@\relax
  \Err@{There's nothing here to be labelled}%
 \else
  {\noexpands@
  \expandafter\ifx\csname#1@L\endcsname\relax
   \expandafter\xdef\csname#1@L\endcsname{\thelabels@0}%
   \immediate\write\laxwrite@{@#1~\thelabels@1}%
  \else
   \edef\next@{@~\csname#1@L\endcsname}%
    \expandafter\getparts@\next@
    \ifodd\count@
    \expandafter\xdef\csname#1@L\endcsname{\thelabels@0}%
    \immediate\write\laxwrite@{@#1~\thelabels@1}%
   \else
    \Err@{Label #1 already used}%
   \fi
  \fi
  }%
 \fi
 \postvanish@}
\rightadd@\label\to\vanishlist@
\def\thepages@{\page@N{\number\page@C}~%
 \page@S{\page@P\page@N{\number\page@C}\page@Q}~%
 \number\page@C ~\page@P\page@N{\number\page@C}\page@Q ~}
\def\pagelabel#1{\prevanish@
 \expandafter\ifx\csname#1@L\endcsname\relax
  {\noexpands@
  \expandafter\xdef\csname#1@L\endcsname{\thepages@2}}%
  \write\laxwrite@{@#1~\thepages@3}%
 \else
  {\noexpands@
  \edef\next@{@~\csname#1@L\endcsname}%
  \expandafter\getparts@\next@
  \ifodd\count@
   \ifnum\count@=\@ne
    \expandafter\xdef\csname#1@L\endcsname{\thelabels@2}%
   \fi
   \write\laxwrite@{@#1~\thepages@3}%
  \else
   \Err@{Label #1 already used}%
  \fi
  }%
 \fi
 \postvanish@}
\rightadd@\pagelabel\to\vanishlist@
\newif\ifreferr@
\referr@true
\def\RefErrors{\global\referr@true}
\def\RefWarnings{\global\referr@false}
\setbox\z@\hbox{\global\count@=`^^30}
\ifnum\count@=48 \let\versionthree@\relax\fi
\def\nolabel@#1#2#3{\expandafter\ifx\csname#2@L\endcsname\relax
 \ifreferr@\Err@{No \noexpand\label found for #2}\else
 \W@{Warning: No \noexpand\label found for #2.}%
 \ifx\versionthree@\relax\W@{l.\number\inputlineno\space ... \string#1{#2}}\fi
 \fi#3\else}
\def\csL@#1{{\noexpands@\xdef\Next@{\csname#1@L\endcsname}}}
\def\ref#1{\nolabel@\ref{#1}\relax
 \DNii@##1~##2\nextii@{##1}%
 \csL@{#1}\expandafter\nextii@\Next@\nextii@\fi}
\def\Ref#1{\nolabel@\Ref{#1}\relax
 \DNii@##1~##2~##3\nextii@{##2}%
 \csL@{#1}\expandafter\nextii@\Next@\nextii@\fi}
\def\nref#1{\nolabel@\nref{#1}\relax
 \DNii@##1~##2~##3~##4\nextii@{##3}%
 \csL@{#1}\expandafter\nextii@\Next@\nextii@\fi}
\def\pref#1{\nolabel@\pref{#1}\relax
 \DNii@##1~##2~##3~##4~##5\nextii@{##4}%
 \csL@{#1}\expandafter\nextii@\Next@\nextii@\fi}
\let\pref@\pref
\def\Evaluatenref#1{\nolabel@\Evaluatenref{#1}{\gdef\Nref{-10000 }}%
 \DNii@##1~##2~##3~##4\nextii@{\DNii@{##3}}%
 \csL@{#1}\expandafter\nextii@\Next@\nextii@
 \xdef\Nref{\nextii@}\fi}
\def\Evaluatepref#1{\nolabel@\Evaluatepref{#1}{\global\let\Pref\empty}%
 \DNii@##1~##2~##3~##4~##5\nextii@{\DNii@{##4}}%
 \csL@{#1}\expandafter\nextii@\Next@\nextii@
 \xdef\Pref{\nextii@}\fi}
\def\readlax#1{\immediate\openin\laxread@=#1.lax\relax
 \ifeof\laxread@\W@{}\W@{File #1.lax not found.}\W@{}\fi
 {\endlinechar\m@ne\noexpands@\catcode`\@=11 \catcode`\~=11
  \loop\ifeof\laxread@\else
   \read\laxread@ to\nextv@
   \ifx\nextv@\empty
   \else
    \expandafter\Finit@\nextv@\Finit@
    \ifx\nextii@ F%
    \else
     \expandafter\getparts@\nextv@
     \expandafter\ifx\csname\nextiv@ @L\endcsname\relax
      \edef\next@{\gdef\csname\nextiv@ @L\endcsname
       {\nextiii@\ifnum\count@=\@ne0\else2\fi}}%
      \next@
     \else
      \Err@{Label \nextiv@\space in #1.lax already used}%
     \fi
    \fi
   \fi
  \repeat}%
 \immediate\closein\laxread@}
\catcode`\~=\active

\def\FNSSP@{\FNSS@\pretendspace@}
\everydisplay{\csname displaymath \endcsname}
\expandafter\def\csname displaymath \endcsname#1$${#1$$\FNSSP@}
\def\locallabel@{\let\thelabel@\Thelabel@\let\thelabel@@\Thelabel@@
 \let\thelabel@@@\Thelabel@@@\let\thelabel@@@@\Thelabel@@@@}
\newcount\tag@C
\tag@C\z@
\let\tag@P\empty
\let\tag@Q\empty
\def\tag@S#1{{\rm(}{#1\/}{\rm)}}
\let\tag@N\arabic
\def\tag@F{\rm}
\def\maketag@{\FN@\maketag@@}
\def\maketag@@{\ifx\next\relax\DN@\relax{\FN@\maketag@@}\else
 \ifx\next"\let\next@\maketag@@@\else
 \let\next@\maketag@@@@\fi\fi\next@}
\def\xdefThelabel@#1{\xdef\Thelabel@{#1{\Thelabel@@@}}}
\def\xdefThelabel@@#1{\xdef\Thelabel@@{#1{\Thelabel@@@@}}}
\def\maketag@@@@#1\maketag@{\global\advance\tag@C\@ne
 {\noexpands@
  \xdef\Thelabel@@@{\number\tag@C}%
  \xdefThelabel@\tag@N
  \xdef\Thelabel@@@@{\ifmathtags@$\tag@P\Thelabel@\tag@Q$\else
   \tag@P\Thelabel@\tag@Q\fi}%
  \xdefThelabel@@\tag@S
  }%
 \locallabel@
 \hbox{\tag@F\thelabel@@}%
 #1}
\def\Qlabel@#1{{\noexpands@\xdef\Thelabel@@{#1}%
 \let\style\empty\xdef\Thelabel@@@@{#1}%
 \let\pre\empty\let\post\empty\xdef\Thelabel@{#1}%
 \let\numstyle\empty\xdef\Thelabel@@@{#1}}}
\def\maketag@@@"#1"#2\maketag@{%
 {\let\pre\tag@P\let\post\tag@Q\let\style\tag@S\let\numstyle\tag@N
  \hbox{\tag@F#1}%
  \noexpands@
  \Qlabel@{#1}%
  }%
 \locallabel@
 #2}
\def\align@{\inalign@true\inany@true
 \vspace@\allowdisplaybreak@\displaybreak@\intertext@
 \def\tag{\global\tag@true\ifnum\and@=\z@
  \DN@{&\omit\global\rwidth@\z@&\relax}\else
  \DN@{&\relax}\fi\next@}%
 \iftagsleft@\DN@{\csname align \endcsname}\else
  \DN@{\csname align \space\endcsname}\fi\next@}
\def\noset@{\def\Offset##1##2{\prevanish@\postvanish@}%
 \def\Reset##1##2{\prevanish@\postvanish@}}
\def\measure@#1\endalign{\global\lwidth@\z@\global\rwidth@\z@
 \global\maxlwidth@\z@\global\maxrwidth@\z@
 \global\and@\z@
 \setbox\z@\vbox
  {\noset@\everycr{\noalign{\global\tag@false\global\and@\z@}}\Let@
  \halign{\setboxz@h{$\m@th\displaystyle{\@lign##}$}%
   \global\lwidth@\wdz@
   \ifdim\lwidth@>\maxlwidth@\global\maxlwidth@\lwidth@\fi
   \global\advance\and@\@ne
   &\setboxz@h{$\m@th\displaystyle{{}\@lign##}$}\global\rwidth@\wdz@
   \ifdim\rwidth@>\maxrwidth@\global\maxrwidth@\rwidth@\fi
   \global\advance\and@\@ne
   &\Tag@\eat@{##}\crcr#1\crcr}}%
 \totwidth@\maxlwidth@\advance\totwidth@\maxrwidth@}
\def\prepost@{\global\let\tag@P@\tag@P\global\let\tag@Q@\tag@Q}
\def\reprepost@{\let\tag@P\tag@P@\let\tag@Q\tag@Q@}
\expandafter\def\csname align \space\endcsname#1\endalign
 {\measure@#1\endalign\global\and@\z@
 \ifingather@\everycr{\noalign{\global\and@\z@}}\else\displ@y@\fi
 \Let@\tabskip\centering@
 \halign to\displaywidth
  {\hfil\strut@\setboxz@h{$\m@th\displaystyle{\@lign##\prepost@}$}%
  \boxz@\global\advance\and@\@ne
  \tabskip\z@skip
  &\setboxz@h{$\m@th\displaystyle{{}\@lign##\prepost@}$}%
  \global\rwidth@\wdz@\boxz@\hfil\global\advance\and@\@ne
  \tabskip\centering@
  &\setboxz@h{\@lign\strut@\reprepost@\maketag@##\maketag@}%
  \dimen@\displaywidth\advance\dimen@-\totwidth@
  \divide\dimen@\tw@\advance\dimen@\maxrwidth@\advance\dimen@-\rwidth@
  \ifdim\dimen@<\tw@\wdz@\llap{\vtop{\normalbaselines\null\boxz@}}%
  \else\llap{\boxz@}\fi
  \tabskip\z@skip
  \crcr#1\crcr
  \black@\totwidth@}}
\expandafter\def\csname align \endcsname#1\endalign{\measure@#1\endalign
 \global\and@\z@
 \ifdim\totwidth@>\displaywidth\let\displaywidth@\totwidth@\else
  \let\displaywidth@\displaywidth\fi
 \ifingather@\everycr{\noalign{\global\and@\z@}}\else\displ@y@\fi
 \Let@\tabskip\centering@\halign to\displaywidth
  {\hfil\strut@\setboxz@h{$\m@th\displaystyle{\@lign##\prepost@}$}%
  \global\lwidth@\wdz@\global\lineht@\ht\z@
  \boxz@\global\advance\and@\@ne
  \tabskip\z@skip&\setboxz@h{$\m@th\displaystyle{{}\@lign##\prepost@}$}%
  \ifdim\ht\z@>\lineht@\global\lineht@\ht\z@\fi
  \boxz@\hfil\global\advance\and@\@ne
  \tabskip\centering@&\kern-\displaywidth@
  \setboxz@h{\@lign\strut@\reprepost@\maketag@##\maketag@}%
  \dimen@\displaywidth\advance\dimen@-\totwidth@
  \divide\dimen@\tw@\advance\dimen@\maxlwidth@\advance\dimen@-\lwidth@
  \ifdim\dimen@<\tw@\wdz@
   \rlap{\vbox{\normalbaselines\boxz@\vbox to\lineht@{}}}\else
   \rlap{\boxz@}\fi
  \tabskip\displaywidth@\crcr#1\crcr\black@\totwidth@}}
\def\attag@#1{\let\Maketag@\maketag@\let\TAG@\Tag@
 \let\Prepost@\prepost@\let\Reprepost@\reprepost@
 \let\Tag@\relax\let\maketag@\relax
 \let\prepost@\relax\let\reprepost@\relax
 \ifmeasuring@
  \def\llap@##1{\setboxz@h{##1}\hbox to\tw@\wdz@{}}%
  \def\rlap@##1{\setboxz@h{##1}\hbox to\tw@\wdz@{}}%
 \else\let\llap@\llap\let\rlap@\rlap\fi
 \toks@{\hfil\strut@
  $\m@th\displaystyle{\@lign\the\hashtoks@\prepost@}$%
  \tabskip\z@skip\global\advance\and@\@ne&
  $\m@th\displaystyle{{}\@lign\the\hashtoks@\prepost@}$\hfil
  \ifxat@\tabskip\centering@\fi\global\advance\and@\@ne}%
 \iftagsleft@
  \toks@@{\tabskip\centering@&\Tag@\kern-\displaywidth
   \rlap@{\@lign\reprepost@\maketag@\the\hashtoks@\maketag@}%
   \global\advance\and@\@ne\tabskip\displaywidth}\else
  \toks@@{\tabskip\centering@&\Tag@\llap@{\@lign\reprepost@\maketag@
   \the\hashtoks@\maketag@}\global\advance\and@\@ne\tabskip\z@skip}\fi
 \atcount@#1\relax\advance\atcount@\m@ne
 \loop\ifnum\atcount@>\z@
  \toks@\expandafter{\the\toks@&\hfil$\m@th\displaystyle{\@lign
  \the\hashtoks@\prepost@}$\global\advance\and@\@ne
  \tabskip\z@skip
  &$\m@th\displaystyle{{}\@lign\the\hashtoks@\prepost@}$\hfil\ifxat@
  \tabskip\centering@\fi\global\advance\and@\@ne}\advance\atcount@\m@ne
 \repeat
 \edef\preamble@{\the\toks@\the\toks@@}%
 \edef\preamble@@{\preamble@}%
 \let\maketag@\Maketag@\let\Tag@\TAG@
 \let\prepost@\Prepost@\let\reprepost@\Reprepost@}
\def\unlabel@{\def\label##1{\prevanish@\postvanish@}%
 \def\pagelabel##1{\prevanish@\postvanish@}}
\newcount\tag@CC
\expandafter\def\csname alignat \endcsname#1#2\endalignat
 {\inany@true\xat@false
 \def\tag{\global\tag@true
  \count@#1\relax\multiply\count@\tw@\advance\count@\m@ne
  \gdef\tag@{&}%
  \loop\ifnum\count@>\and@\xdef\tag@{&\omit\tag@}%
  \advance\count@\m@ne\repeat
  \tag@\relax}%
 \vspace@\allowdisplaybreak@\displaybreak@\intertext@
 \displ@y@\measuring@true\tag@CC\tag@C
 \setbox\savealignat@\hbox{\noset@\unlabel@$\m@th\displaystyle\Let@
  \attag@{#1}\vbox{\halign{\span\preamble@@\crcr#2\crcr}}$}%
 \measuring@false
 \Let@\attag@{#1}\tag@C\tag@CC
 \tabskip\centering@\halign to\displaywidth
  {\span\preamble@@\crcr#2\crcr\black@{\wd\savealignat@}}}
\expandafter\def\csname xalignat \endcsname#1#2\endxalignat
 {\inany@true\xat@true
 \def\tag{\global\tag@true
  \count@#1\relax\multiply\count@\tw@\advance\count@\m@ne
  \gdef\tag@{&}%
  \loop\ifnum\count@>\and@\xdef\tag@{&\omit\tag@}%
  \advance\count@\m@ne\repeat
  \tag@\relax}%
 \vspace@\allowdisplaybreak@\displaybreak@\intertext@
 \displ@y@\measuring@true\tag@CC\tag@C
 \setbox\savealignat@\hbox{\noset@\unlabel@$\m@th\displaystyle\Let@
  \attag@{#1}\vbox{\halign{\span\preamble@@\crcr#2\crcr}}$}%
 \measuring@false\Let@\attag@{#1}\tag@C\tag@CC
 \tabskip\centering@\halign to\displaywidth
 {\span\preamble@@\crcr#2\crcr\black@{\wd\savealignat@}}}
\def\gather{\RIfMIfI@\DN@{\onlydmatherr@\gather}\else
 \ingather@true\inany@true\def\tag{&\relax}%
 \vspace@\allowdisplaybreak@\displaybreak@\intertext@
 \displ@y\Let@
 \iftagsleft@\DN@{\csname gather \endcsname}\else
  \DN@{\csname gather \space\endcsname}\fi\fi
 \else\DN@{\onlydmatherr@\gather}\fi\next@}
\def\exstring@{\expandafter\eat@\string}
\def\newcounter#1{\define#1{}%
 \edef\next@{\def\noexpand#1{\futurelet\noexpand\next
  \csname\exstring@#1@Z\endcsname}}\next@
 \edef\next@{\def\csname\exstring@#1@Z\endcsname
  {\global\advance\csname\exstring@#1@C\endcsname\@ne
  {\csname\exstring@#1@F\endcsname\csname\exstring@#1@S\endcsname
   {\csname\exstring@#1@P\endcsname\csname\exstring@#1@N\endcsname
   {\noexpand\number\csname\exstring@#1@C\endcsname}%
   \csname\exstring@#1@Q\endcsname}}%
  \noexpand\ifx\noexpand\next\noexpand\label
   \def\noexpand\next@\noexpand\label########1{{\noexpand\noexpands@
    \xdef\noexpand\Thelabel@{\csname\exstring@#1@N\endcsname
     {\noexpand\number\csname\exstring@#1@C\endcsname}}%
    \xdef\noexpand\Thelabel@@@{\noexpand\number
     \csname\exstring@#1@C\endcsname}%
    \xdef\noexpand\Thelabel@@{\csname\exstring@#1@S\endcsname
     {\csname\exstring@#1@P\endcsname
     \csname\exstring@#1@N\endcsname
     {\noexpand\number\csname\exstring@#1@C\endcsname}%
     \csname\exstring@#1@Q\endcsname}}%
    \xdef\noexpand\Thelabel@@@@{\csname\exstring@#1@P\endcsname
     \csname\exstring@#1@N\endcsname
     {\noexpand\number\csname\exstring@#1@C\endcsname}%
     \csname\exstring@#1@Q\endcsname}}%
    {\noexpand\locallabel@\noexpand\label{########1}}}%
   \noexpand\else\let\noexpand\next@\relax\noexpand\fi\noexpand\next@}}\next@
 \expandafter\newcount@\csname\exstring@#1@C\endcsname
 \expandafter\let\csname\exstring@#1@N\endcsname\arabic
 \expandafter\def\csname\exstring@#1@S\endcsname##1{##1\/}%
 \expandafter\let\csname\exstring@#1@P\endcsname\empty
 \expandafter\let\csname\exstring@#1@Q\endcsname\empty
 \expandafter\def\csname\exstring@#1@F\endcsname{\rm}%
 }
\def\HASH@#1#2{\ifnum#2=\z@\else
 \edef\next@{\toks@{\the\toks@\the\hashtoks@#2}%
 \toks@@{\the\toks@@{\the\hashtoks@#2}}}\next@\expandafter\HASH@\fi}
\def\HASH@@{\toks@{}\toks@@{}\expandafter\HASH@\macpar@00}
\def\usecounter#1#2{\expandafter\ifx\csname\exstring@#1@Z\endcsname
 \relax\Err@{\noexpand#1not created with \string\newcounter}\fi
 \expandafter\let\csname\exstring@#1@@Z\endcsname\relax
 \expandafter\let\csname\exstring@#1@@Z@\endcsname\relax
 \expandafter\let\csname\exstring@#1@@Z@@\endcsname\relax
 \edef\next@{\def\noexpand#2{\futurelet\noexpand\next
  \csname\exstring@#1@@Z\endcsname}}\next@
 \edef\next@{\def\csname\exstring@#1@@Z\endcsname{\noexpand\ifx
  \noexpand\next\noexpand\label\def\noexpand\next@\noexpand\label
   ########1{\csname\exstring@#1@@Z@\endcsname
   {\noexpand#1\noexpand\label{########1}}}%
   \noexpand\else\noexpand\ifx\noexpand\next
   \noexpand"\def\noexpand\next@\noexpand"########1\noexpand"%
   {\csname\exstring@#1@@Z@\endcsname{{\expandafter\noexpand
   \csname\exstring@#1@F\endcsname
   \let\noexpand\pre\expandafter\noexpand\csname\exstring@#1@P\endcsname
   \let\noexpand\post\expandafter\noexpand\csname\exstring@#1@Q\endcsname
   \let\noexpand\style\expandafter\noexpand\csname\exstring@#1@S\endcsname
   \let\noexpand\numstyle\expandafter\noexpand\csname\exstring@#1@N\endcsname
   ########1}}}\noexpand\else
   \def\noexpand\next@{\csname\exstring@#1@@Z@\endcsname{\noexpand#1}}%
   \noexpand\fi\noexpand\fi\noexpand\next@}}\next@
 \def\next@{\expandafter\expandafter\expandafter\unmacro@\expandafter
  \meaning\csname\exstring@#1@@Z@@\endcsname\unmacro@
  \HASH@@
  \edef\next@{\def\csname\exstring@#1@@Z@\endcsname\the\toks@{%
   \expandafter\noexpand\csname\exstring@#1@@Z@@\endcsname\the\toks@@
   \noexpand\FNSSP@}}\next@}%
 \afterassignment\next@
 \expandafter\def\csname\exstring@#1@@Z@@\endcsname}
\def\listbi@{\penalty50 \medskip}
\def\listbii@{\penalty100 \smallskip}
\let\listbiii@\relax
\let\listbiv@\relax
\let\listbv@\relax
\def\listmi@{\advance\leftskip30\p@\relax}
\let\listmii@\listmi@
\let\listmiii@\listmi@
\let\listmiv@\listmi@
\let\listmv@\listmi@
\def\itemi@#1{\noindent@@\llap{#1\hskip5\p@}}
\let\itemii@\itemi@
\let\itemiii@\itemi@
\let\itemiv@\itemi@
\let\itemv@\itemi@
\def\liste@{\penalty-50 \medskip}
\def\listei@{\penalty-100 \smallskip}
\let\listeii@\relax
\let\listeiii@\relax
\let\listeiv@\relax
\expandafter\newcount\csname list@C1\endcsname
\csname list@C1\endcsname\z@
\expandafter\newcount\csname list@C2\endcsname
\csname list@C2\endcsname\z@
\expandafter\newcount\csname list@C3\endcsname
\csname list@C3\endcsname\z@
\expandafter\newcount\csname list@C4\endcsname
\csname list@C4\endcsname\z@
\expandafter\newcount\csname list@C5\endcsname
\csname list@C5\endcsname\z@
\expandafter\let\csname list@P1\endcsname\empty
\expandafter\let\csname list@P2\endcsname\empty
\expandafter\let\csname list@P3\endcsname\empty
\expandafter\let\csname list@P4\endcsname\empty
\expandafter\let\csname list@P5\endcsname\empty
\expandafter\let\csname list@Q1\endcsname\empty
\expandafter\let\csname list@Q2\endcsname\empty
\expandafter\let\csname list@Q3\endcsname\empty
\expandafter\let\csname list@Q4\endcsname\empty
\expandafter\let\csname list@Q5\endcsname\empty
\expandafter\def\csname list@S1\endcsname#1{{\rm(}{#1\/}{\rm)}}
\expandafter\def\csname list@S2\endcsname#1{{\rm(}{#1\/}{\rm)}}
\expandafter\def\csname list@S3\endcsname#1{{\rm(}{#1\/}{\rm)}}
\expandafter\def\csname list@S4\endcsname#1{{\rm(}{#1\/}{\rm)}}
\expandafter\def\csname list@S5\endcsname#1{{\rm(}{#1\/}{\rm)}}
\expandafter\let\csname list@N1\endcsname\arabic
\expandafter\let\csname list@N2\endcsname\arabic
\expandafter\let\csname list@N3\endcsname\arabic
\expandafter\let\csname list@N4\endcsname\arabic
\expandafter\let\csname list@N5\endcsname\arabic
\expandafter\def\csname list@F1\endcsname{\rm}
\expandafter\def\csname list@F2\endcsname{\rm}
\expandafter\def\csname list@F3\endcsname{\rm}
\expandafter\def\csname list@F4\endcsname{\rm}
\expandafter\def\csname list@F5\endcsname{\rm}
\newcount\listlevel@
\listlevel@\z@
\def\list@@C{\csname list@C\number\listlevel@\endcsname}
\def\list@@P{\csname list@P\number\listlevel@\endcsname}
\def\list@@Q{\csname list@Q\number\listlevel@\endcsname}
\def\list@@S{\csname list@S\number\listlevel@\endcsname}
\def\list@@N{\csname list@N\number\listlevel@\endcsname}
\def\list@@F{\csname list@F\number\listlevel@\endcsname}
\newif\iffirstitemi@
\newif\iffirstitemii@
\newif\iffirstitemiii@
\newif\iffirstitemiv@
\newif\iffirstitemv@
\def\Firstitem@true{\csname firstitem\romannumeral\listlevel@
 @true\endcsname}
\def\Firstitem@false{\csname firstitem\romannumeral\listlevel@
 @false\endcsname}
\def\Listm@{\csname listm\romannumeral\listlevel@ @\endcsname}
\def\Item@{\csname item\romannumeral\listlevel@ @\endcsname}
\def\Liste@{\csname liste\romannumeral\listlevel@ @\endcsname}
\newif\iflistcontinue@
\def\keepitem{\listcontinue@true}
\newcount\list@C@
\def\list{%
 \iflistcontinue@\csname list@C1\endcsname\csname list@C@\endcsname\fi
 \global\csname list@C2\endcsname\z@
 \global\csname list@C3\endcsname\z@
 \global\csname list@C4\endcsname\z@
 \global\csname list@C5\endcsname\z@
 \begingroup
 \firstitemi@true
 \listlevel@\@ne
 \def\item{\FN@\item@}%
 \FN@\list@}
\Invalid@\runinitem
\def\list@{\ifx\next\par
 \DN@\par{\FN@\list@}\else
 \ifx\next\runinitem
  \DN@\runinitem{\FN@\runinitem@}\else
  \DN@{\par\dimen@\parskip\parskip\dimen@}\fi\fi\next@}
\newif\ifoutlevel@
\newif\ifrunin@
\def\item@{%
 \ifoutlevel@\Liste@\outlevel@false\fi
 \ifrunin@\runin@false\par
  \dimen@\parskip\parskip\dimen@
  \Listm@\fi
 \iffirstitemi@\listbi@\listmi@\firstitemi@false\else\par\fi
 \iffirstitemii@\listbii@\listmii@\firstitemii@false\else\par\fi
 \iffirstitemiii@\listbiii@\listmiii@\firstitemiii@false\else\par\fi
 \iffirstitemiv@\listbiv@\listmiv@\firstitemiv@false\else\par\fi
 \iffirstitemv@\listbv@\listmv@\firstitemv@false\else\par\fi
 \DN@"##1"{{\let\pre\list@@P\let\post\list@@Q
  \let\style\list@@S\let\numstyle\list@@N
  \vskip-\parskip
  \Item@{\list@@F##1}%
  \noexpands@
  \Qlabel@{##1}}%
  \locallabel@
  \FNSSP@}%
 \DNii@{\global\advance\list@@C\@ne
  {\noexpands@
   \xdef\Thelabel@@@{\number\list@@C}%
   \xdefThelabel@\list@@N
   \xdef\Thelabel@@@@{\list@@P\Thelabel@\list@@Q}%
   \xdefThelabel@@\list@@S
  }%
  \locallabel@
  \vskip-\parskip
  \Item@{\list@@F\thelabel@@}%
  \FN@\pretendspace@}%
 \ifx\next"\expandafter\next@\else\expandafter\nextii@\fi}
\def\runinitem@{%
  \runin@true
  \Firstitem@false
  \DN@"##1"{{\let\pre\list@@P\let\post\list@@Q
   \let\style\list@@S\let\numstyle\list@@N
   \unskip\space{\list@@F##1} %
   \noexpands@
   \Qlabel@{##1}}%
   \locallabel@
   \ignorespaces}%
  \DNii@{\global\advance\list@@C\@ne
   {\noexpands@
    \xdef\Thelabel@@@{\number\list@@C}%
    \xdefThelabel@\list@@N
    \xdef\Thelabel@@@@{\list@@P\Thelabel@\list@@Q}%
    \xdefThelabel@@\list@@S
   }%
   \locallabel@
   \unskip\space{\list@@F\thelabel@@} }%
  \ifx\next"\expandafter\next@\else\expandafter\nextii@\fi}
\def\inlevel{\ifnum\listlevel@=5
 \DN@{\Err@{Already 5 levels down}}\else
 \DN@{\begingroup\advance\listlevel@\@ne
 \Firstitem@true\FN@\inlevel@}\fi\next@}
\def\inlevel@{\ifx\next\par
 \DN@\par{\FN@\inlevel@}\else
 \ifx\next\runinitem
  \DN@\runinitem{\FN@\runinitem@}\else
  \let\next@\relax\fi\fi\next@}
\def\outlevel{\ifnum\listlevel@=\@ne
 \Err@{At top level}\else
 \par\global\list@@C\z@\endgroup\outlevel@true\fi}
\def\endlist{%
 \expandafter\global\csname list@C@\endcsname\csname list@C1\endcsname
 \par
 \global\toks\@ne{}\count@\listlevel@
 {\loop
  \ifnum\count@>\z@\global\toks\@ne\expandafter{\the\toks\@ne\endgroup}%
  \advance\count@\m@ne
  \repeat}%
 \the\toks\@ne
 \liste@
 \listcontinue@false\global\csname list@C1\endcsname\z@
 \vskip-\parskip
 \noindent@@
 \FN@\pretendspace@}
\newif\iffirstdescribe@
\def\describe{\par
 \begingroup\firstdescribe@true
 \def\item##1{%
  \iffirstdescribe@\penalty50 \medskip\vskip-\parskip
  \firstdescribe@false\else\par\fi
  \noindent@@\hangindent2pc\hangafter\@ne
  {\bf##1}\hskip.5em}}

\Invalid@\pullin
\Invalid@\pullinmore
\newif\iffirstpull@
\def\margins{\par\begingroup\firstpull@true
 \def\pullin##1##2{\par
  \iffirstpull@\firstpull@false\else\endgroup\fi
  \begingroup\DN@{##1}%
  \ifx\next@\empty\leftskip\z@\else\ifx\next@\space\leftskip\z@
  \else\leftskip##1\fi\fi
  \DN@{##2}\ifx\next@\empty\rightskip\z@\else\ifx\next@\space
  \rightskip\z@\else\rightskip##2\fi\fi\ignorespaces}%
 \def\pullinmore##1##2{\par
  \xdef\Next@{\leftskip\the\leftskip\relax\rightskip\the\rightskip\relax}%
  \iffirstpull@\firstpull@false\else\endgroup\fi
  \begingroup\Next@
  \DN@{##1}%
  \ifx\next@\empty\else\ifx\next@\space\else\advance\leftskip##1\fi\fi
  \DN@{##2}\ifx\next@\empty\else\ifx\next@\space\else
  \advance\rightskip##2\fi\fi\ignorespaces}}

\newif\ifnopunct@
\newif\ifnospace@
\newif\ifoverlong@
\let\nofrillslist@\empty
\let\overlonglist@\empty
\def\nopunct{\nopunct@true\FN@\nopunct@}
\def\nospace{\nospace@true\FN@\nospace@}
\def\overlong{\overlong@true\FN@\overlong@}
\def\nopunct@{\ifx\next\nospace
 \DN@\nospace{\nospace@true\FN@\nopnos@}\else\ifx\next\overlong
 \DN@\overlong{\overlong@true\FN@\nopol@}\else
 \let\next@\nopunct@@\fi\fi\next@}
\def\nopunct@@#1{\ismember@\nofrillslist@#1%
 \iftest@\let\next@#1\else
 \DN@{\nopunct@false\Err@{\noexpand\nopunct can't be used with
 \string#1}#1}\fi\next@}
\def\nospace@{\ifx\next\nopunct
 \DN@\nopunct{\nopunct@true\FN@\nopnos@}\else\ifx\next\overlong
 \DN@\overlong{\overlong@true\FN@\nosol@}\else
 \let\next@\nospace@@\fi\fi\next@}
\def\nospace@@#1{\ismember@\nofrillslist@#1%
 \iftest@\let\next@#1\else
 \DN@{\nospace@false\Err@{\noexpand\nospace can't be used with
 \string#1}#1}\fi\next@}
\def\overlong@{\ifx\next\nopunct
 \DN@\nopunct{\nopunct@true\FN@\nopol@}\else\ifx\next\nospace
 \DN@\nospace{\nospace@true\FN@\nosol@}\else
 \let\next@\overlong@@\fi\fi\next@}
\def\overlong@@#1{\ismember@\overlonglist@#1%
 \iftest@\let\next@#1\else
 \DN@{\overlong@false\Err@{\noexpand\overlong can't be used with
 \string#1}#1}\fi\next@}
\def\nopnos@{\ifx\next\overlong
 \DN@\overlong{\overlong@true\nopnosol@}\else
 \let\next@\nopnos@@\fi\next@}
\def\nopol@{\ifx\next\nospace
 \DN@\nospace{\nospace@true\nopnosol@}\else
 \let\next@\nopol@@\fi\next@}
\def\nosol@{\ifx\next\nopunct
 \DN@\nopunct{\nopunct@true\nopnosol@}\else
 \let\next@\nosol@@\fi\next@}
\def\nopnos@@#1{\ismember@\nofrillslist@#1%
 \iftest@\let\next@#1\else
 \DN@{\nopunct@false\nospace@false
  \Err@{\noexpand\nopunct\noexpand\nospace
   can't be used with \string#1}#1}\fi\next@}
\def\testii@#1{\ismember@\nofrillslist@#1%
 \iftest@\let\nextiii@ T\else\let\nextiii@ F\fi
 \ismember@\overlonglist@#1%
 \iftest@\let\nextiv@ T\else\let\nextiv@ F\fi
 \test@false\if\nextiii@ T\if\nextiv@ T\test@true\fi\fi}
\def\nopol@@#1{\testii@{#1}%
 \iftest@\let\next@#1%
 \else\DN@{\if\nextiii@ T\else\nopunct@false\fi
  \if\nextiv@ T\else\overlong@false\fi
  \Err@{\if\nextiii@ T\else\noexpand\nopunct\fi
  \if\nextiv@ T\else\noexpand\overlong\fi can't be used
  with \string#1}#1}\fi\next@}
\def\nosol@@#1{\testii@{#1}%
 \iftest@\let\next@#1%
 \else\DN@{\if\nextiii@ T\else\nospace@false\fi
  \if\nextiv@ T\else\overlong@false\fi
  \Err@{\if\nextiii@ T\else\noexpand\nospace\fi
  \if\nextiv@ T\else\noexpand\overlong\fi can't be used
  with \string#1}#1}\fi\next@}
\def\nopnosol@#1{\testii@{#1}%
 \iftest@\let\next@#1%
 \else\DN@{\if\nextiii@ T\else\nopunct@false\nospace@false\fi
  \if\nextiv@ T\else\overlong@false\fi
  \Err@{\if\nextiii@ T\else\noexpand\nopunct\noexpand\nospace\fi
  \if\nextiv@ T\else\noexpand\overlong\fi can't be used
  with \string#1}#1}\fi\next@}
\def\punct@#1{\ifnopunct@\else#1\fi}
\def\addspace@#1{\ifnospace@\else#1\fi}
\def\hss@{\ifoverlong@\z@ plus\@m\p@ minus\@m\p@
 \else \z@ plus\@m\p@\fi}
\rightadd@\demo\to\nofrillslist@
\newif\ifclaim@
\def\exxx@{\expandafter\expandafter\expandafter\eat@\expandafter\string}
\let\colon@:
\def\demo#1{\ifclaim@
 \Err@{Previous \expandafter\noexpand\claimtype@ has
  no matching \string\end\exxx@\claimtype@}%
 \let\next@\relax
 \else
  \par
  \ifdim\lastskip<\smallskipamount\removelastskip\smallskip\fi
  \begingroup
  \noindent@@{\smc\ignorespaces#1\unskip
   \punct@{\null\colon@}\addspace@\enspace}%
  \nopunct@false\nospace@false
  \rm
  \DN@{\FNSSP@}%
 \fi
 \next@}
\def\enddemo{\par\endgroup\nopunct@false\nospace@false\smallskip}
\rightadd@\claim\to\nofrillslist@
\def\claim@F{\smc}
\def\claim@@@F{\csname\exxx@\claimtype@ @F\endcsname}
\def\claimformat@#1#2#3{%
 \medbreak\noindent@@{\smc#1 {\claim@@@F#2} #3%
 \punct@{\null.}\addspace@\enspace}\sl}
\def\claimformat@@#1#2{\claimformat@{\ignorespaces#1\unskip}%
 {\ifx\thelabel@@\empty\unskip\else\thelabel@@\fi}%
 {\ignorespaces#2\unskip}%
 \let\Claimformat@@\claimformat@@\FNSSP@}
\let\Claimformat@@\claimformat@@
\def\claim@@@P{\csname\exxx@\claimtype@ @P\endcsname}
\def\claim@@@Q{\csname\exxx@\claimtype@ @Q\endcsname}
\def\claim@@@S{\csname\exxx@\claimtype@ @S\endcsname}
\def\claim@@@N{\csname\exxx@\claimtype@ @N\endcsname}
\def\claim@@@C{\csname claim@C\claimclass@\endcsname}
\newcount\claim@C
\claim@C\z@
\let\claim@P\empty
\let\claim@Q\empty
\def\claim@S#1{#1\/}
\let\claim@N\arabic
\def\claim{\claim@true\let\claimclass@\empty
 \def\claimtype@{\claim}\FN@\claim@}
\def\claim@{%
 \ifx\next\c
  \let\next@\claim@c
 \else
  \ifx\next"%
   \let\next@\claim@q
  \else
   \begingroup\global\advance\claim@C\@ne
   {\noexpands@
    \xdef\Thelabel@@@{\number\claim@C}%
    \xdefThelabel@\claim@N
    \xdef\Thelabel@@@@{\claim@P\Thelabel@\claim@Q}%
    \xdefThelabel@@\claim@S
   }%
   \locallabel@
   \let\next@\Claimformat@@
  \fi
 \fi
 \next@}
\def\claim@c\c#1{\claim@true\begingroup
 \expandafter
 \ifx\csname claim@C#1\endcsname\relax
  \expandafter\newcount@\csname claim@C#1\endcsname
  \global\csname claim@C#1\endcsname\@ne
 \else
  \global\advance\csname claim@C#1\endcsname\@ne
 \fi
 \def\claimclass@{#1}%
 {\noexpands@
  \xdef\Thelabel@@@{\number\claim@@@C}%
  \xdefThelabel@\claim@@@N
  \xdef\Thelabel@@@@{\claim@@@P\Thelabel@\claim@@@Q}%
  \xdefThelabel@@\claim@@@S
 }%
 \locallabel@
 \FNSS@\claim@c@}
\def\claim@q"#1"{\begingroup
 {\let\pre\claim@@@P\let\post\claim@@@Q
  \let\style\claim@@@S\let\numstyle\claim@@@N
  \noexpands@
  \Qlabel@{#1}}%
 \locallabel@
 \FNSS@\claim@q@}
\def\claim@c@{\ifx\next"%
 \global\advance\claim@@@C\m@ne\let\next@\claim@cq
 \else\let\next@\Claimformat@@\fi\next@}
\def\claim@cq"#1"{{\let\pre\claim@@@P\let\post\claim@@@Q
 \let\style\claim@@@S\let\numstyle\claim@@@N
 \noexpands@
 \Qlabel@{#1}}%
 \locallabel@
 \FNSS@\Claimformat@@}
\def\claim@q@{\ifx\next\c\expandafter\claim@qc
 \else\expandafter\Claimformat@@\fi}
\def\claim@qc\c#1{\expandafter\ifx\csname claim@C#1\endcsname\relax
 \expandafter\newcount@\csname claim@C#1\endcsname
 \global\csname claim@C#1\endcsname\z@\fi
 \FNSS@\Claimformat@@}
\def\endclaim{\endgroup\claim@false\nopunct@false\nospace@false
 \let\Claimformat@@\claimformat@@\medbreak}
\Invalid@\claimclause
\def\newclaim{\FN@\newclaim@}
\def\newclaim@{\ifx\next\claimclause
 \DN@\claimclause##1{\newclaim@@{##1}}\else
 \DN@{\newclaim@@\relax}\fi\next@}
\def\claimlist@{\\\claim}
\newtoks\claim@i
\newtoks\claim@v
\let\noclaimclause@=F
\def\newclaim@@#1#2#3\c#4#5{\define#2{}%
 \rightadd@#2\to\claimlist@\rightadd@#2\to\nofrillslist@%
 \expandafter\def\csname\exstring@#2@P\endcsname{\claim@P}%
 \expandafter\def\csname\exstring@#2@Q\endcsname{\claim@Q}%
 \expandafter\def\csname\exstring@#2@S\endcsname{\claim@S}%
 \expandafter\def\csname\exstring@#2@N\endcsname{\claim@N}%
 \expandafter\def\csname\exstring@#2@F\endcsname{\claim@F}%
 \expandafter\def\csname end\exstring@#2\endcsname{\endclaim}%
 \expandafter\ifx\csname claim@C#4\endcsname\relax
  \expandafter\newcount@\csname claim@C#4\endcsname
  \global\csname claim@C#4\endcsname\z@\fi
 \edef\next@{\let\csname\exstring@#2@C\endcsname
   \csname claim@C#4\endcsname}\next@
 \def#2{\ifx\noclaimclause@ T\else#1\fi
  \global\claim@i{#1}\gdef\claim@iv{#4}\global\claim@v{#5}%
  \def\claimtype@{#2}\def\Claimformat@@{\claimformat@@{#5}}\claim@c\c{#4}}}
\def\shortenclaim#1#2{\define#2{}%
 \ismember@\claimlist@#1%
 \iftest@
  \rightadd@#2\to\nofrillslist@%
  \expandafter\def\csname\exstring@#2@P\endcsname
   {\csname\exstring@#1@P\endcsname}%
  \expandafter\def\csname\exstring@#2@Q\endcsname
   {\csname\exstring@#1@Q\endcsname}%
  \expandafter\def\csname\exstring@#2@S\endcsname
   {\csname\exstring@#1@S\endcsname}%
  \expandafter\def\csname\exstring@#2@N\endcsname
   {\csname\exstring@#1@N\endcsname}%
  \expandafter\def\csname\exstring@#2@F\endcsname
   {\csname\exstring@#1@F\endcsname}%
  \expandafter\def\csname end\exstring@#2\endcsname{\endclaim}%
  \edef\next@{\let\csname\exstring@#2@C\endcsname
    \csname claim\exstring@#1C\endcsname}\next@
  \setbox\z@\vbox{\let\noclaimclause@ T#1""\relax\endgroup}%
  \edef#2{\the\claim@i
   \def\noexpand\claimtype@{\noexpand#2}%
   \def\noexpand\Claimformat@@{\noexpand\claimformat@@{\the\claim@v}\relax}%
   \noexpand\claim@c\noexpand\c{\claim@iv}}%
 \else
  \Err@{\noexpand#1not yet created by \string\newclaim}%
 \fi}
\def\classtest@#1{\DN@{#1}\ifx\next@\claimclass@
 \test@true\else\test@false\fi}
\def\typetest@#1{\DN@{#1}\ifx\next@\claimtype@\test@true\else
  \test@false\fi}
\newif\iftoc@
\def\tocfile{\iftoc@\else\alloc@@7\write\chardef\sixt@@n\toc@
 \immediate\openout\toc@=\jobname.toc
 \alloc@@7\write\chardef\sixt@@n\tic@
 \immediate\openout\tic@=\jobname.tic
 \global\toc@true\fi}
\rightadd@\hl\to\nofrillslist@
\rightadd@\HL\to\overlonglist@
\def\HL@@C{\csname HL@C\HLlevel@\endcsname}
\def\HL@@P{\csname HL@P\HLlevel@\endcsname}
\def\HL@@Q{\csname HL@Q\HLlevel@\endcsname}
\def\HL@@S{\csname HL@S\HLlevel@\endcsname}
\def\HL@@N{\csname HL@N\HLlevel@\endcsname}
\def\HL@@F{\csname HL@F\HLlevel@\endcsname}
\def\HL@@@C{\csname\exxx@\HLtype@ @C\endcsname}
\def\HL@@@P{\csname\exxx@\HLtype@ @P\endcsname}
\def\HL@@@Q{\csname\exxx@\HLtype@ @Q\endcsname}
\def\HL@@@S{\csname\exxx@\HLtype@ @S\endcsname}
\def\HL@@@N{\csname\exxx@\HLtype@ @N\endcsname}
\def\HL#1{\expandafter
 \ifx\csname HL@C#1\endcsname\relax
  \DN@{\Err@{\string\HL#1 not defined in this style}}%
 \else
  \DN@{\gdef\HLlevel@{#1}\def\HLname@{\HL{#1}}\let\HLtype@\relax\FNSS@\HL@}%
 \fi
 \next@}%
\newif\ifquoted@
\let\aftertoc@\relax
\def\HL@{%
 \DN@"##1"##2\endHL{\def\entry@{##2}\quoted@true
  {\noexpands@
  \ifx\HLtype@\relax
   \let\pre\HL@@P\let\post\HL@@Q\let\style\HL@@S\let\numstyle\HL@@N
  \else
   \let\pre\HL@@@P\let\post\HL@@@Q\let\style\HL@@@S\let\numstyle\HL@@@N
  \fi
  \Qlabel@{##1}\let\style\relax\xdef\Qlabel@@@@{##1}%
  \xdef\Thepref@{\Thelabel@@@@}}%
  \csname HL@\HLlevel@\endcsname##2\endHL
  \let\pref\Thepref@
  \csname HL@I\HLlevel@\endcsname
  \csname HL@J\HLlevel@\endcsname
  \let\pref\pref@
  \HLtoc@	
  \aftertoc@
  \let\aftertoc@\relax\overlong@false}%
 \DNii@##1\endHL{\def\entry@{##1}\quoted@false
  {\noexpands@
  \ifx\HLtype@\relax
   \global\advance\HL@@C\@ne
   \xdef\Thelabel@@@{\number\HL@@C}%
   \xdefThelabel@{\HL@@N}%
   \xdef\Thelabel@@@@{\HL@@P\Thelabel@\HL@@Q}%
   \xdefThelabel@@{\HL@@S}%
  \else
   \global\advance\HL@@@C\@ne
   \xdef\Thelabel@@@{\number\HL@@@C}%
   \xdefThelabel@{\HL@@@N}%
   \xdef\Thelabel@@@@{\HL@@@P\Thelabel@\HL@@@Q}%
   \xdefThelabel@@{\HL@@@S}%
  \fi
  \xdef\Thepref@{\Thelabel@@@@}}%
  \csname HL@\HLlevel@\endcsname##1\endHL
  \let\pref\Thepref@
  \csname HL@I\HLlevel@\endcsname
  \csname HL@J\HLlevel@\endcsname
  \let\pref\pref@
  \HLtoc@
  \aftertoc@
  \let\aftertoc@\relax\overlong@false}%
 \ifx\next"\expandafter\next@\else\expandafter\nextii@\fi}%
\Invalid@\endHL
\def\hl@@C{\csname hl@C\hllevel@\endcsname}
\def\hl@@P{\csname hl@P\hllevel@\endcsname}
\def\hl@@Q{\csname hl@Q\hllevel@\endcsname}
\def\hl@@S{\csname hl@S\hllevel@\endcsname}
\def\hl@@N{\csname hl@N\hllevel@\endcsname}
\def\hl@@F{\csname hl@F\hllevel@\endcsname}
\def\hl@@@C{\csname\exxx@\hltype@ @C\endcsname}
\def\hl@@@P{\csname\exxx@\hltype@ @P\endcsname}
\def\hl@@@Q{\csname\exxx@\hltype@ @Q\endcsname}
\def\hl@@@S{\csname\exxx@\hltype@ @S\endcsname}
\def\hl@@@N{\csname\exxx@\hltype@ @N\endcsname}
\def\hl#1{\expandafter
 \ifx\csname hl@C#1\endcsname\relax
  \DN@{\Err@{\string\hl#1 not defined in this style}}%
 \else
  \DN@{\gdef\hllevel@{#1}\def\hlname@{\hl{#1}}\let\hltype@\relax\FNSS@\hl@}%
 \fi
 \next@}
\def\hl@{%
 \DN@"##1"##2{\def\entry@{##2}\quoted@true
  {\noexpands@
  \ifx\hltype@\relax
   \let\pre\hl@@P\let\post\hl@@Q\let\style\hl@@S\let\numstyle\hl@@N
  \else
   \let\pre\hl@@@P\let\post\hl@@@Q\let\style\hl@@@S\let\numstyle\hl@@@N
  \fi
  \Qlabel@{##1}\let\style\relax\xdef\Qlabel@@@@{##1}%
  \xdef\Thepref@{\Thelabel@@@@}}%
  \csname hl@\hllevel@\endcsname{##2}%
  \let\pref\Thepref@
  \csname hl@I\hllevel@\endcsname
  \csname hl@J\hllevel@\endcsname
  \let\pref\pref@
  \hltoc@
  \aftertoc@
  \let\aftertoc@\relax\nopunct@false\nospace@false\FNSSP@}%
 \DNii@##1{\def\entry@{##1}\quoted@false
  {\noexpands@
  \ifx\hltype@\relax
   \global\advance\hl@@C\@ne
   \xdef\Thelabel@@@{\number\hl@@C}%
   \xdefThelabel@{\hl@@N}%
   \xdef\Thelabel@@@@{\hl@@P\Thelabel@\hl@@Q}%
   \xdefThelabel@@{\hl@@S}%
  \else
   \global\advance\hl@@@C\@ne
   \xdef\Thelabel@@@{\number\hl@@@C}%
   \xdefThelabel@{\hl@@@N}%
   \xdef\Thelabel@@@@{\hl@@@P\Thelabel@\hl@@@Q}%
   \xdefThelabel@@{\hl@@@S}%
  \fi
  \xdef\Thepref@{\Thelabel@@@@}}%
  \csname hl@\hllevel@\endcsname{##1}%
  \let\pref\Thepref@
  \csname hl@I\hllevel@\endcsname
  \csname hl@J\hllevel@\endcsname
  \let\pref\pref@
  \hltoc@
  \aftertoc@
  \let\aftertoc@\relax\nopunct@false\nospace@false\FNSSP@}%
 \ifx\next"\expandafter\next@\else\expandafter\nextii@\fi}%
\def\six@#1#2 #3 #4 #5 #6 #7 {\DN@{#2}\ifx\next@\empty
 \DN@##1\six@{}\else
 \write#1{ #2 #3 #4 #5 #6 #7}\DN@{\six@#1}\fi
 \next@}
\def\Sixtoc@{\ifx\macdef@\empty\else
 \DN@##1##2\next@{\def\macdef@{##1##2}}%
 \expandafter\next@\macdef@\next@
 \edef\next@
  {\noexpand\six@\toc@\macdef@
  \space\space\space\space\space\space\space\space\space\space\space\space
  \noexpand\six@}%
 \next@\let\macdef@\relax\fi}
\def\QorThelabel@@@@{\ifquoted@
 \noexpand\noexpand\noexpand"\Qlabel@@@@\noexpand\noexpand\noexpand"\else
 \Thelabel@@@@\fi}
\def\HLtoc@{%
 \iftoc@
 \expandafter\expandafter\expandafter\unmacro@
  \expandafter\meaning\csname HL@W\HLlevel@\endcsname\unmacro@
  {\noexpands@\let\style\relax
   \edef\next@{\write\toc@{\noexpand\noexpand\expandafter\noexpand\HLname@
   {\macdef@}{\QorThelabel@@@@}}}%
  \next@}%
  \expandafter\unmacro@\meaning\entry@\unmacro@
  \Sixtoc@
  \write\toc@{\noexpand\Page{\number\pageno}{\page@N}%
   {\page@P}{\page@Q}^^J}%
 \fi}
\def\hltoc@{%
 \iftoc@
 \expandafter\expandafter\expandafter\unmacro@
  \expandafter\meaning\csname hl@W\hllevel@\endcsname\unmacro@
  {\noexpands@\let\style\relax
  \edef\next@{\write\toc@{%
   \ifnopunct@\noexpand\noexpand\noexpand\nopunct\fi
   \ifnospace@\noexpand\noexpand\noexpand\nospace\fi
   \noexpand\noexpand\expandafter\noexpand\hlname@
   {\macdef@}{\QorThelabel@@@@}}}%
  \next@}%
  \expandafter\unmacro@\meaning\entry@\unmacro@
  \Sixtoc@
  \write\toc@{\noexpand\Page{\number\pageno}{\page@N}%
   {\page@P}{\page@Q}^^J}%
 \fi}
\def\mainfile#1{\def\mainfile@{#1}}
\def\checkmainfile@{\ifx\mainfile@\undefined
 \Err@{No \noexpand\mainfile specified}\fi}
\expandafter\newcount@\csname HL@C1\endcsname
\csname HL@C1\endcsname\z@
\expandafter\def\csname HL@S1\endcsname#1{#1\null.}
\expandafter\let\csname HL@N1\endcsname\arabic
\expandafter\let\csname HL@P1\endcsname\empty
\expandafter\let\csname HL@Q1\endcsname\empty
\expandafter\def\csname HL@F1\endcsname{\bf}
\expandafter\let\csname HL@W1\endcsname\empty
\expandafter\newcount@\csname hl@C1\endcsname
\csname hl@C1\endcsname\z@
\expandafter\def\csname hl@S1\endcsname#1{#1\/}
\expandafter\let\csname hl@N1\endcsname\arabic
\expandafter\let\csname hl@P1\endcsname\empty
\expandafter\let\csname hl@Q1\endcsname\empty
\expandafter\def\csname hl@F1\endcsname{\bf}
\expandafter\let\csname hl@W1\endcsname\empty
\expandafter\def\csname HL@1\endcsname#1\endHL{\bigbreak
 {\locallabel@
  \global\setbox\@ne\vbox{\Let@\tabskip\hss@
  \halign to\hsize{\bf\hfil\ignorespaces##\unskip\hfil\cr
  \expandafter\ifx\csname HL@W1\endcsname\empty\else
   \csname HL@W1\endcsname\space\fi
  {\HL@@F\ifx\thelabel@@\empty\else\thelabel@@\space\fi}%
  \ignorespaces#1\crcr}}%
  }%
 \unvbox\@ne\nobreak\medskip}
\expandafter\def\csname hl@1\endcsname#1{\medbreak\noindent@@
 {\locallabel@
 \bf{\hl@@F\ifx\thelabel@@\empty\else\thelabel@@\space\fi}%
 \ignorespaces#1\unskip\punct@{\null.}\addspace@\enspace}}
\expandafter\def\csname HL@I1\endcsname{\Reset\hl1{1}%
 \ifx\pref\empty\newpre\hl1{}\else\newpre\hl1{\pref.}\fi}
\def\NameHL#1#2{\define#2{}%
 \expandafter\ifx\csname HL@R#1\endcsname\relax
 \else
  \def\nextiv@{\let\nextiii@}%
  \expandafter\nextiv@\csname HL@R#1\endcsname
  \expandafter\let\nextiii@\undefined
  \expandafter\let\csname\exxx@\nextiii@ @C\endcsname\relax
  \expandafter\let\csname\exxx@\nextiii@ @P\endcsname\relax
  \expandafter\let\csname\exxx@\nextiii@ @Q\endcsname\relax
  \expandafter\let\csname\exxx@\nextiii@ @S\endcsname\relax
  \expandafter\let\csname\exxx@\nextiii@ @N\endcsname\relax
  \expandafter\let\csname\exxx@\nextiii@ @F\endcsname\relax
  \expandafter\let\csname\exxx@\nextiii@ @W\endcsname\relax
  \expandafter\let\csname end\exxx@\nextiii@\endcsname\undefined
 \fi
 \expandafter\gdef\csname HL@R#1\endcsname{#2}%
 \expandafter\gdef\csname\exstring@#2@R\endcsname{{HL}{#1}}%
 \iftoc@\write\toc@{\noexpand\NameHL#1\noexpand#2^^J}\fi
 \rightadd@#2\to\overlonglist@
 \edef\next@{\let\csname\exstring@#2@C\endcsname\expandafter\noexpand
  \csname HL@C#1\endcsname}\next@
 \edef\next@{\let\csname\exstring@#2@P\endcsname\expandafter\noexpand
  \csname HL@P#1\endcsname}\next@
 \edef\next@{\let\csname\exstring@#2@Q\endcsname\expandafter\noexpand
  \csname HL@Q#1\endcsname}\next@
 \edef\next@{\let\csname\exstring@#2@S\endcsname\expandafter\noexpand
  \csname HL@S#1\endcsname}\next@
 \edef\next@{\let\csname\exstring@#2@N\endcsname\expandafter\noexpand
  \csname HL@N#1\endcsname}\next@
 \edef\next@{\let\csname\exstring@#2@F\endcsname\expandafter\noexpand
  \csname HL@F#1\endcsname}\next@
 \edef\next@{\let\csname\exstring@#2@W\endcsname\expandafter\noexpand
  \csname HL@W#1\endcsname}\next@
 \edef\next@{\def\noexpand#2####1\expandafter\noexpand
  \csname end\exstring@#2\endcsname
  {\def\noexpand\HLtype@{\noexpand#2}%
   \def\noexpand\HLname@{\noexpand#2}%
   \gdef\noexpand\HLlevel@{#1}%
   \noexpand\FNSS@\noexpand\HL@####1\noexpand\endHL}}%
  \next@
 \edef\next@{\noexpand\Invalid@\expandafter\noexpand
  \csname end\exstring@#2\endcsname}%
 \next@}
\def\Namehl#1#2{\define#2{}%
 \expandafter\ifx\csname hl@R#1\endcsname\relax
 \else
  \def\nextiv@{\let\nextiii@}%
  \expandafter\nextiv@\csname hl@R#1\endcsname
  \expandafter\let\nextiii@\undefined
  \expandafter\let\csname\exxx@\nextiii@ @C\endcsname\relax
  \expandafter\let\csname\exxx@\nextiii@ @P\endcsname\relax
  \expandafter\let\csname\exxx@\nextiii@ @Q\endcsname\relax
  \expandafter\let\csname\exxx@\nextiii@ @S\endcsname\relax
  \expandafter\let\csname\exxx@\nextiii@ @N\endcsname\relax
  \expandafter\let\csname\exxx@\nextiii@ @F\endcsname\relax
  \expandafter\let\csname\exxx@\nextiii@ @W\endcsname\relax
 \fi
 \expandafter\gdef\csname hl@R#1\endcsname{#2}%
 \expandafter\gdef\csname\exstring@#2@R\endcsname{{hl}{#1}}%
 \iftoc@\write\toc@{\noexpand\Namehl#1\noexpand#2^^J}\fi
 \rightadd@#2\to\nofrillslist@%
 \edef\next@{\let\csname\exstring@#2@C\endcsname\expandafter\noexpand
  \csname hl@C#1\endcsname}\next@
 \edef\next@{\let\csname\exstring@#2@P\endcsname\expandafter\noexpand
  \csname hl@P#1\endcsname}\next@
 \edef\next@{\let\csname\exstring@#2@Q\endcsname\expandafter\noexpand
  \csname hl@Q#1\endcsname}\next@
 \edef\next@{\let\csname\exstring@#2@S\endcsname\expandafter\noexpand
  \csname hl@S#1\endcsname}\next@
 \edef\next@{\let\csname\exstring@#2@N\endcsname\expandafter\noexpand
  \csname hl@N#1\endcsname}\next@
 \edef\next@{\let\csname\exstring@#2@F\endcsname\expandafter\noexpand
  \csname hl@F#1\endcsname}\next@
 \edef\next@{\let\csname\exstring@#2@W\endcsname\expandafter\noexpand
  \csname hl@W#1\endcsname}\next@
 \edef\next@{\def\noexpand#2{%
  \def\noexpand\hltype@{\noexpand#2}%
  \def\noexpand\hlname@{\noexpand#2}%
  \gdef\noexpand\hllevel@{#1}%
  \noexpand\FNSS@\noexpand\hl@}}%
 \next@}%
\def\Initialize{\FN@\Init@}
\def\Init@{\ifx\next\HL\let\next@\InitH@\else\ifx\next\hl\let\next@\InitH@
  \else\let\next@\InitS@\fi\fi\next@}
\def\InitH@#1#2{\expandafter\ifx\csname\exstring@#1@C#2\endcsname\relax
 \DN@{\Err@{\noexpand#1level #2 not defined in this style}}\else
 \DN@{\expandafter\gdef\csname\exstring@#1@J#2\endcsname}\fi\next@}
\def\InitC@#1#2{\edef\nextii@{\expandafter\noexpand\csname#1\endcsname{#2}}}
\def\InitS@#1{\expandafter\ifx\csname\exstring@#1@R\endcsname\relax
 \Err@{\noexpand#1not defined in this style}\let\next@\relax\else
 \DN@{\let\next@}\expandafter\next@\csname\exstring@#1@R\endcsname
 \expandafter\InitC@\next@
 \DN@{\expandafter\InitH@\nextii@}\fi\next@}
\def\value#1{\expandafter
 \ifx\csname\exstring@#1@C\endcsname\relax
  \expandafter\ifx\csname\exstring@#1@C1\endcsname\relax
   \DN@{\Err@{\noexpand\value can't be used with \string#1}}%
  \else
   \DN@{\value@#1}%
  \fi
 \else
  \DN@{\number\csname\exstring@#1@C\endcsname\relax}%
 \fi
 \next@}
\def\value@#1#2{\expandafter
 \ifx\csname\exstring@#1@C#2\endcsname\relax
  \DN@{\Err@{\string\value\string#1 can't be followed by \string#2}}%
 \else
  \DN@{\number\csname\exstring@#1@C#2\endcsname\relax}%
 \fi
 \next@}
\newcount\Value
\def\Evaluate#1{\expandafter
 \ifx\csname\exstring@#1@C\endcsname\relax
  \expandafter\ifx\csname\exstring@#1@C1\endcsname\relax
   \DN@{\Err@{\noexpand\Evaluate can't be used with \string#1}}%
  \else
   \DN@{\Evaluate@#1}%
  \fi
 \else
  \DN@{\global\Value\csname\exstring@#1@C\endcsname}%
 \fi
 \next@}
\def\Evaluate@#1#2{\expandafter
 \ifx\csname\exstring@#1@C#2\endcsname\relax
  \DN@{\Err@{\string\Evaluate\string#1 can't be followed by \string#2}}%
 \else
  \DN@{\global\Value\csname\exstring@#1@C#2\endcsname}%
 \fi\next@}
\def\pre#1{\expandafter
 \ifx\csname\exstring@#1@P\endcsname\relax
  \expandafter\ifx\csname\exstring@#1@P1\endcsname\relax
   \DN@{\Err@{\noexpand\pre can't be used with \string#1}}%
  \else
   \DN@{\pre@#1}%
  \fi
 \else
  \DN@{{\csname\exstring@#1@P\endcsname}}%
 \fi
 \next@}
\def\pre@#1#2{\expandafter
 \ifx\csname\exstring@#1@P#2\endcsname\relax
  \DN@{\Err@{\string\pre\string#1 can't be followed by \string#2}}%
 \else
  \DN@{{\csname\exstring@#1@P#2\endcsname}}%
 \fi
 \next@}
\def\post#1{\expandafter
 \ifx\csname\exstring@#1@Q\endcsname\relax
  \expandafter\ifx\csname\exstring@#1@Q1\endcsname\relax
   \DN@{\Err@{\noexpand\post can't be used with \string#1}}%
  \else
   \DN@{\post@#1}%
  \fi
 \else
  \DN@{{\csname\exstring@#1@Q\endcsname}}%
 \fi
 \next@}
\def\post@#1#2{\expandafter
 \ifx\csname\exstring@#1@Q#2\endcsname\relax
  \DN@{\Err@{\string\post\string#1 can't be followed by \string#2}}%
 \else
  \DN@{{\csname\exstring@#1@Q#2\endcsname}}%
 \fi
 \next@}
\def\style#1{\expandafter
 \ifx\csname\exstring@#1@S\endcsname\relax
  \expandafter\ifx\csname\exstring@#1@S1\endcsname\relax
   \DN@{\Err@{\noexpand\style can't be used with \string#1}}%
  \else
   \DN@{\style@#1}%
  \fi
 \else
  \DN@{\csname\exstring@#1@S\endcsname}%
 \fi
 \next@}
\def\style@#1#2{\expandafter
 \ifx\csname\exstring@#1@S#2\endcsname\relax
  \DN@{\Err@{\string\style\string#1 can't be followed by \string#2}}%
 \else
  \DN@{\csname\exstring@#1@S#2\endcsname}%
 \fi
 \next@}
\def\fontstyle#1{\expandafter
 \ifx\csname\exstring@#1@F\endcsname\relax
  \expandafter\ifx\csname\exstring@#1@F1\endcsname\relax
   \DN@{\Err@{\noexpand\fontstyle can't be used with \string#1}}%
  \else
   \DN@{\fontstyle@#1}%
  \fi
 \else
  \DN@##1{{\csname\exstring@#1@F\endcsname##1}}%
 \fi
 \next@}
\def\fontstyle@#1#2{\expandafter
 \ifx\csname\exstring@#1@F#2\endcsname\relax
  \DN@{\Err@{\string\fontstyle\string#1 can't be followed by \string#2}}%
 \else
  \DN@##1{{\csname\exstring@#1@F#2\endcsname##1}}%
 \fi
 \next@}
\def\Reset#1{\expandafter
 \ifx\csname\exstring@#1@C\endcsname\relax
  \expandafter\ifx\csname\exstring@#1@C1\endcsname\relax
   \DN@{\Err@{\noexpand\Reset can't be used with \string#1}}%
  \else
   \DN@{\Reset@#1}%
  \fi
 \else
  \DN@##1{\count@##1\relax\ifx#1\page\else\advance\count@\m@ne\fi
   \global\csname\exstring@#1@C\endcsname\count@}%
 \fi
 \next@}
\def\Reset@#1#2{\expandafter
 \ifx\csname\exstring@#1@C#2\endcsname\relax
  \DN@{\Err@{\string\Reset\string#1 can't be followed by \string#2}}%
 \else
  \DN@##1{\count@##1\relax\advance\count@\m@ne
   \global\csname\exstring@#1@C#2\endcsname\count@}%
 \fi
 \next@}
\def\Offset#1{\expandafter
 \ifx\csname\exstring@#1@C\endcsname\relax
  \expandafter\ifx\csname\exstring@#1@C1\endcsname\relax
   \DN@{\Err@{\noexpand\Offset can't be used with \string#1}}%
  \else
   \DN@{\Offset@#1}%
  \fi
 \else
  \DN@##1{\count@##1\relax\advance\count@\m@ne\global\advance
   \csname\exstring@#1@C\endcsname\count@}%
 \fi
 \next@}
\def\Offset@#1#2{\expandafter
 \ifx\csname\exstring@#1@C#2\endcsname\relax
  \DN@{\Err@{\string\Offset\string#1 can't be followed by \string#2}}%
 \else
  \DN@##1{\count@##1\relax\advance\count@\m@ne
   \global\advance\csname\exstring@#1@C#2\endcsname\count@}%
 \fi
 \next@}
\def\getR@#1#2{\def\nextiv@{\let\nextiii@}\expandafter\nextiv@
 \csname\exstring@#1@R#2\endcsname}
\def\letR@#1#2#3{\expandafter\let\csname#1@#3#2\endcsname\Next@}
\def\letR@@#1#2{\expandafter\let\csname\exstring@#1@#2\endcsname\Next@}
\def\newpre#1{\expandafter
 \ifx\csname\exstring@#1@P\endcsname\relax
  \expandafter\ifx\csname\exstring@#1@P1\endcsname\relax
   \DN@{\Err@{\noexpand\newpre can't be used with \string#1}}%
  \else
   \DN@{\newpre@#1}%
  \fi
 \else
  \DN@{%
   \DNii@{%
    \endgroup
    \expandafter\let\csname\exstring@#1@P\endcsname\Next@
    \expandafter\ifx\csname\exstring@#1@R\endcsname\relax\else
    \getR@#1{}\expandafter\letR@\nextiii@ P\fi
    }%
   \begingroup\noexpands@\afterassignment\nextii@\xdef\Next@}%
 \fi
 \next@}
\def\newpre@#1#2{\expandafter
 \ifx\csname\exstring@#1@P#2\endcsname\relax
  \DN@{\Err@{\string\newpre\string#1 can't be followed by \string#2}}%
 \else
  \DN@{%
   \DNii@{%
    \endgroup
    \expandafter\let\csname\exstring@#1@P#2\endcsname\Next@
    \expandafter\ifx\csname\exstring@#1@R#2\endcsname\relax\else
    \getR@#1{#2}\expandafter\letR@@\nextiii@ P\fi
    }%
   \begingroup\noexpands@\afterassignment\nextii@\xdef\Next@}%
 \fi
 \next@}
\def\newpost#1{\expandafter
 \ifx\csname\exstring@#1@Q\endcsname\relax
  \expandafter\ifx\csname\exstring@#1@Q1\endcsname\relax
   \DN@{\Err@{\noexpand\newpost can't be used with \string#1}}%
  \else
   \DN@{\newpost@#1}%
  \fi
 \else
  \DN@{%
   \DNii@{%
    \endgroup
    \expandafter\let\csname\exstring@#1@Q\endcsname\Next@
    \expandafter\ifx\csname\exstring@#1@R\endcsname\relax\else
    \getR@#1{}\expandafter\letR@\nextiii@ Q\fi
    }%
   \begingroup\noexpands@\afterassignment\nextii@\xdef\Next@}%
 \fi
 \next@}
\def\newpost@#1#2{\expandafter
 \ifx\csname\exstring@#1@Q#2\endcsname\relax
  \DN@{\Err@{\string\newpost\string#1 can't be followed by \string#2}}%
 \else
  \DN@{%
   \DNii@{%
    \endgroup
    \expandafter\let\csname\exstring@#1@Q#2\endcsname\Next@
    \expandafter\ifx\csname\exstring@#1@R#2\endcsname\relax\else
    \getR@#1{#2}\expandafter\letR@@\nextiii@ Q\fi
    }%
   \begingroup\noexpands@\afterassignment\nextii@\xdef\Next@}%
 \fi
 \next@}
\def\newstyle#1{\expandafter
 \ifx\csname\exstring@#1@S\endcsname\relax
  \expandafter\ifx\csname\exstring@#1@S1\endcsname\relax
   \DN@{\Err@{\noexpand\newstyle can't be used
    with \string#1}}%
  \else
   \DN@{\newstyle@#1}%
  \fi
 \else
  \DN@{%
   \DNii@{%
    \expandafter\let\csname\exstring@#1@S\endcsname\Next@
    \expandafter\ifx\csname\exstring@#1@R\endcsname\relax\else
    \getR@#1{}\expandafter\letR@\nextiii@ S\fi
    }%
   \afterassignment\nextii@\gdef\Next@}%
 \fi
 \next@}
\def\newstyle@#1#2{\expandafter
 \ifx\csname\exstring@#1@S#2\endcsname\relax
  \DN@{\Err@{\string\newstyle\string#1 can't be followed by
   \string#2}}%
 \else
  \DN@{%
   \DNii@{%
    \expandafter\let\csname\exstring@#1@S#2\endcsname\Next@
    \expandafter\ifx\csname\exstring@#1@R#2\endcsname\relax\else
    \getR@#1{#2}\expandafter\letR@@\nextiii@ S\fi
    }%
   \afterassignment\nextii@\gdef\Next@}%
 \fi
 \next@}
\def\newnumstyle#1{\expandafter
 \ifx\csname\exstring@#1@N\endcsname\relax
  \expandafter\ifx\csname\exstring@#1@N1\endcsname\relax
   \DN@{\Err@{\noexpand\newnumstyle can't be used with
    \string#1}}%
  \else
   \DN@{\newnumstyle@#1}%
  \fi
 \else
  \DN@##1{%
   \gdef\Next@{##1}%
    \expandafter\let\csname\exstring@#1@N\endcsname\Next@
    \expandafter\ifx\csname\exstring@#1@R\endcsname\relax\else
    \getR@#1{}\expandafter\letR@\nextiii@ N\fi
    }%
 \fi
 \next@}
\def\newnumstyle@#1#2{\expandafter
 \ifx\csname\exstring@#1@N#2\endcsname\relax
  \DN@{\Err@{\string\newnumstyle\string#1 can't be followed by
   \string#2}}%
 \else
  \DN@##1{%
   \gdef\Next@{##1}%
    \expandafter\let\csname\exstring@#1@N#2\endcsname\Next@
    \expandafter\ifx\csname\exstring@#1@R#2\endcsname\relax\else
    \getR@#1{#2}\expandafter\letR@@\nextiii@ N\fi
    }%
  \fi
 \next@}
\def\newfontstyle#1{\expandafter
 \ifx\csname\exstring@#1@F\endcsname\relax
  \expandafter\ifx\csname\exstring@#1@F1\endcsname\relax
   \DN@{\Err@{\noexpand\newfontstyle can't be used with
    \string#1}}%
  \else
   \DN@{\newfontstyle@#1}%
  \fi
 \else
  \DN@##1{%
   \gdef\Next@{##1}%
    \expandafter\let\csname\exstring@#1@F\endcsname\Next@
    \expandafter\ifx\csname\exstring@#1@R\endcsname\relax\else
    \getR@#1{}\expandafter\letR@\nextiii@ F\fi
    }%
 \fi
 \next@}
\def\newfontstyle@#1#2{\expandafter
 \ifx\csname\exstring@#1@F#2\endcsname\relax
  \DN@{\Err@{\string\newfontstyle\string#1 can't be followed by
   \string#2}}%
 \else
  \DN@##1{%
   \gdef\Next@{##1}%
    \expandafter\let\csname\exstring@#1@F#2\endcsname\Next@
    \expandafter\ifx\csname\exstring@#1@R#2\endcsname\relax\else
    \getR@#1{#2}\expandafter\letR@@\nextiii@ F\fi
    }%
 \fi
 \next@}
\def\word#1{\expandafter
 \ifx\csname\exstring@#1@W\endcsname\relax
  \expandafter\ifx\csname\exstring@#1@W1\endcsname\relax
   \DN@{\Err@{\noexpand\word can't be used with \string#1}}%
  \else
   \DN@{\word@#1}%
  \fi
 \else
  \DN@{{\csname\exstring@#1@W\endcsname}}%
 \fi
 \next@}
\def\word@#1#2{\expandafter
 \ifx\csname\exstring@#1@W#2\endcsname\relax
  \DN@{\Err@{\string\word\noexpand#1can't be followed by \string#2}}%
 \else
  \DN@{{\csname\exstring@#1@W#2\endcsname}}%
 \fi
 \next@}
\def\newword#1{\expandafter
 \ifx\csname\exstring@#1@W\endcsname\relax
  \expandafter\ifx\csname\exstring@#1@W1\endcsname\relax
   \DN@{\Err@{\noexpand\newword can't be used  with \string#1}}%
  \else
   \DN@{\newword@#1}%
  \fi
 \else
  \DN@{%
   \DNii@{%
    \expandafter\let\csname\exstring@#1@W\endcsname\Next@
    \expandafter\ifx\csname\exstring@#1@R\endcsname\relax\else
     \getR@#1{}\expandafter\letR@\nextiii@ W\fi
    }%
   \afterassignment\nextii@\gdef\Next@}%
 \fi
 \next@}
\def\newword@#1#2{\expandafter
 \ifx\csname\exstring@#1@W#2\endcsname\relax
  \DN@{\Err@{\string\newword\noexpand#1can't be followed by \string#2}}%
 \else
  \DN@{%
   \DNii@{%
    \expandafter\let\csname\exstring@#1@W#2\endcsname\Next@
    \expandafter\ifx\csname\exstring@#1@R#2\endcsname\relax\else
     \getR@#1{#2}\expandafter\letR@@\nextiii@ W\fi
    }%
   \afterassignment\nextii@\gdef\Next@}%
 \fi
 \next@}
\newif\iffn@
\newcount\footmark@C
\footmark@C\z@
\def\footmark@S#1{$^{#1}$}
\let\footmark@N\arabic
\def\footmark@F{\rm}
\def\foottext@S#1{$^{#1}$}
\def\foottext@F{\rm}
\let\modifyfootnote@\relax
\def\modifyfootnote#1{\def\modifyfootnote@{#1}}
\def\vfootnote@#1{\insert\footins
 \bgroup
 \floatingpenalty\@MM\interlinepenalty\interfootnotelinepenalty
 \leftskip\z@\rightskip\z@\spaceskip\z@\xspaceskip\z@
 \rm\splittopskip\ht\strutbox\splitmaxdepth\dp\strutbox
 \locallabel@\noindent@@{\foottext@F#1}\modifyfootnote@
 \footstrut\FN@\fo@t}
\def\fo@t{\ifcat\bgroup\noexpand\next\expandafter\f@@t\else
 \expandafter\f@t\fi}
\def\f@t#1{#1\@foot}
\def\f@@t{\bgroup\aftergroup\@foot\afterassignment\FNSSP@\let\next@}
\def\@foot{\unskip\lower\dp\strutbox\vbox to\dp\strutbox{}\egroup
 \iffn@\expandafter\fn@false\else
 \expandafter\postvanish@\fi}
\newif\ifplainfn@
\plainfn@true
\def\fancyfootnotes{\plainfn@false}
\newcount\fancyfootmarkcount@
\fancyfootmarkcount@\z@
\newcount\lastfnpage@
\lastfnpage@-\@M
\let\justfootmarklist@\empty
\def\footmark{\let\@sf\empty
 \ifhmode\edef\@sf{\spacefactor\the\spacefactor}\/\fi
 \DN@{\ifx"\next\expandafter\nextii@\else\expandafter\footmark@\fi}%
 \DNii@"##1"{%
  \iffirstchoice@
   {\let\style\footmark@S\let\numstyle\footmark@N
   \footmark@F##1%
   \noexpands@
   \let\style\foottext@S
   \Qlabel@{##1}%
   }%
   \iffn@\else
    {\noexpands@
    \xdef\Next@{{\Thelabel@}{\Thelabel@@}{\Thelabel@@@}{\Thelabel@@@@}}%
    }%
    \expandafter\rightappend@\Next@\to\justfootmarklist@
   \fi
  \fi
  \@sf\relax}%
 \FN@\next@}
\def\footmark@{%
 \iffirstchoice@
  \global\advance\footmark@C\@ne
  \ifplainfn@
   \xdef\adjustedfootmark@{\number\footmark@C}%
  \else
   {\let\\\or\xdef\Next@{\ifcase\number\footmark@C\fnpages@\else
     -\@M\fi}}%
   \ifnum\Next@=-\@M
    \xdef\adjustedfootmark@{\number\footmark@C}%
   \else
    \ifnum\Next@=\lastfnpage@
     \global\advance\fancyfootmarkcount@\@ne
    \else
     \global\fancyfootmarkcount@\@ne
     \global\lastfnpage@\Next@
    \fi
    \xdef\adjustedfootmark@{\number\fancyfootmarkcount@}%
   \fi
  \fi
  {\noexpands@
  \xdef\Thelabel@@@{\adjustedfootmark@}%
  \xdefThelabel@\footmark@N
  \xdef\Thelabel@@@@{\Thelabel@}%
  \xdefThelabel@@\foottext@S
  }%
  \iffn@\else
   {\noexpands@
   \xdef\Next@{{\Thelabel@}{\Thelabel@@}{\Thelabel@@@}{\Thelabel@@@@}}%
   }%
   \expandafter\rightappend@\Next@\to\justfootmarklist@
  \fi
  \ifplainfn@
  \else
   \edef\next@{\write\laxwrite@{F\noexpand\the\pageno}}\next@
  \fi
 \fi
 \footmark@S{\footmark@N{\adjustedfootmark@}}%
 \@sf\relax}
\def\foottext{\prevanish@
 \ifx\justfootmarklist@\empty
  \Err@{There is no \noexpand\footmark for this \string\foottext}\fi
 \DN@\\##1##2\next@{\DN@{##1}\gdef\justfootmarklist@{##2}}%
 \expandafter\next@\justfootmarklist@\next@
 \expandafter\foottext@\next@}
\def\foottext@#1#2#3#4{{\noexpands@
  \xdef\Thelabel@{#1}\xdef\Thelabel@@{#2}%
  \xdef\Thelabel@@@{#3}\xdef\Thelabel@@@@{#4}}%
  \vfootnote@{\thelabel@@}}
\rightadd@\foottext\to\vanishlist@
\def\footnote{\fn@true
 \let\@sf\empty
 \ifhmode\edef\@sf{\spacefactor\the\spacefactor}\/\fi
 \DN@{\ifx"\next\expandafter\nextii@\else\expandafter\nextiii@\fi}%
 \DNii@"##1"{\footmark"##1"\vfootnote@{\let\style\foottext@S
  \let\numstyle\footmark@N##1}}%
 \def\nextiii@{\footmark\vfootnote@{\foottext@S{\footmark@N
  {\adjustedfootmark@}}}}%
 \FN@\next@}
\newdimen\litindent
\litindent20\p@
\newbox\litbox@
\newbox\Litbox@
\newcount\interlitpenalty@
\interlitpenalty@\@M
\newcount\litlines@
{\obeyspaces\gdef\defspace@{\def {\allowbreak\hskip.5emminus.15em}}}
{\obeylines\gdef\letM@{\let^^M\CtrlM@}}
\def\CtrlM@{\egroup
 \ifcase\litlines@\advance\litlines@\@ne\or
 \box\litbox@\advance\litlines@\@ne\else
 \penalty\interlitpenalty@\box\litbox@\fi
 \Lit@}
\def\Lit@{\setbox\litbox@\hbox\bgroup\litdefs@\hskip\litindent}
\newcount\littab@
\littab@8
\def\littab#1{\littab@#1\relax}
{\catcode`\^^I=\active\gdef\letTAB@{\let^^I\TAB@}}
\def\TAB@{\egroup
 \dimen@\wd\litbox@
 \advance\dimen@-\litindent
 \setboxz@h{\tt0}%
 \dimen@ii\littab@\wdz@
 \divide\dimen@\dimen@ii
 \multiply\dimen@\dimen@ii
 \advance\dimen@\littab@\wdz@
 \advance\dimen@\litindent
 \setbox\litbox@\hbox\bgroup\litdefs@\hbox to\dimen@{\unhbox\litbox@\hfil}}
{\catcode`\`=\active\gdef`{\relax\lq}}
\let\litbs@\relax
\let\litbs@@\relax
\def\litbackslash#1{%
 \edef\litbs@{\catcode`\string#1=\z@
 \def\noexpand\litbs@@{\def\expandafter\noexpand\csname\string#1\endcsname
  {\char`\string#1}}}}
\def\litcodes@{\catcode`\\=12
 \catcode`\{=12 \catcode`\}=12
 \catcode`\$=12 \catcode`\&=12
 \catcode`\#=12
 \catcode`\^=12 \catcode`\_=12
 \catcode`\@=12 \catcode`\~=12 \catcode`\"=12
 \catcode`\;=12 \catcode`\:=12 \catcode`\!=12 \catcode`\?=12
 \catcode`\%=12 \litbs@\catcode`\`=\active\obeyspaces\defspace@}
\def\activate@#1#2{{\lccode`\~=`#2%
 \lowercase{%
  \if0#1%
  \gdef\Next@{\def~{\egroup\endgroup\bigskip\vskip-\parskip
   \def\next@{\noindent@@\FN@\pretendspace@}\FNSS@\next@}}\else
  \gdef\Next@{\def~{\egroup\egroup\endgroup}}\fi
  }%
 }}
\def\litdefs@{\let\0\empty\let\1\litdelim@\def\ {\char32 }\litbs@@}%
\def\litdelimiter#1{%
 \edef\litdelim@{\char`#1}%
 \def\lit#1{\leavevmode\begingroup\litcodes@\litdefs@
  \tt\hyphenchar\tentt\m@ne\lit@}%
 \def\lit@##1#1{##1\endgroup\null}%
 \def\Lit#1{\ifhmode$$\abovedisplayskip\bigskipamount
  \abovedisplayshortskip\bigskipamount
  \belowdisplayskip\z@\belowdisplayshortskip\z@
  \postdisplaypenalty\@M
  $$\vskip-\baselineskip\else\bigskip\fi
  \begingroup\litlines@\z@
  \catcode`#1=\active\activate@0#1\Next@
  \def\displaybreak{\egroup\break\litlines@\z@\Lit@}%
  \def\allowdisplaybreak{\egroup\allowbreak\litlines@\z@\Lit@}%
  \def\allowdisplaybreaks{\egroup\allowbreak\interlitpenalty@\z@
   \litlines@\z@\Lit@}%
  \litcodes@\tt\catcode`\^^I=\active\letTAB@
  \obeylines\letM@\Lit@}%
 \def\Litbox##1=#1{\begingroup\ifodd##1\relax\aftergroup\global\fi
  \aftergroup\setbox\aftergroup##1\aftergroup\box\aftergroup\Litbox@
  \def\allowdisplaybreak{\egroup\allowbreak\litlines@\z@\Lit@}%
  \def\allowdisplaybreaks{\egroup\allowbreak\interlitpenalty@\z@
   \litlines@\z@\Lit@}%
  \catcode`#1=\active\activate@1#1\Next@
  \litcodes@\tt\catcode`\^^I=\active\letTAB@
  \obeylines\letM@\global\setbox\Litbox@\vbox\bgroup\litindent\z@%
  \litlines@\z@\Lit@}%
}
\newbox\titlebox@
\setbox\titlebox@\vbox{}
\rightadd@\title\to\overlonglist@
\def\title{\begingroup\Let@
 \global\setbox\titlebox@\vbox\bgroup\tabskip\hss@
 \halign to\hsize\bgroup\bf\hfil\ignorespaces##\unskip\hfil\cr}
\def\endtitle{\crcr\egroup\egroup\endgroup\overlong@false}
\newbox\authorbox@
\rightadd@\author\to\overlonglist@
\def\author{\begingroup\Let@
 \global\setbox\authorbox@\vbox\bgroup\tabskip\hss@
 \halign to\hsize\bgroup\rm\hfil\ignorespaces##\unskip\hfil\cr}
\def\endauthor{\crcr\egroup\egroup\endgroup\overlong@false}
\newbox\affilbox@
\def\affil{\begingroup\Let@
 \global\setbox\affilbox@\vbox\bgroup\tabskip\hss@
 \halign to\hsize\bgroup\rm\hfil\ignorespaces##\unskip\hfil\cr}%
\def\endaffil{\crcr\egroup\egroup\endgroup\overlong@false}
\let\date@\relax
\def\date#1{\gdef\date@{\ignorespaces#1\unskip}}
\def\today{\ifcase\month\or January\or February\or March\or April\or May\or
 June\or July\or August\or September\or October\or November\or December\fi
 \space\number\day, \number\year}
\def\maketitle{\hrule\height\z@\vskip-\topskip
 \vskip24\p@ plus12\p@ minus12\p@
 \unvbox\titlebox@
 \ifvoid\authorbox@\else\vskip12\p@ plus6\p@ minus3\p@\unvbox\authorbox@\fi
 \ifvoid\affilbox@\else\vskip10\p@ plus5\p@ minus2\p@\unvbox\affilbox@\fi
 \ifx\date@\relax\else\vskip6\p@ plus2\p@ minus\p@\centerline{\rm\date@}\fi
 \vskip18\p@ plus12\p@ minus6\p@}
\def\cite{%
 \DNii@(##1)##2{{\rm[}{##2}, {##1\/}{\rm]}}%
 \def\nextiii@##1{{\rm[}{##1\/}{\rm]}}%
 \DN@{\ifx\next(\expandafter\nextii@\else\expandafter\nextiii@\fi}%
 \FN@\next@}
\def\makebib@W{Bibliography}
\def\makebib{\begingroup\rm\bigbreak\centerline{\smc\makebib@W}%
 \nobreak\medskip
 \sfcode`\.=\@m\everypar{}\parindent\z@
 \def\nopunct{\nopunct@true}\def\nospace{\nospace@true}%
 \nopunct@false\nospace@false
 \def\lkerns@{\null\kern\m@ne sp\kern\@ne sp}%
 \def\nkerns@{\null\kern-\tw@ sp\kern\tw@ sp}%
}
\let\endmakebib\endgroup
\newif\ifnoprepunct@
\newif\ifnoprespace@
\newif\ifnoquotes@
\def\noprepunct{\noprepunct@true}
\def\noprespace{\noprespace@true}
\def\noquotes{\noquotes@true}
\newbox\nobox@
\newbox\keybox@
\newbox\bybox@
\newbox\paperbox@
\newbox\paperinfobox@
\newbox\jourbox@
\newbox\volbox@
\newbox\issuebox@
\newbox\yrbox@
\newbox\pgbox@
\newbox\ppbox@
\newbox\bookbox@
\newbox\inbookbox@
\newbox\bookinfobox@
\newbox\publbox@
\newbox\publaddrbox@
\newbox\edbox@
\newbox\edsbox@
\newbox\langbox@
\newbox\translbox@
\newbox\finalinfobox@
\def\setbibinfo@#1{\edef\next@{\ifnopunct@1\else0\fi
 \ifnospace@1\else0\fi\ifnoprepunct@1\else0\fi\ifnoprespace@1\else0\fi
 \ifnoquotes@1\else0\fi}%
 \DNii@{00000}%
 \ifx\next@\nextii@\else\xdef\bibinfo@{\bibinfo@\the#1,\next@}%
 \fi}
\def\getbibinfo@#1{\ifx\bibinfo@\empty
 \let\next@0\let\nextii@0\let\nextiii@0\let\nextiv@0\let\nextv@0\else
 \edef\next@{\def
  \noexpand\next@####1\the#1,####2####3####4####5####6####7\noexpand\next@
  {\let\noexpand\next@####2\let\noexpand\nextii@####3%
  \let\noexpand\nextiii@####4\let\noexpand\nextiv@####5%
  \let\noexpand\nextv@####6}%
  \noexpand\next@\bibinfo@\the#1,00000\noexpand\next@}\next@
 \fi}
\newif\ifbookinquotes@
\def\bookinquotes{\bookinquotes@true}
\newif\ifpaperinquotes@
\def\paperinquotes{\paperinquotes@true}
\newif\ifininbook@
\def\ininbook{\ininbook@true}
\newif\ifopenquotes@
\def\closequotes@{\ifopenquotes@''\openquotes@false\fi}
\newif\ifbeginbib@
\newif\ifendbib@
\newif\ifprevjour@
\newif\ifprevbook@
\newdimen\bibindent@
\bibindent@20\p@
\def\bib{\global\let\bibinfo@\empty\global\let\translinfo@\relax\beginbib@true
 \begingroup\noindent@
 \hangindent\bibindent@\hangafter\@ne\bib@}
\def\v@id#1{\setbox#1\box\voidb@x}
\def\bib@{\v@id\nobox@\v@id\keybox@\v@id\bybox@\v@id\paperbox@
 \v@id\paperinfobox@\v@id\jourbox@\v@id\volbox@\v@id\issuebox@
 \v@id\yrbox@\v@id\pgbox@\v@id\ppbox@\v@id\bookbox@\v@id\inbookbox@
 \v@id\bookinfobox@\v@id\publbox@\v@id\publaddrbox@\v@id\edbox@
 \v@id\edsbox@\v@id\langbox@\v@id\translbox@\v@id\finalinfobox@
 \bgroup}
\def\Setnonemptybox@#1#2{\unskip\setbibinfo@#1\egroup#2%
 \def\aftergroup@{\ifdim\wd#1=\z@\setbox#1\box\voidb@x\fi}%
 \setbox#1\vbox\bgroup\aftergroup\aftergroup@\hsize\maxdimen\leftskip\z@
 \rightskip\z@\hbadness\@M\hfuzz\maxdimen\noindent}
\def\setnonemptybox@#1{\Setnonemptybox@#1\relax}
\def\no{\setnonemptybox@\nobox@}
\def\key{\setnonemptybox@\keybox@\bf}
\def\by{\setnonemptybox@\bybox@}
\def\bysame{\setnonemptybox@\bybox@\leaders\hrule\hskip3em\null}
\def\paper{\setnonemptybox@\paperbox@
 \ifpaperinquotes@\getbibinfo@\paperbox@
 \if\nextv@1\else``\fi\else\it\fi}
\def\paperinfo{\setnonemptybox@\paperinfobox@}
\def\jour{\Setnonemptybox@\jourbox@\prevjour@true}
\def\vol{\setnonemptybox@\volbox@\bf}
\def\issue{\setnonemptybox@\issuebox@}
\def\yr{\setnonemptybox@\yrbox@}

\def\pg{\setnonemptybox@\pgbox@}
\def\pp{\setnonemptybox@\ppbox@}
\def\book{\Setnonemptybox@\bookbox@\prevbook@true
 \ifbookinquotes@\getbibinfo@\bookbox@
 \if\nextv@1\else``\fi\else\it\fi}
\def\inbook{\Setnonemptybox@\inbookbox@\prevbook@true
 \ifininbook@ in \fi\ifbookinquotes@\getbibinfo@\inbookbox@
 \if\nextv@1\else``\fi\fi}
\def\bookinfo{\setnonemptybox@\bookinfobox@}
\def\publ{\setnonemptybox@\publbox@}
\def\publaddr{\setnonemptybox@\publaddrbox@}
\def\ed{\setnonemptybox@\edbox@}
\def\eds{\setnonemptybox@\edsbox@}
\def\lang{\setnonemptybox@\langbox@}
\def\finalinfo{\setnonemptybox@\finalinfobox@}
\def\setboxzl@{\setbox\z@\lastbox}
\def\getbox@#1{\setbox\z@\vbox{\vskip-\@M\p@
 \unvbox#1%
 \setboxzl@
 \global\setbox\@ne\hbox{\unhbox\z@\unskip\unskip\unpenalty}%
 \ifdim\lastskip=-\@M\p@\else
 \loop\ifdim\lastskip=-\@M\p@
 \else\unskip\unpenalty\setboxzl@
 \global\setbox\@ne\hbox{\unhbox\z@\unhbox\@ne}%
 \repeat\fi}%
 \unhbox\@ne}
\def\adjustpunct@#1{\count@\lastkern
 \ifnum\count@=\z@#1\closequotes@\else
 \ifnum\count@>\tw@#1\closequotes@\else
 \ifnum\count@<-\tw@#1\closequotes@\else
  \unkern\unkern\setboxzl@
  \skip@\lastskip\unskip
  \count@@\lastpenalty\unpenalty
  \ifnum\count@=\tw@\unskip\setboxzl@\fi
  \ifdim\skip@=\z@\else\hskip\skip@\fi
  #1\closequotes@
  \ifnum\count@=\tw@\null\hfill\fi
  \penalty\count@@
 \fi\fi\fi}
\def\prepunct@#1#2{\getbibinfo@#2%
 \ifnopunct@
 \else
  \if\nextiii@0\adjustpunct@#1\fi
 \fi
 \closequotes@
 \ifnospace@
 \else
  \if\nextiv@0\space\else\fi
 \fi
 \nopunct@false\nospace@false
 \if\next@1\nopunct@true\fi
 \if\nextii@1\nospace@true\fi}
\def\ppunbox@#1#2{\prepunct@{#1}#2%
 \getbox@#2}
\let\semicolon@;
\def\endbib@{%
 \ifbeginbib@
  \ifvoid\nobox@
   \ifvoid\keybox@\else\hbox to\bibindent@{[\getbox@\keybox@]\hss}\fi
  \else\hbox to\bibindent@{\hss\getbox@\nobox@. }\fi
  \ifvoid\bybox@\else\getbox@\bybox@\fi
 \else
  \nopunct@true
  \ifvoid\bybox@\else\ppunbox@\relax\bybox@\fi
 \fi
 \ifvoid\translbox@\else\ppunbox@,\translbox@\fi
 \ifvoid\paperbox@\else\ppunbox@,\paperbox@\ifpaperinquotes@
  \if\nextv@1\else\openquotes@true\fi\fi
 \fi
 \ifvoid\paperinfobox@\else\ppunbox@,\paperinfobox@\fi
 \test@false
 \ifvoid\jourbox@\else\test@true\ppunbox@,\jourbox@\fi
 \ifprevjour@\test@true\fi
 \iftest@
  \ifvoid\volbox@\else\ppunbox@\relax\volbox@\fi
  \ifvoid\issuebox@
   \else\prepunct@\relax\issuebox@ no.~\getbox@\issuebox@\fi
  \ifvoid\yrbox@\else\prepunct@\relax\yrbox@(\getbox@\yrbox@)\fi
  \ifvoid\ppbox@\else\ppunbox@,\ppbox@\fi
  \ifvoid\pgbox@\else\prepunct@,\pgbox@ p.~\getbox@\pgbox@\fi
 \fi
 \test@false
 \ifvoid\bookbox@\else\test@true\ppunbox@,\bookbox@\ifbookinquotes@
  \if\nextv@1\else\openquotes@true\fi\fi\fi
 \ifvoid\inbookbox@\else\test@true\ppunbox@,\inbookbox@\ifbookinquotes@
  \if\nextv@1\else\openquotes@true\fi\fi\fi
 \ifprevbook@\test@true\fi
 \iftest@
  \ifvoid\edbox@\else\prepunct@\relax\edbox@(\getbox@\edbox@, ed.)\fi
  \ifvoid\edsbox@\else\prepunct@\relax\edsbox@(\getbox@\edsbox@, eds.)\fi
  \ifvoid\bookinfobox@\else\ppunbox@,\bookinfobox@\fi
  \ifvoid\publbox@\else\ppunbox@,\publbox@\fi
  \ifvoid\publaddrbox@\else\ppunbox@,\publaddrbox@\fi
  \ifvoid\yrbox@\else\ppunbox@,\yrbox@\fi
  \ifvoid\ppbox@\else\prepunct@,\ppbox@ pp.~\getbox@\ppbox@\fi
  \ifvoid\pgbox@\else\prepunct@,\pgbox@ p.~\getbox@\pgbox@\fi
 \fi
 \ifvoid\finalinfobox@
  \ifendbib@
   \ifnopunct@\else.\closequotes@\fi
  \else
  \ifvoid\langbox@\else\space(\getbox@\langbox@)\fi
   \/\semicolon@\closequotes@
  \fi
 \else
  \ifendbib@
   \ppunbox@{.\spacefactor3000\relax}\finalinfobox@
    \ifnopunct@\else.\fi
  \else
   \ppunbox@,\finalinfobox@\/\semicolon@\fi
 \fi
 \ifvoid\langbox@\else\space(\getbox@\langbox@)\fi
}
\def\endbib{\unskip\egroup\endbib@true\endbib@\par\endgroup}
\def\morebib{\unskip\egroup
 \endbib@false\endbib@
 \global\let\bibinfo@\empty\beginbib@false
 \bib@}
\def\anotherbib{\unskip\egroup
 \endbib@false\endbib@
 \global\let\bibinfo@\empty\beginbib@false
 \prevjour@false\prevbook@false\bib@}
\def\transl{\unskip
 \xdef\translinfo@{\the\translbox@,\ifnopunct@1\else0\fi
 \ifnospace@1\else0\fi\ifnoprepunct@1\else0\fi\ifnoprespace@1\else0\fi0}%
 \egroup\endbib@false\endbib@
 \global\let\bibinfo@\translinfo@\beginbib@false
 \bib@
 \egroup
 \def\aftergroup@{\ifdim\wd\translbox@=\z@\setbox\translbox@\box\voidb@x\fi}%
 \setbox\translbox@\vbox\bgroup\aftergroup\aftergroup@
 \hsize\maxdimen\leftskip\z@\rightskip\z@\hbadness\@M\hfuzz\maxdimen
 \noindent}
\newwrite\auxwrite@
\newread\bbl@
\def\UseBibTeX{\immediate\openout\auxwrite@=\jobname.aux
 \let\cite\BTcite@
 \def\nocite##1{\immediate\write\auxwrite@{\string\citation{##1}}}%
 \def\bibliographystyle##1{\immediate\write\auxwrite@{\string
  \bibstyle{##1}}}%
 \def\bibliography@W{Bibliography}%
 \def\bibliography##1{\immediate\write\auxwrite@{\string\bibdata{##1}}%
  \immediate\openin\bbl@=\jobname.bbl
  \ifeof\bbl@
   \W@{No .bbl file}%
  \else
   \immediate\closein\bbl@
   \begingroup\input bibtex \input\jobname.bbl \endgroup
  \fi}%
 }
\def\BTcite@{%
 \DNii@(##1)##2{{\rm[}\BTcite@@##2,\BTcite@@{\rm, }{##1\/}{\rm]}%
  \immediate\write\auxwrite@{\string\citation{##2}}}%
 \def\nextiii@##1{{\rm[}\BTcite@@##1,\BTcite@@\/{\rm]}%
  \immediate\write\auxwrite@{\string\citation{##1}}}%
 \DN@{\ifx\next(\expandafter\nextii@\else\expandafter\nextiii@\fi}%
 \FN@\next@}%
\def\BTcite@@#1,{\BTcite@@@{#1}\FN@\BTcite@@@@}
\def\BTcite@@@@{\ifx\next\BTcite@@
 \expandafter\eat@\else{\rm, }\expandafter\BTcite@@\fi}
\catcode`\~=11
\def\BTcite@@@#1{\nolabel@\cite{#1}\relax
 \DNii@##1~##2\nextii@{##1}%
 \csL@{#1}\expandafter\nextii@\Next@\nextii@\fi}
\catcode`\~=\active

\def\beginthebibliography@#1{\rm\setboxz@h{#1\ }\bibindent@\wdz@
 \bigbreak\centerline{\smc\bibliography@W}\nobreak\medskip
 \sfcode`\.=\@m\everypar{}\parindent\z@}

\newif\iffigproofing@
\def\Figureproofing{\figproofing@true}
\def\noFigureproofing{\figproofing@false}
\newif\ifHby@
\def\Hbyw#1{\global\Hby@true\hbyw\vsize{#1}}
\def\hbyw#1#2{%
 \hbox{%
  \ifHby@
  \else
   \iffigproofing@
    \setbox\z@\vbox{\hrule\width5\p@}\ht\z@\z@
    \vbox to#1{\hrule\height5\p@\width.4\p@\vfil\hrule\height5\p@\width.4\p@}%
    \kern-.4\p@\rlap{\copy\z@}\raise#1\hbox{\rlap{\copy\z@}}%
   \fi
  \fi
  \vbox to#1{\hbox to#2{}\vfil}%
  \ifHby@
  \else
   \iffigproofing@
    \vbox to#1{\hrule\height5\p@\width.4\p@\vfil\hrule\height5\p@\width.4\p@}%
    \kern-.4\p@\llap{\copy\z@}\raise#1\hbox{\llap{\boxz@}}%
   \fi
  \fi}}
\newcount\island@C
\let\island@P\empty
\let\island@Q\empty
\def\island@S#1{#1\null.}
\let\island@N\arabic
\def\island@F{\rm}
\def\island@@@P{\csname\exxx@\islandtype@ @P\endcsname}
\def\island@@@Q{\csname\exxx@\islandtype@ @Q\endcsname}
\def\island@@@S{\csname\exxx@\islandtype@ @S\endcsname}
\def\island@@@N{\csname\exxx@\islandtype@ @N\endcsname}
\def\island@@@F{\csname\exxx@\islandtype@ @F\endcsname}
\def\island@@@C{\csname island@C\islandclass@\endcsname}
\newif\ifplace@
\newif\ifisland@
\def\island{%
 \ifplace@
  \DN@{\let\islandclass@\empty\def\islandtype@{\island}\FN@\island@}%
 \else
  \long\DN@##1\endisland{\Err@{\noexpand\island must be used after some
   type of \string\...place}}%
 \fi
 \next@}
\def\island@{\ifx\next\c\let\next@\island@c\else
 \DN@{\FN@\island@@}\fi\next@}
\def\island@@{\ifcat\bgroup\noexpand\next\let\next@\island@@@\else
 \DN@{\Err@{\noexpand\island must be followed by a {prefix} for
 \string\caption's}}\fi\next@}
\newbox\islandbox@
\newcount\captioncount@
\def\island@@@#1{\def\captionprefix@{#1}\captioncount@\z@
 \global\setbox\islandbox@\vbox\bgroup}
\def\island@c\c#1{%
 \ifplace@
 \DN@{\def\islandclass@{#1}%
  \expandafter\ifx\csname island@C#1\endcsname\relax
  \expandafter\newcount@\csname island@C#1\endcsname
   \global\csname island@C#1\endcsname\z@\fi
  \FNSS@\island@c@}%
 \else
 \DN@{\edef\next@{\long\def\noexpand\next@########1\expandafter\noexpand
  \csname end\exxx@\islandtype@\endcsname{\noexpand\Err@{\noexpand\noexpand
  \expandafter\noexpand
  \islandtype@ must be used after some type of \noexpand\string
   \noexpand\...place}}}\next@\next@}%
 \fi
 \next@}
\def\island@c@{%
 \ifcat\bgroup\noexpand\next
  \let\next@\island@c@@
 \else
  \DN@{\Err@{\noexpand\island\string\c{\expandafter\string\islandclass@} must
   be followed by a {prefix} for \string\caption's}}%
 \fi\next@}
\def\island@c@@#1{\def\captionprefix@{#1}%
 \captioncount@\z@\global\setbox\islandbox@\vbox\bgroup}
\rightadd@\caption\to\nofrillslist@
\newbox\captionbox@
\newbox\Captionbox@
\def\caption{%
 \ifnum\captioncount@=\z@
  \ifnopunct@
   \DN@{\egroup\nopunct@true}%
  \else
   \let\next@\egroup
  \fi
 \else
  \let\next@\relax
 \fi
 \next@
 \advance\captioncount@\@ne
 \FN@\caption@}
\def\caption@{\ifx\next"\expandafter\caption@q\else\expandafter\caption@@\fi}
\def\caption@q"#1"{\quoted@true
 {\noexpands@
 \let\pre\island@@@P\let\post\island@@@Q
 \let\style\island@@@S\let\numstyle\island@@@N
 \Qlabel@{#1}\let\style\relax\xdef\Qlabel@@@@{#1}}%
 \finishcaption@}
\def\caption@@{\quoted@false
 \global\advance\island@@@C\@ne
 {\noexpands@
 \xdef\Thelabel@@@{\number\island@@@C}%
 \xdefThelabel@\island@@@N
 \xdef\Thelabel@@@@{\island@@@P\Thelabel@\island@@@Q}%
 \xdefThelabel@@\island@@@S
 \xdef\Thepref@{\Thelabel@@@@}}%
 \finishcaption@}
\long\def\captionformat@#1#2#3{\rm\strut#1 {\island@@@F#2} #3%
 \punct@.\strut}
\long\def\widerthanisland@#1#2#3{\test@true\setbox\z@\vbox{\hsize\maxdimen
 \noindent@@\captionformat@{#1}{#2}{#3}\par\setboxzl@}%
 \ifdim\wdz@=\z@
  \global\setbox\captionbox@\hbox{\noset@\unlabel@
   \captionformat@{#1}{#2}{#3}}%
  \ifdim\wd\captionbox@>\wd\islandbox@\else\test@false\fi
 \fi}
\long\def\captionformat@@#1#2#3{\widerthanisland@{#1}{#2}{#3}%
 \iftest@
  \global\setbox\captionbox@\vbox{\hsize\wd\islandbox@
   \vskip-\parskip\noindent@@\noset@\unlabel@
   \captionformat@{#1}{#2}{#3}\par}%
 \else
  \global\setbox\captionbox@
   \hbox to\wd\islandbox@{\hfil\box\captionbox@\hfil}%
 \fi}
\long\def\finishcaption@#1{\def\entry@{#1}%
 {\locallabel@
 \captionformat@@
  {\expandafter\ignorespaces\captionprefix@\unskip}%
  {\ifx\thelabel@@\empty\unskip\else\thelabel@@\fi}%
  {\ignorespaces#1\unskip}%
 \ifnum\captioncount@=\@ne
  \global\setbox\islandbox@\vbox{\ticwrite@\vbox{\box\islandbox@}}%
  \global\setbox\Captionbox@\vbox{\box\captionbox@}%
 \else
  \global\setbox\islandbox@\vbox{\unvbox\islandbox@\setboxzl@
   \ticwrite@\boxz@}%
  \global\setbox\Captionbox@\vbox{\unvbox\Captionbox@
   \smallskip\box\captionbox@}%
 \fi}%
 \nopunct@false\nospace@false\ignorespaces}
\def\Sixtic@{\ifx\macdef@\empty\else
 \DN@##1##2\next@{\def\macdef@{##1##2}}%
 \expandafter\next@\macdef@\next@
 \edef\next@
  {\noexpand\six@\tic@\macdef@
  \space\space\space\space\space\space\space\space\space\space\space\space
  \noexpand\six@}%
 \next@\let\macdef@\relax\fi}
\def\ticwrite@{%
 \iftoc@
  {\noexpands@\let\style\relax
  \DN@{\island}%
  \edef\next@{\write\tic@{%
   \ifnopunct@\noexpand\noexpand\noexpand\nopunct\fi
   \ifx\islandtype@\next@\noexpand\noexpand\noexpand\island
    \noexpand\string\noexpand\c{\islandclass@}{\captionprefix@}%
     {\QorThelabel@@@@}\else\noexpand\noexpand\expandafter\noexpand
     \islandtype@{\QorThelabel@@@@}}\fi}%
  \next@}%
  \expandafter\unmacro@\meaning\entry@\unmacro@
  \Sixtic@
  \write\tic@{\noexpand\Page{\number\pageno}{\page@N}{\page@P}{\page@Q}^^J}%
 \fi}
\def\Htrim@#1{%
 \ifHby@
  \dimen@\vsize
  \ifnum\captioncount@=\z@
  \else
   \advance\dimen@-\ht\Captionbox@
   \advance\dimen@-#1%
  \fi
  \global\Hby@false
  \dimen@ii\wd\islandbox@
  \global\setbox\islandbox@\vbox
   {\unvbox\islandbox@\setboxzl@
   \vbox to\z@{\vss\boxz@}\nointerlineskip\hbyw\dimen@\dimen@ii}%
  \global\Hby@true
 \fi}
\newif\ifdata@
\def\iclasstest@#1{\DN@{#1}\ifx\next@\islandclass@
 \test@true\else\test@false\fi}
\skipdef\skipi@=1
\def\endisland{\ifnum\captioncount@=\z@\expandafter\egroup\fi
 \ifdata@
 \else
  \iclasstest@{T}%
  \iftest@
   {\rm\global\skipi@-\dp\strutbox}\global\advance\skipi@\bigskipamount
   \Htrim@\skipi@
   \global\setbox\islandbox@\vbox
    {\ifnum\captioncount@=\z@\else
     \box\Captionbox@
     \nointerlineskip
     \vskip\skipi@\fi
     \box\islandbox@}%
  \else
   {\rm\global\skipi@\dp\strutbox}\global\advance\skipi@\medskipamount
   \Htrim@\skipi@
   \global\setbox\islandbox@\vbox
    {\box\islandbox@
     \ifnum\captioncount@=\z@\else
     \nointerlineskip
     \vskip\skipi@
     \box\Captionbox@
     \fi}%
  \fi
  \ifHby@
  \else
   \dimen@\ht\islandbox@\advance\dimen@\dp\islandbox@
   \ifdim\dimen@>\vsize
    \DN@{\island}%
    \Err@{%
     \ifx\islandtype@\next@\noexpand\island\else
      \expandafter\noexpand\islandtype@\fi
     \ifnum\captioncount@=\z@\else
       with \noexpand\caption\fi
      is larger than page}%
     \ht\islandbox@=\vsize
   \fi
  \fi
 \fi
 \global\Hby@false\island@true}
\def\newisland#1\c#2#3{\define#1{}%
 \iftoc@\immediate\write\tic@{\noexpand\newisland\noexpand#1%
  \string\c{#2}{#3}^^J}\fi
 \expandafter\def\csname\exstring@#1@S\endcsname{\island@S}%
 \expandafter\def\csname\exstring@#1@N\endcsname{\island@N}%
 \expandafter\def\csname\exstring@#1@P\endcsname{\island@P}%
 \expandafter\def\csname\exstring@#1@Q\endcsname{\island@Q}%
 \expandafter\def\csname\exstring@#1@F\endcsname{\island@F}%
 \expandafter\def\csname end\exstring@#1\endcsname{\endisland}%
 \expandafter
 \ifx\csname island@C#2\endcsname\relax
  \expandafter\newcount@\csname island@C#2\endcsname
  \global\csname island@C#2\endcsname\z@
 \fi
 \edef\next@{\noexpand\expandafter\noexpand\let\noexpand
  \csname\exstring@#1@C\noexpand\endcsname
  \csname island@C#2\endcsname}%
 \next@
 \def#1{\def\islandtype@{#1}\island@c\c{#2}{#3}}}
\newisland\Figure\c{F}{Figure}
\newisland\Table\c{T}{Table}
\newbox\islandboxi
\newbox\islandboxii
\newbox\islandboxiii
\newbox\captionboxi
\newbox\captionboxii
\newbox\captionboxiii
\long\def\islandpairdata#1#2{{\data@true
 \place@true
 #1%
 \global\setbox\islandboxi\box\islandbox@
 \global\setbox\captionboxi\box\Captionbox@
 #2%
 \global\setbox\islandboxii\box\islandbox@
 \global\setbox\captionboxii\box\Captionbox@
 }}
\long\def\islandpairbox#1#2{\islandpairdata{#1}{#2}%
 \dimen@\ht\captionboxi
 \ifdim\ht\captionboxii>\dimen@\dimen@\ht\captionboxii\fi
 \ifdim\dimen@>\z@
  \ifdim\ht\captionboxi<\dimen@
   \global\setbox\captionboxi\vbox to\dimen@{\unvbox\captionboxi\vfil}\fi
  \ifdim\ht\captionboxii<\dimen@
   \global\setbox\captionboxii\vbox to\dimen@{\unvbox\captionboxii\vfil}\fi
 \fi
 \global\setbox\islandbox@\vbox
 {\hbox to\hsize{\hfil\box\islandboxi\hfil\box\islandboxii\hfil}%
 \ifdim\dimen@>\z@\nointerlineskip
 {\rm\global\skipi@\dp\strutbox}\global\advance\skipi@\medskipamount
  \vskip\skipi@
  \hbox to\hsize{\hfil\box\captionboxi\hfil\box\captionboxii\hfil}\fi}}	
\long\def\islandpairboxa#1#2{\islandpairdata{#1}{#2}%
 \dimen@\ht\captionboxi
 \ifdim\ht\captionboxii>\dimen@\dimen@\ht\captionboxii\fi
 \ifdim\dimen@>\z@
  \ifdim\ht\captionboxi<\dimen@
   \global\setbox\captionboxi\vbox to\dimen@{\vfil\unvbox\captionboxi}\fi
  \ifdim\ht\captionboxii<\dimen@
   \global\setbox\captionboxii\vbox to\dimen@{\vfil\unvbox\captionboxii}\fi
 \fi
 \dimen@ii\ht\islandboxi
 \ifdim\ht\islandboxii>\dimen@ii \dimen@ii\ht\islandboxii\fi
 \ifdim\dimen@ii>\z@
  \ifdim\ht\islandboxi<\dimen@ii
   \global\setbox\islandboxi\vbox to\dimen@ii{\box\islandboxi\vfil}\fi
  \ifdim\ht\islandboxii<\dimen@ii
   \global\setbox\islandboxii\vbox to\dimen@ii{\box\islandboxii\vfil}\fi
 \fi
 \global\setbox\islandbox@\vbox{\ifdim\dimen@>\z@
  \hbox to\hsize{\hfil\box\captionboxi\hfil\box\captionboxii\hfil}%
  \nointerlineskip{\rm\global\skipi@-\dp\strutbox}%
  \global\advance\skipi@\bigskipamount\vskip\skipi@\fi
  \hbox to\hsize{\hfil\box\islandboxi\hfil\box\islandboxii\hfil}}}
\long\def\islandtripledata#1#2#3{{\data@true\place@true
 #1%
 \global\setbox\islandboxi\box\islandbox@
 \global\setbox\captionboxi\box\Captionbox@
 #2%
 \global\setbox\islandboxii\box\islandbox@
 \global\setbox\captionboxii\box\Captionbox@
 #3%
 \global\setbox\islandboxiii\box\islandbox@
 \global\setbox\captionboxiii\box\Captionbox@
 }}
\long\def\islandtriplebox#1#2#3{\islandtripledata{#1}{#2}{#3}%
 \dimen@\ht\captionboxi
 \ifdim\ht\captionboxii>\dimen@ \dimen@\ht\captionboxii\fi
 \ifdim\ht\captionboxiii>\dimen@ \dimen@\ht\captionboxiii\fi
 \ifdim\dimen@>\z@
  \ifdim\ht\captionboxi<\dimen@
   \global\setbox\captionboxi\vbox to\dimen@{\unvbox\captionboxi\vfil}\fi
  \ifdim\ht\captionboxii<\dimen@
   \global\setbox\captionboxii\vbox to\dimen@{\unvbox\captionboxii\vfil}\fi
  \ifdim\ht\captionboxiii<\dimen@
   \global\setbox\captionboxiii\vbox to\dimen@{\unvbox\captionboxiii\vfil}\fi
 \fi
 \global\setbox\islandbox@\vbox
  {\hbox to\hsize{\hfil\box\islandboxi\hfil\box\islandboxii\hfil
   \box\islandboxiii\hfil}%
 \ifdim\dimen@>\z@\nointerlineskip
  {\rm\global\skipi@\dp\strutbox}\global\advance\skipi@\medskipamount
  \vskip\skipi@
  \hbox to\hsize{\hfil\box\captionboxi\hfil\box\captionboxii\hfil
   \box\captionboxiii\hfil}\fi}}
\def\islandtripleboxa#1#2#3{\islandtripledata{#1}{#2}{#3}%
 \dimen@\ht\captionboxi
 \ifdim\ht\captionboxii>\dimen@ \dimen@\ht\captionboxii\fi
 \ifdim\ht\captionboxiii>\dimen@ \dimen@\ht\captionboxiii\fi
 \ifdim\dimen@>\z@
  \ifdim\ht\captionboxi<\dimen@
   \global\setbox\captionboxi\vbox to\dimen@{\vfil\unvbox\captionboxi}\fi
  \ifdim\ht\captionboxii<\dimen@
   \global\setbox\captionboxii\vbox to\dimen@{\vfil\unvbox\captionboxii}\fi
  \ifdim\ht\captionboxiii<\dimen@
   \global\setbox\captionboxiii\vbox to\dimen@{\vfil\unvbox\captionboxiii}\fi
 \fi
 \dimen@ii\ht\islandboxi
 \ifdim\ht\islandboxii>\dimen@ii \dimen@ii\ht\islandboxii\fi
 \ifdim\ht\islandboxiii>\dimen@ii \dimen@ii\ht\islandboxiii\fi
 \ifdim\dimen@ii>\z@
  \ifdim\ht\islandboxi<\dimen@ii
   \global\setbox\islandboxi\vbox to\dimen@ii{\box\islandboxi\vfil}\fi
  \ifdim\ht\islandboxii<\dimen@ii
   \global\setbox\islandboxii\vbox to\dimen@ii{\box\islandboxii\vfil}\fi
  \ifdim\ht\islandboxiii<\dimen@ii
   \global\setbox\islandboxiii\vbox to\dimen@ii{\box\islandboxiii\vfil}\fi
 \fi
 \global\setbox\islandbox@\vbox
  {\ifdim\dimen@>\z@
  \hbox to\hsize{\hfil\box\captionboxi\hfil\box\captionboxii\hfil
   \box\captionboxiii\hfil}%
  \nointerlineskip{\rm\global\skipi@-\dp\strutbox}%
  \global\advance\skipi@\bigskipamount\vskip\skipi@\fi
  \hbox to\hsize{\hfil\box\islandboxi\hfil\box\islandboxii\hfil
   \box\islandboxiii\hfil}}}
\def\Figurepair#1\and#2\endFigurepair{\island@true
 \islandpairbox{\Figure#1\endFigure}{\Figure#2\endFigure}}
\def\Figuretriple#1\and#2\and#3\endFiguretriple{\island@true
 \islandtriplebox{\Figure#1\endFigure}{\Figure#2\endFigure}%
  {\Figure#3\endFigure}}
\def\Tablepair#1\and#2\endTablepair{\island@true
 \islandpairboxa{\Table#1\endTable}{\Table#2\endTable}}
\def\Tabletriple#1\and#2\and#3\endTabletriple{\island@true
 \islandtripleboxa{\Table#1\endTable}{\Table#2\endTable}%
 {\Table#3\endTable}}
\def\place#1{\place@true\island@false
 #1%
 \ifisland@
  \box\islandbox@
 \else
  \Err@{Whoa ... there's no \string\Figure, \string\Table,
   etc., here}%
 \fi
 \place@false}
\newskip\belowtopfigskip
\belowtopfigskip 15\p@ plus 5\p@ minus5\p@
\newskip\abovebotfigskip
\abovebotfigskip 18\p@ plus 6\p@ minus6\p@
\newdimen\minpagesize
\minpagesize 5pc
\dimen@\belowtopfigskip
\advance\dimen@-\abovebotfigskip
\skip\topins\dimen@
\dimen\topins\z@
\newcount\topinscount@
\newbox\topinsdims@
\def\storedim@{\global\setbox\topinsdims@
 \vbox{\hbox to\dimen@{}\unvbox\topinsdims@}}
\def\advancedimtopins@{%
 \ifnum\pageno=\@ne
 \else
   \advance\dimen@\dimen\topins
   \global\dimen\topins\dimen@
 \fi}
\newcount\flipcount@
\def\fliptopins@{%
 \global\flipcount@\z@
 \ifvoid\topins\else
 \setbox\z@\vbox
  {\vskip\p@
   \unvbox\topins
   \global\setbox\topins\vbox{}%
   \loop
    \test@false
    \ifdim\lastskip=\z@\unskip
     \ifdim\lastskip=\z@
      \test@true\fi\fi
    \iftest@
    \global\advance\flipcount@\@ne
    \setboxzl@
    \global\setbox\topins\vbox{\unvbox\topins\boxz@}%
    \unpenalty
   \repeat}\fi}
\newif\ifPar@
\newcount\Parcount@
\newbox\Parbox@
\expandafter\newbox\csname Parfigbox1\endcsname
\expandafter\newbox\csname Parfigbox2\endcsname
\expandafter\newbox\csname Parfigbox3\endcsname
\expandafter\newbox\csname Parfigbox4\endcsname
\expandafter\newbox\csname Parfigbox5\endcsname
\expandafter\newdimen\csname Parprev1\endcsname
\expandafter\newdimen\csname Parprev2\endcsname
\expandafter\newdimen\csname Parprev3\endcsname
\expandafter\newdimen\csname Parprev4\endcsname
\expandafter\newdimen\csname Parprev5\endcsname
\expandafter\newdimen\csname Parprev6\endcsname
\def\Par{\par\global\csname Parprev1\endcsname\prevdepth
 \global\Parcount@\@ne
 \global\Par@true\global\let\Parlist@\empty
 \global\setbox\Parbox@\vbox\bgroup\break}
\def\place@#1#2{%
 \ifisland@
  \ifhmode
   \ifPar@
    \ifnum\Parcount@>5
     \Err@{Only 5 \string\place's allowed per
      \string\Par...\noexpand\endPar paragraph}%
    \else
     \expandafter\expandafter\expandafter
      \global\expandafter\setbox
       \csname Parfigbox\number\Parcount@\endcsname\box\islandbox@
     \global\advance\Parcount@\@ne
     \xdef\Parlist@{\Parlist@#1}%
    \fi
   \else
    \vadjust{#2}%
   \fi
  \else
   #2%
  \fi
 \else
  \Err@{Whoa ... there's no \string\Figure,
   \string\Table, etc., here}%
 \fi
 \place@false}
\long\def\Aplace#1{\prevanish@
 \place@true\island@false
 #1%
 \place@ a\Aplace@
 \postvanish@}
\long\def\AAplace#1{\prevanish@\place@true\island@false
 #1%
 \place@ A\AAplace@
 \postvanish@}
\newif\ifAA@
\def\AAplace@{\AA@true\Aplace@\AA@false}
\let\AAlist@\empty
\def\Aplace@{\allowbreak
 \dimen@=\ht\islandbox@
 \advance\dimen@\abovebotfigskip
 \ht\islandbox@\dimen@
 \advance\dimen@\dp\islandbox@
 \storedim@
 \ifAA@
  \xdef\AAlist@{\AAlist@1}%
  \advancedimtopins@
 \else
  \xdef\AAlist@{\AAlist@0}%
  \ifnum\topinscount@>\@ne\else\advancedimtopins@\fi
 \fi
 \insert\topins{\penalty\z@\splittopskip\z@\floatingpenalty\z@
  \box\islandbox@}%
 \global\advance\topinscount@\@ne}
\long\def\Bplace#1{\prevanish@\place@true\island@false
 #1%
 \place@ b\Bplace@
 \postvanish@}
\def\Bplace@{\allowbreak
 \ifnum\topinscount@=\z@
  \setbox\z@\vbox{\vbox to-\belowtopfigskip{}}%
  \dimen@-\skip\topins
  \ht\z@\dimen@
  \storedim@
  \advancedimtopins@
  \insert\topins{\boxz@}%
  \global\advance\topinscount@\@ne
  \xdef\AAlist@{\AAlist@0}%
 \fi
 \dimen@\ht\islandbox@
 \advance\dimen@\abovebotfigskip
 \ht\islandbox@\dimen@
 \advance\dimen@\dp\islandbox@
 \storedim@
 \xdef\AAlist@{\AAlist@0}%
 \ifnum\topinscount@>\@ne\else\advancedimtopins@\fi
 \insert\topins{\penalty\z@\splittopskip\z@
  \floatingpenalty\z@
  \box\islandbox@}%
 \global\advance\topinscount@\@ne}
\def\breakisland@{\global\setbox\@ne\lastbox\global\skipi@\lastskip\unskip
 \global\setbox\thr@@\lastbox}%
\def\printisland@{\centerline{\box\thr@@}\nobreak\nointerlineskip
 \vskip\skipi@
 \ifdim\ht\@ne<\z@\box\@ne\else\centerline{\box\@ne}\fi}
\def\bottomfigs@{%
 \count@\@ne
 \loop
  \ifnum\count@<\flipcount@
  \nointerlineskip
  \vskip\abovebotfigskip
  \global\setbox\topins\vbox{\unvbox\topins\setboxzl@
   \unvbox\z@
   \breakisland@}%
  \printisland@
  \advance\count@\@ne
  \repeat}
\def\resetdimtopins@{%
 \global\advance\topinscount@-\flipcount@
 \global\setbox\topinsdims@\vbox
  {\unvbox\topinsdims@
   \count@\z@
   \DN@##1##2\next@{\gdef\AAlist@{##2}}%
   \loop
    \ifnum\count@<\flipcount@\setboxzl@
    \expandafter\next@\AAlist@\next@
    \advance\count@\@ne
    \repeat
   \dimen@\z@
   \count@\z@
   \setbox\tw@\vbox{}%
   \edef\nextiii@{\AAlist@}%
   \DN@##1##2\next@{\DNii@{##1}\def\nextiii@{##2}}%
   \loop
    \test@false
    \ifnum\count@<\topinscount@
    \expandafter\next@\nextiii@\next@
     \ifnum\count@<\tw@
      \test@true
     \else
      \if\nextii@ 1\test@true\fi
     \fi
    \fi
    \iftest@
     \setboxzl@
     \advance\dimen@\wdz@
     \setbox\tw@\vbox{\boxz@\unvbox\tw@}%
     \advance\count@\@ne
    \repeat
    \unvbox\tw@
    \global\dimen\topins\dimen@}}
\def\Place@#1#2{%
 \ifisland@
  \ifhmode
   \ifPar@
    \ifnum\Parcount@>5
     \Err@{Only 5 \string\place's allowed per
       \string\Par...\noexpand\endPar paragraph}%
    \else
     \expandafter\expandafter\expandafter\global\expandafter\setbox
      \csname Parfigbox\number\Parcount@\endcsname\box\islandbox@
     \global\advance\Parcount@\@ne
     \xdef\Parlist@{\Parlist@#1}%
     \vadjust{\break}%
    \fi
   \else
    \Err@{\noexpand#2allowed only in a \string\Par...\noexpand\endPar
     paragraph}%
   \fi
  \else
   #2%
  \fi
 \else
  \Err@{Who ... there's no \string\Figure, \string\Table,
   etc., here}%
 \fi
 \place@false}
\newif\ifC@
\newdimen\Cdim@
\long\def\Cplace#1{\prevanish@\place@true\island@false
 #1%
 \Place@ c\Cplace@
 \postvanish@}
\def\Cplace@{\allowbreak
 \ifnum\topinscount@>\z@\else
  \global\C@true\global\Cdim@\pagetotal\fi
 \Aplace@}
\long\def\Mplace#1{\prevanish@\place@true\island@false
 #1%
 \Place@ m\Mplace@
 \postvanish@}
\long\def\MXplace#1{\prevanish@\place@true\island@false
 #1%
 \Place@ M\MXplace@
 \postvanish@}
\newif\ifMX@
\def\MXplace@{\MX@true\Mplace@\MX@false}
\def\Mplace@{\allowbreak
 \dimen@\ht\islandbox@\advance\dimen@\dp\islandbox@
 \ifdim\pagetotal=\z@\else
  \ifdim\lastskip<\abovebotfigskip\advance\dimen@\abovebotfigskip
  \advance\dimen@-\lastskip\fi
 \fi
 \advance\dimen@\pagetotal
 \ifdim\dimen@>\pagegoal
  \Aplace@
 \else
  \nointerlineskip
  \ifdim\lastskip<\abovebotfigskip\removelastskip\vskip\abovebotfigskip\fi
  \setbox\z@\vbox{\unvbox\islandbox@
   \breakisland@}%
  \printisland@
  \ifnum\topinscount@=\z@
   \setbox\z@\vbox{\vbox to-\belowtopfigskip{}}%
   \dimen@-\skip\topins
   \ht\z@\dimen@
   \storedim@
   \advancedimtopins@
   \insert\topins{\boxz@}%
   \global\advance\topinscount@\@ne
   \xdef\AAlist@{\AAlist@0}%
  \fi
  \ifMX@
   \ifnum\topinscount@=\@ne
    \setbox\z@\vbox{\vbox to-\abovebotfigskip{}}%
    \ht\z@\z@
    \dimen@\z@
    \storedim@
    \advancedimtopins@
    \insert\topins{\boxz@}%
    \global\advance\topinscount@\@ne
    \xdef\AAlist@{\AAlist@0}%
   \fi
  \fi
  \nointerlineskip
  \vskip\belowtopfigskip
 \fi}
\expandafter\newbox\csname Parbox1\endcsname
\expandafter\newbox\csname Parbox2\endcsname
\expandafter\newbox\csname Parbox3\endcsname
\expandafter\newbox\csname Parbox4\endcsname
\expandafter\newbox\csname Parbox5\endcsname
\def\endPar{\egroup
 \count@\@ne
 {\vbadness\@M\vfuzz\maxdimen\splitmaxdepth\maxdimen\splittopskip\ht\strutbox
 \setbox\z@\vsplit\Parbox@ to\ht\Parbox@
 \loop
  \ifnum\count@<\Parcount@
  \expandafter\expandafter\expandafter\global\expandafter\setbox
   \csname Parbox\number\count@\endcsname\vsplit\Parbox@ to\ht\Parbox@
  \count@@\count@\advance\count@@\@ne
  \global\csname Parprev\number\count@@\endcsname
   \dp\csname Parbox\number\count@\endcsname
  \advance\count@\@ne
  \repeat}%
 \vskip\parskip
 \count@\@ne
 \def\nextv@##1##2\nextv@{\DN@{##1}\gdef\Parlist@{##2}}%
 \loop
  \ifnum\count@<\Parcount@
   \dimen@\csname Parprev\number\count@\endcsname
   \advance\dimen@\ht\strutbox
   \ifdim\dimen@<\baselineskip
    \advance\dimen@-\baselineskip\vskip-\dimen@
   \else
    \vskip\lineskip
   \fi
   \unvbox\csname Parbox\number\count@\endcsname
   \global\setbox\islandbox@\box\csname Parfigbox\number\count@\endcsname
   \expandafter\nextv@\Parlist@\nextv@
   \if a\next@\Aplace@\else
   \if A\next@\AAplace@\else
   \if b\next@\Bplace@\else
   \if c\next@\Cplace@\else
   \if m\next@\Mplace@\else
   \if M\next@\MXplace@\fi\fi\fi\fi\fi\fi
  \advance\count@\@ne
  \repeat
 \global\Par@false
 \ifvoid\Parbox@
  \prevdepth\csname Parprev\number\count@\endcsname
 \else
  \dimen@\csname Parprev\number\count@\endcsname\advance\dimen@\ht\strutbox
  \ifdim\dimen@<\baselineskip
    \advance\dimen@-\baselineskip\vskip-\dimen@
  \else
    \vskip\lineskip
  \fi
  \dimen@\dp\Parbox@
  \unvbox\Parbox@
  \prevdepth\dimen@
 \fi}
\def\folio{{\page@F\page@S{\page@P\page@N{\number\page@C}\page@Q}}}
\def\advancepageno{\global\advance\pageno\@ne}
\newif\ifspecialsplit@
\newbox\outbox@
\let\shipout@\shipout
\def\plainoutput{\specialsplit@false\ifvoid\topins\else\ifdim\ht\topins=\z@
 \specialsplit@true\advance\minpagesize-\skip\topins\fi\fi
 \fliptopins@
 \setbox\outbox@\vbox{\makeheadline\pagebody\makefootline}%
 {\noexpands@\let\style\relax
 \shipout@\box\outbox@}%
 \advancepageno
 \resetdimtopins@
 \ifvoid\@cclv\else\unvbox\@cclv\penalty\outputpenalty\fi
 \ifnum\outputpenalty>-\@MM\else\dosupereject\fi}
\def\pagebody{\vbox to\vsize{\boxmaxdepth\maxdepth
 \ifvoid\margin@\else
 \rlap{\kern\hsize\vbox to\z@{\kern4\p@\box\margin@\vss}}\fi
 \pagecontents}}
\newif\ifonlytop@
\def\pagecontents{%
 \onlytop@false
 \ifdim\ht\@cclv<\minpagesize\ifnum\flipcount@<\tw@\ifvoid\footins
  \onlytop@true\fi\fi\fi
 \test@false
 \ifC@
  \ifnum\flipcount@=\@ne
   \global\multiply\Cdim@\tw@
   \ifdim\Cdim@>\ht\@cclv
    \test@true
   \fi
  \fi
 \fi
 \global\C@false
 \iftest@
  \dimen@\ht\@cclv
  \advance\dimen@\skip\topins
  {\vfuzz\maxdimen\vbadness\@M
  \splitmaxdepth\maxdepth\splittopskip\topskip
  \setbox\z@\vsplit\@cclv to\dimen@
  \unvbox\z@}%
  \global\setbox\topins\vbox{\unvbox\topins
   \global\setbox\@ne\lastbox}%
  \setbox\z@\vbox{\unvbox\@ne
   \breakisland@}%
  \nointerlineskip
  \vskip\abovebotfigskip
  \printisland@
 \else
  \ifnum\flipcount@>\z@
   \global\setbox\topins\vbox{\unvbox\topins\global\setbox\@ne\lastbox}%
   \setbox\z@\vbox{\unvbox\@ne
    \breakisland@}%
   \printisland@
   \ifonlytop@\kern-\prevdepth\vfill\else\vskip\belowtopfigskip\fi
  \fi
 \fi
 \ifdim\ht\@cclv<\minpagesize
  \ifonlytop@\else\vfill\fi
 \else
  \ifspecialsplit@
   {\vfuzz\maxdimen\vbadness\@M
   \splitmaxdepth\maxdepth\splittopskip\topskip
   \dimen@ii\ht\@cclv \advance\dimen@ii\skip\topins
   \setbox\z@\vsplit\@cclv to\dimen@ii
   \unvbox\z@}%
  \else
   \unvbox\@cclv
  \fi
 \fi
 \bottomfigs@
 \ifvoid\footins\else\vskip\skip\footins\footnoterule\unvbox\footins\fi}
\newread\readdata@
\def\readthedata@#1{\expandafter
 \ifx\csname#1@D\endcsname\relax
  \immediate\openin\readdata@=#1.dat
  \ifeof\readdata@
   \Err@{No file #1.dat}%
  \else
   {\endlinechar\m@ne\gdef\Next@{}%
   \DNii@##1 ##2 ##3pt{\global\data@ht##1\global\data@dp##2%
    \global\data@wd##3pt}%
   \loop
    \ifeof\readdata@
    \else
    \read\readdata@ to\next@
    \ifx\next@\empty\else
     \edef\next@{\expandafter\nextii@\next@}%
     \expandafter\rightadd@\next@\to\Next@
    \fi
    \repeat}%
   \immediate\closein\readdata@
   \expandafter\expandafter\expandafter\global\expandafter
    \let\csname#1@D\endcsname\Next@\global\let\Next@\relax
  \fi
 \fi}
\newdimen\data@ht
\newdimen\data@dp
\newdimen\data@wd
\newif\ifgetdata@
\def\getdata@#1#2{\global\getdata@true\count@#2\relax
 {\let\\\or\xdef\Next@{\ifcase\number\count@#1\else
 \global\noexpand\getdata@false\fi}}\Next@}
\def\paste#1#2{\readthedata@{#1}%
 \getdata@{\csname#1@D\endcsname}{#2}%
 \ifgetdata@
 \dimen@\data@ht \advance\dimen@\data@dp
  \hbox{\special{dvipaste: #1 #2}%
   \lower\data@dp\vbox to\dimen@{\hbox to\data@wd{}\vfil}}%
 \else
  {\lccode`\Z=`\#\lccode`\N=`\N\lccode`\F=`\F%
   \lowercase{\Err@{No data for File [#1], Z#2}}}%
 \fi}
\newdimen\httable
\newdimen\dptable
\newdimen\wdtable
\def\measuretable#1#2{\readthedata@{#1}%
 \getdata@{\csname#1@D\endcsname}{#2}%
 \ifgetdata@
  \httable\data@ht \dptable\data@dp \wdtable\data@wd
 \else
  {\lccode`\Z=`\#\lccode`\N=`\N\lccode`\F=`\F%
  \lowercase{\Err@{No data for File [#1], Z#2}}}%
 \fi}
\def\East#1#2{\setboxz@h{$\m@th\ssize\;{#1}\;\;$}%
 \setbox\tw@\hbox{$\m@th\ssize\;{#2}\;\;$}\setbox4=\hbox{$\m@th#2$}%
 \dimen@\minaw@
 \ifdim\wdz@>\dimen@\dimen@\wdz@\fi\ifdim\wd\tw@>\dimen@\dimen@\wd\tw@\fi
 \ifdim\wd4 >\z@
  \mathrel{\mathop{\hbox to\dimen@{\rightarrowfill}}\limits^{#1}_{#2}}%
 \else
  \mathrel{\mathop{\hbox to\dimen@{\rightarrowfill}}\limits^{#1}}%
 \fi}
\def\West#1#2{\setboxz@h{$\m@th\ssize\;\;{#1}\;$}%
 \setbox\tw@\hbox{$\m@th\ssize\;\;{#2}\;$}\setbox4=\hbox{$\m@th#2$}%
 \dimen@\minaw@
 \ifdim\wdz@>\dimen@\dimen@\wdz@\fi\ifdim\wd\tw@>\dimen@\dimen@\wd\tw@\fi
 \ifdim\wd4 >\z@
  \mathrel{\mathop{\hbox to\dimen@{\leftarrowfill}}\limits^{#1}_{#2}}%
 \else
  \mathrel{\mathop{\hbox to\dimen@{\leftarrowfill}}\limits^{#1}}%
 \fi}
\font\arrow@i=lams1
\font\arrow@ii=lams2
\font\arrow@iii=lams3
\font\arrow@iv=lams4
\font\arrow@v=lams5
\newdimen\standardcgap
\standardcgap40\p@
\newdimen\hunit
\hunit\tw@\p@
\newdimen\standardrgap
\standardrgap32\p@
\newdimen\vunit
\vunit1.6\p@
\def\Cgaps#1{\RIfM@
 \standardcgap#1\standardcgap\relax\hunit#1\hunit\relax
 \else\nonmatherr@\Cgaps\fi}
\def\Rgaps#1{\RIfM@
 \standardrgap#1\standardrgap\relax\vunit#1\vunit\relax
 \else\nonmatherr@\Rgaps\fi}
\newdimen\getdim@
\def\getcgap@#1{\ifcase#1\or\getdim@\z@\else\getdim@\standardcgap\fi}
\def\getrgap@#1{\ifcase#1\getdim@\z@\else\getdim@\standardrgap\fi}
\def\cgaps{\RIfM@\expandafter\cgaps@\else\expandafter\nonmatherr@
 \expandafter\cgaps\fi}
\def\cgaps@{\ifnum\catcode`\;=\active\expandafter\cgapsA@\else
 \expandafter\cgapsO@\fi}
\def\cgapsO@#1{\toks@{\ifcase\i@\or\getdim@=\z@}%
 \gaps@@\standardcgap#1;\gaps@@\gaps@@
 \edef\next@{\the\toks@\noexpand\else\noexpand\getdim@\noexpand\standardcgap
  \noexpand\fi}%
 \toks@=\expandafter{\next@}%
 \edef\getcgap@##1{\i@##1\relax\the\toks@}\toks@{}}
{\catcode`\;=\active
 \gdef\cgapsA@#1{\toks@{\ifcase\i@\or\getdim@=\z@}%
 \gaps@@\standardcgap#1;\gaps@@\gaps@@
 \edef\next@{\the\toks@\noexpand\else\noexpand\getdim@\noexpand\standardcgap
  \noexpand\fi}%
 \toks@=\expandafter{\next@}%
 \edef\getcgap@##1{\i@##1\relax\the\toks@}\toks@{}}
}
\def\Gaps@@{\gaps@@}
\def\gaps@@#1#2;#3{\mgaps@#1#2\mgaps@
 \edef\next@{\the\toks@\noexpand\or\noexpand\getdim@
  \noexpand#1\the\mgapstoks@@}%
 \toks@\expandafter{\next@}%
 \DN@{#3}%
 \ifx\next@\Gaps@@\def\next@##1\gaps@@{}\else
  \def\next@{\gaps@@#1#3}\fi\next@}
{\catcode`\;=\active
 \gdef\rgaps#1{\RIfM@{\ifnum\catcode`\;=\active\def;{\string;}\fi
   \xdef\Next@{\noexpand\rgaps@{#1}}}%
  \Next@\edef\getrgap@##1{\i@##1\relax\the\toks@}\toks@{}\else
  \nonmatherr@\rgaps\fi}
}
\def\rgaps@#1{\toks@{\ifcase\i@\getdim@=\z@}%
 \gaps@@\standardrgap#1;\gaps@@\gaps@@
 \edef\next@{\the\toks@\noexpand\else\noexpand\getdim@\noexpand\standardrgap
  \noexpand\fi}%
 \toks@=\expandafter{\next@}}
\newbox\ZER@
\def\mgaps@#1{\let\mgapsnext@#1\FNSS@\mgaps@@}
\def\mgaps@@{\ifx\next\w\expandafter\mgaps@@@\else
 \expandafter\mgaps@@@@\fi}
\newtoks\mgapstoks@@
\def\mgaps@@@@#1\mgaps@{\getdim@\mgapsnext@\getdim@#1\getdim@
 \edef\next@{\noexpand\getdim@\the\getdim@}%
 \mgapstoks@@\expandafter{\next@}}
\def\mgaps@@@\w#1#2\mgaps@{\mgaps@@@@#2\mgaps@
 \setbox\ZER@\hbox{$\m@th\hskip15\p@\tsize@#1$}%
 \dimen@\wd\ZER@
 \ifdim\dimen@>\getdim@\getdim@\dimen@\fi
 \edef\next@{\noexpand\getdim@\the\getdim@}%
 \mgapstoks@@\expandafter{\next@}}
\def\changewidth#1#2{\setbox\ZER@{$\m@th#2}%
 \hbox to\wd\ZER@{\hss$\m@th#1$\hss}}
\atdef@({\FN@\ARROW@}
\def\ARROW@{\ifx\next)\let\next@\OPTIONS@\else
 \DN@{\csname\string @(\endcsname}\fi\next@}
\newif\ifoptions@
\def\OPTIONS@){\ifoptions@\let\next@\relax\else
 \DN@{\global\options@true\begingroup\optioncodes@}\fi\next@}
\newif\ifN@
\newif\ifE@
\newif\ifNESW@
\newif\ifH@
\newif\ifV@
\newif\ifHshort@
\expandafter\def\csname\string @(\endcsname #1,#2){%
 \ifoptions@\expandafter\endgroup\fi
 \N@false\E@false\H@false\V@false\Hshort@false
 \ifnum#1>\z@\E@true\fi
 \ifnum#1=\z@\V@true\global\tX@false\global\tY@false\global\a@false\fi
 \ifnum#2>\z@\N@true\fi
 \ifnum#2=\z@\H@true\global\tX@false\global\tY@false\global\a@false
  \ifshort@\Hshort@true\fi\fi
 \NESW@false
 \ifN@\ifE@\NESW@true\fi\else\ifE@\else\NESW@true\fi\fi
 \arrow@{#1}{#2}%
 \global\options@false
 \global\scount@\z@\global\tcount@\z@\global\arrcount@\z@
 \global\s@false\global\sxdimen@\z@\global\sydimen@\z@
 \global\tX@false\global\tXdimen@i\z@\global\tXdimen@ii\z@
 \global\tY@false\global\tYdimen@i\z@\global\tYdimen@ii\z@
 \global\a@false\global\exacount@\z@
 \global\x@false\global\xdimen@\z@
 \global\X@false\global\Xdimen@\z@
 \global\y@false\global\ydimen@\z@
 \global\Y@false\global\Ydimen@\z@
 \global\p@false\global\pdimen@\z@
 \global\label@ifalse\global\label@iifalse
 \global\dl@ifalse\global\ldimen@i\z@
 \global\dl@iifalse\global\ldimen@ii\z@
 \global\short@false\global\unshort@false}
\newif\iflabel@i
\newif\iflabel@ii
\newcount\scount@
\newcount\tcount@
\newcount\arrcount@
\newif\ifs@
\newdimen\sxdimen@
\newdimen\sydimen@
\newif\iftX@
\newdimen\tXdimen@i
\newdimen\tXdimen@ii
\newif\iftY@
\newdimen\tYdimen@i
\newdimen\tYdimen@ii
\newif\ifa@
\newcount\exacount@
\newif\ifx@
\newdimen\xdimen@
\newif\ifX@
\newdimen\Xdimen@
\newif\ify@
\newdimen\ydimen@
\newif\ifY@
\newdimen\Ydimen@
\newif\ifp@
\newdimen\pdimen@
\newif\ifdl@i
\newif\ifdl@ii
\newdimen\ldimen@i
\newdimen\ldimen@ii
\newif\ifshort@
\newif\ifunshort@
\def\zero@#1{\ifnum\scount@=\z@
 \if#1e\global\scount@\m@ne\else
 \if#1t\global\scount@\tw@\else
 \if#1h\global\scount@\thr@@\else
 \if#1'\global\scount@6 \else
 \if#1`\global\scount@7 \else
 \if#1(\global\scount@8 \else
 \if#1)\global\scount@9 \else
 \if#1s\global\scount@12 \else
 \if#1H\global\scount@13 \else
 \Err@{\Invalid@@ option \string\0}\fi\fi\fi\fi\fi\fi\fi\fi\fi
 \fi}
\def\one@#1{\ifnum\tcount@=\z@
 \if#1e\global\tcount@\m@ne\else
 \if#1h\global\tcount@\tw@\else
 \if#1t\global\tcount@\thr@@\else
 \if#1'\global\tcount@4 \else
 \if#1`\global\tcount@5 \else
 \if#1(\global\tcount@\ten@ \else
 \if#1)\global\tcount@11 \else
 \if#1s\global\tcount@12 \else
 \if#1H\global\tcount@13 \else
 \Err@{\Invalid@@ option \string\1}\fi\fi\fi\fi\fi\fi\fi\fi\fi
 \fi}
\def\a@#1{\ifnum\arrcount@=\z@
 \if#10\global\arrcount@\m@ne\else
 \if#1+\global\arrcount@\@ne\else
 \if#1-\global\arrcount@\tw@\else
 \if#1=\global\arrcount@\thr@@\else
 \Err@{\Invalid@@ option \string\a}\fi\fi\fi\fi
 \fi}
\def\ds@{\ifnum\catcode`\;=\active\expandafter\dsA@\else
 \expandafter\dsO@\fi}
\def\dsO@(#1;#2){\ds@@{#1}{#2}}
\def\ds@@#1#2{\ifs@\else
 \global\s@true
 \global\sxdimen@\hunit\global\sxdimen@#1\sxdimen@\relax
 \global\sydimen@\vunit\global\sydimen@#2\sydimen@\relax
 \fi}
\def\dtX@{\ifnum\catcode`\;=\active\expandafter\dtXA@\else
 \expandafter\dtXO@\fi}
\def\dtXO@(#1;#2){\dtX@@{#1}{#2}}
\def\dtX@@#1#2{\iftX@\else
 \global\tX@true
 \global\tXdimen@i\hunit\global\tXdimen@i#1\tXdimen@i\relax
 \global\tXdimen@ii\vunit\global\tXdimen@ii#2\tXdimen@ii\relax
 \fi}
\def\dtY@{\ifnum\catcode`\;=\active\expandafter\dtYA@\else
 \expandafter\dtYO@\fi}
\def\dtYO@(#1;#2){\dtY@@{#1}{#2}}
\def\dtY@@#1#2{\iftY@\else
 \global\tY@true
 \global\tYdimen@i\hunit\global\tYdimen@i#1\tYdimen@i\relax
 \global\tYdimen@ii\vunit\global\tYdimen@ii#2\tYdimen@ii\relax
 \fi}
{\catcode`\;=\active
 \gdef\dsA@(#1;#2){\ds@@{#1}{#2}}
 \gdef\dtXA@(#1;#2){\dtX@@{#1}{#2}}
 \gdef\dtYA@(#1;#2){\dtY@@{#1}{#2}}
}
\def\da@#1{\ifa@\else\global\a@true\global\exacount@#1\relax\fi}
\def\dx@#1{\ifx@\else
 \global\x@true
 \global\xdimen@\hunit\global\xdimen@#1\xdimen@\relax
 \fi}
\def\dX@#1{\ifX@\else
 \global\X@true
 \global\Xdimen@\hunit\global\Xdimen@#1\Xdimen@\relax
 \fi}
\def\dy@#1{\ify@\else
 \global\y@true
 \global\ydimen@\vunit\global\ydimen@#1\ydimen@\relax
 \fi}
\def\dY@#1{\ifY@\else
 \global\Y@true
 \global\Ydimen@\vunit\global\Ydimen@#1\Ydimen@\relax
 \fi}
\def\p@@#1{\ifp@\else
 \global\p@true
 \global\pdimen@\hunit\global\divide\pdimen@\tw@
 \global\pdimen@#1\pdimen@\relax
 \fi}
\def\L@#1{\iflabel@i\else
 \global\label@itrue\gdef\label@i{#1}%
 \fi}
\def\l@#1{\iflabel@ii\else
 \global\label@iitrue\gdef\label@ii{#1}%
 \fi}
\def\dL@#1{\ifdl@i\else
 \global\dl@itrue\global\ldimen@i\hunit\global\ldimen@i#1\ldimen@i\relax
 \fi}
\def\dl@#1{\ifdl@ii\else
 \global\dl@iitrue\global\ldimen@ii\hunit\global\ldimen@ii#1\ldimen@ii\relax
 \fi}
\def\s@{\ifunshort@\else\global\short@true\fi}
\def\uns@{\ifshort@\else\global\unshort@true\global\short@false\fi}
\def\optioncodes@{\let\0\zero@\let\1\one@\let\a\a@\let\ds\ds@\let\dtX\dtX@
 \let\dtY\dtY@\let\da\da@\let\dx\dx@\let\dX\dX@\let\dY\dY@\let\dy\dy@
 \let\p\p@@\let\L\L@\let\l\l@\let\dL\dL@\let\dl\dl@\let\s\s@\let\uns\uns@}
\def\slopes@{\\161\\152\\143\\134\\255\\126\\357\\238\\349\\45{10}\\56{11}%
 \\11{12}\\65{13}\\54{14}\\43{15}\\32{16}\\53{17}\\21{18}\\52{19}\\31{20}%
 \\41{21}\\51{22}\\61{23}}
\newcount\tan@i
\newcount\tan@ip
\newcount\tan@ii
\newcount\tan@iip
\newdimen\slope@i
\newdimen\slope@ip
\newdimen\slope@ii
\newdimen\slope@iip
\newcount\angcount@
\newcount\extracount@
\def\slope@{{\slope@i\secondy@\advance\slope@i-\firsty@
 \ifN@\else\multiply\slope@i\m@ne\fi
 \slope@ii\secondx@\advance\slope@ii-\firstx@
 \ifE@\else\multiply\slope@ii\m@ne\fi
 \ifdim\slope@ii<\z@
  \global\tan@i6 \global\tan@ii\@ne\global\angcount@23
 \else
  \dimen@\slope@i\multiply\dimen@6
  \ifdim\dimen@<\slope@ii
   \global\tan@i\@ne\global\tan@ii6 \global\angcount@\@ne
  \else
   \dimen@\slope@ii\multiply\dimen@6
   \ifdim\dimen@<\slope@i
    \global\tan@i6 \global\tan@ii\@ne\global\angcount@23
   \else
    \global\tan@ip\z@\global\tan@iip\@ne
    \def\\##1##2##3{\global\angcount@##3\relax
     \slope@ip\slope@i\slope@iip\slope@ii
     \multiply\slope@iip##1\relax\multiply\slope@ip##2\relax
     \ifdim\slope@iip<\slope@ip
      \global\tan@ip##1\relax\global\tan@iip##2\relax
     \else
      \global\tan@i##1\relax\global\tan@ii##2\relax
      \def\\####1####2####3{}%
     \fi}%
    \slopes@
    \slope@i\secondy@\advance\slope@i-\firsty@
    \ifN@\else\multiply\slope@i\m@ne\fi
    \multiply\slope@i\tan@ii\multiply\slope@i\tan@iip\multiply\slope@i\tw@
    \count@\tan@i\multiply\count@\tan@iip
    \extracount@\tan@ip\multiply\extracount@\tan@ii
    \advance\count@\extracount@
    \slope@ii\secondx@\advance\slope@ii-\firstx@
    \ifE@\else\multiply\slope@ii\m@ne\fi
    \multiply\slope@ii\count@
    \ifdim\slope@i<\slope@ii
     \global\tan@i\tan@ip\global\tan@ii\tan@iip
     \global\advance\angcount@\m@ne
    \fi
   \fi
  \fi
 \fi}%
}
\def\slope@a#1{{\def\\##1##2##3{\ifnum##3=#1\global\tan@i##1\relax
 \global\tan@ii##2\relax\fi}\slopes@}}
\newcount\i@
\newcount\j@
\newcount\colcount@
\newcount\Colcount@
\newcount\tcolcount@
\newdimen\rowht@
\newdimen\rowdp@
\newcount\rowcount@
\newcount\Rowcount@
\newcount\maxcolrow@
\newtoks\colwidthtoks@
\newtoks\Rowheighttoks@
\newtoks\Rowdepthtoks@
\newtoks\widthtoks@
\newtoks\Widthtoks@
\newtoks\heighttoks@
\newtoks\Heighttoks@
\newtoks\depthtoks@
\newtoks\Depthtoks@
\newif\iffirstCDcr@
\def\dotoks@i{%
 \global\widthtoks@\expandafter{\the\widthtoks@\else\getdim@\z@\fi}%
 \global\heighttoks@\expandafter{\the\heighttoks@\else\getdim@\z@\fi}%
 \global\depthtoks@\expandafter{\the\depthtoks@\else\getdim@\z@\fi}}
\def\dotoks@ii{%
 \global\widthtoks@{\ifcase\j@}%
 \global\heighttoks@{\ifcase\j@}%
 \global\depthtoks@{\ifcase\j@}}
\def\preCD@#1\endCD{\setbox\ZER@
 \vbox{%
  \def\arrow@##1##2{{}}%
  \global\rowcount@\m@ne\global\colcount@\z@\global\Colcount@\z@
  \global\firstCDcr@true\toks@{}%
  \global\widthtoks@{\ifcase\j@}%
  \global\Widthtoks@{\ifcase\i@}%
  \global\heighttoks@{\ifcase\j@}%
  \global\Heighttoks@{\ifcase\i@}%
  \global\depthtoks@{\ifcase\j@}%
  \global\Depthtoks@{\ifcase\i@}%
  \global\Rowheighttoks@{\ifcase\i@}%
  \global\Rowdepthtoks@{\ifcase\i@}%
  \Let@
  \everycr{%
   \noalign{%
    \global\advance\rowcount@\@ne
    \ifnum\colcount@<\Colcount@
    \else
     \global\Colcount@\colcount@\global\maxcolrow@\rowcount@
    \fi
    \global\colcount@\z@
    \iffirstCDcr@
     \global\firstCDcr@false
    \else
     \edef\next@{\the\Rowheighttoks@\noexpand\or\noexpand\getdim@\the\rowht@}%
      \global\Rowheighttoks@\expandafter{\next@}%
     \edef\next@{\the\Rowdepthtoks@\noexpand\or\noexpand\getdim@\the\rowdp@}%
      \global\Rowdepthtoks@\expandafter{\next@}%
     \global\rowht@\z@\global\rowdp@\z@
     \dotoks@i
     \edef\next@{\the\Widthtoks@\noexpand\or\the\widthtoks@}%
      \global\Widthtoks@\expandafter{\next@}%
     \edef\next@{\the\Heighttoks@\noexpand\or\the\heighttoks@}%
      \global\Heighttoks@\expandafter{\next@}%
     \edef\next@{\the\Depthtoks@\noexpand\or\the\depthtoks@}%
      \global\Depthtoks@\expandafter{\next@}%
     \dotoks@ii
    \fi}}%
  \tabskip\z@
  \halign{&\setbox\ZER@\hbox{\vrule\height\ten@\p@\width\z@\depth\z@     
   $\m@th\displaystyle{##}$}\copy\ZER@
   \ifdim\ht\ZER@>\rowht@\global\rowht@\ht\ZER@\fi
   \ifdim\dp\ZER@>\rowdp@\global\rowdp@\dp\ZER@\fi
   \global\advance\colcount@\@ne
   \edef\next@{\the\widthtoks@\noexpand\or\noexpand\getdim@\the\wd\ZER@}%
    \global\widthtoks@\expandafter{\next@}%
   \edef\next@{\the\heighttoks@\noexpand\or\noexpand\getdim@\the\ht\ZER@}%
    \global\heighttoks@\expandafter{\next@}%
   \edef\next@{\the\depthtoks@\noexpand\or\noexpand\getdim@\the\dp\ZER@}%
    \global\depthtoks@\expandafter{\next@}%
   \cr#1\crcr}}%
 \Rowcount@\rowcount@
 \global\Widthtoks@\expandafter{\the\Widthtoks@\fi\relax}%
 \edef\Width@##1##2{\i@##1\relax\j@##2\relax\the\Widthtoks@}%
 \global\Heighttoks@\expandafter{\the\Heighttoks@\fi\relax}%
 \edef\Height@##1##2{\i@##1\relax\j@##2\relax\the\Heighttoks@}%
 \global\Depthtoks@\expandafter{\the\Depthtoks@\fi\relax}%
 \edef\Depth@##1##2{\i@##1\relax\j@##2\relax\the\Depthtoks@}%
 \edef\next@{\the\Rowheighttoks@\noexpand\fi\relax}%
 \global\Rowheighttoks@\expandafter{\next@}%
 \edef\Rowheight@##1{\i@##1\relax\the\Rowheighttoks@}%
 \edef\next@{\the\Rowdepthtoks@\noexpand\fi\relax}%
 \global\Rowdepthtoks@\expandafter{\next@}%
 \edef\Rowdepth@##1{\i@##1\relax\the\Rowdepthtoks@}%
 \global\colwidthtoks@{\fi}%
 \setbox\ZER@\vbox{%
  \unvbox\ZER@
  \count@\rowcount@
  \loop
   \unskip\unpenalty
   \setbox\ZER@\lastbox
   \ifnum\count@>\maxcolrow@\advance\count@\m@ne
   \repeat
  \hbox{%
   \unhbox\ZER@
   \count@\z@
   \loop
    \unskip
    \setbox\ZER@\lastbox
    \edef\next@{\noexpand\or\noexpand\getdim@\the\wd\ZER@\the\colwidthtoks@}%
     \global\colwidthtoks@\expandafter{\next@}%
    \advance\count@\@ne
    \ifnum\count@<\Colcount@
    \repeat}}%
 \edef\next@{\noexpand\ifcase\noexpand\i@\the\colwidthtoks@}%
  \global\colwidthtoks@\expandafter{\next@}%
 \edef\Colwidth@##1{\i@##1\relax\the\colwidthtoks@}%
 \global\colwidthtoks@{}\global\Rowheighttoks@{}\global\Rowdepthtoks@{}%
 \global\widthtoks@{}\global\Widthtoks@{}\global\heighttoks@{}%
 \global\Heighttoks@{}\global\depthtoks@{}\global\Depthtoks@{}%
}
\newcount\xoff@
\newcount\yoff@
\newcount\endcount@
\newcount\rcount@
\newdimen\firstx@
\newdimen\firsty@
\newdimen\secondx@
\newdimen\secondy@
\newdimen\tocenter@
\newdimen\charht@
\newdimen\charwd@
\def\outside@{\Err@{This arrow points outside the \string\CD}}
\newif\ifsvertex@
\newif\iftvertex@
\def\arrow@#1#2{\global\xoff@#1\relax\global\yoff@#2\relax
 \count@\rowcount@\advance\count@-\yoff@
 \ifnum\count@<\@ne\outside@\else\ifnum\count@>\Rowcount@\outside@\fi\fi
 \count@\colcount@\advance\count@\xoff@
 \ifnum\count@<\@ne\outside@\else\ifnum\count@>\Colcount@\outside@\fi\fi
 \tcolcount@\colcount@\advance\tcolcount@\xoff@
 \Width@\rowcount@\colcount@\divide\getdim@\tw@\tocenter@-\getdim@
 \ifdim\getdim@=\z@
  \firstx@\z@\firsty@\mathaxis@\svertex@true
 \else
  \svertex@false
  \ifHshort@
   \Colwidth@\colcount@\divide\getdim@\tw@
   \ifE@ \firstx@\getdim@ \else \firstx@-\getdim@ \fi
  \else
   \ifE@ \firstx@\getdim@ \else \firstx@-\getdim@ \fi
  \fi
  \ifE@
   \ifH@ \advance\firstx@\thr@@\p@ \else \advance\firstx@-\thr@@\p@ \fi  
  \else
   \ifH@ \advance\firstx@-\thr@@\p@ \else \advance\firstx@\thr@@\p@ \fi  
  \fi
  \ifN@
   \Height@\rowcount@\colcount@ \firsty@\getdim@                         
   \ifV@ \advance\firsty@\thr@@\p@ \fi                                   
  \else
   \ifV@
    \Depth@\rowcount@\colcount@ \firsty@-\getdim@                        
    \advance\firsty@-\thr@@\p@                                           
   \else
    \firsty@\z@                                                          
   \fi
  \fi
 \fi
 \ifV@
 \else
  \Colwidth@\colcount@\divide\getdim@\tw@
  \ifE@\secondx@\getdim@\else\secondx@-\getdim@\fi
  \ifE@\else\getcgap@\colcount@\advance\secondx@-\getdim@\fi
  \endcount@\colcount@\advance\endcount@\xoff@
  \count@\colcount@
  \ifE@
   \advance\count@\@ne
   \loop
    \ifnum\count@<\endcount@
    \Colwidth@\count@\advance\secondx@\getdim@
    \getcgap@\count@\advance\secondx@\getdim@
    \advance\count@\@ne
    \repeat
  \else
   \advance\count@\m@ne
   \loop
    \ifnum\count@>\endcount@
    \Colwidth@\count@\advance\secondx@-\getdim@
    \getcgap@\count@\advance\secondx@-\getdim@
    \advance\count@\m@ne
    \repeat
  \fi
  \Colwidth@\count@\divide\getdim@\tw@
  \ifHshort@
  \else
   \ifE@\advance\secondx@\getdim@\else\advance\secondx@-\getdim@\fi
  \fi
  \ifE@\getcgap@\count@\advance\secondx@\getdim@\fi
  \rcount@\rowcount@\advance\rcount@-\yoff@
  \Width@\rcount@\count@\divide\getdim@\tw@
  \tvertex@false
  \ifH@\ifdim\getdim@=\z@\tvertex@true\Hshort@false\fi\fi
  \ifHshort@
  \else
   \ifE@\advance\secondx@-\getdim@\else\advance\secondx@\getdim@\fi
  \fi
  \iftvertex@
   \advance\secondx@.4\p@
  \else
   \ifE@\advance\secondx@-\thr@@\p@\else\advance\secondx@\thr@@\p@\fi    
  \fi
 \fi
 \ifH@
 \else
  \ifN@
   \Rowheight@\rowcount@\secondy@\getdim@
  \else
   \Rowdepth@\rowcount@\secondy@-\getdim@
   \getrgap@\rowcount@\advance\secondy@-\getdim@
  \fi
  \endcount@\rowcount@\advance\endcount@-\yoff@
  \count@\rowcount@
  \ifN@
   \advance\count@\m@ne
   \loop
    \ifnum\count@>\endcount@
    \Rowheight@\count@\advance\secondy@\getdim@
    \Rowdepth@\count@\advance\secondy@\getdim@
    \getrgap@\count@\advance\secondy@\getdim@
    \advance\count@\m@ne
    \repeat
  \else
   \advance\count@\@ne
   \loop
    \ifnum\count@<\endcount@
    \Rowheight@\count@\advance\secondy@-\getdim@
    \Rowdepth@\count@\advance\secondy@-\getdim@
    \getrgap@\count@\advance\secondy@-\getdim@
    \advance\count@\@ne
    \repeat
  \fi
  \tvertex@false
  \ifV@\Width@\count@\colcount@\ifdim\getdim@=\z@\tvertex@true\fi\fi
  \ifN@
   \getrgap@\count@\advance\secondy@\getdim@
   \Rowdepth@\count@\advance\secondy@\getdim@
   \iftvertex@
    \advance\secondy@\mathaxis@
   \else
    \Depth@\count@\tcolcount@\advance\secondy@-\getdim@
    \advance\secondy@-\thr@@\p@                                          
   \fi
  \else
   \Rowheight@\count@\advance\secondy@-\getdim@
   \iftvertex@
    \advance\secondy@\mathaxis@
   \else
    \Height@\count@\tcolcount@\advance\secondy@\getdim@
    \advance\secondy@\thr@@\p@                                           
   \fi
  \fi
 \fi
 \ifV@\else\advance\firstx@\sxdimen@\fi
 \ifH@\else\advance\firsty@\sydimen@\fi
 \iftX@
  \advance\secondy@\tXdimen@ii
  \advance\secondx@\tXdimen@i
  \slope@
 \else
  \iftY@
   \advance\secondy@\tYdimen@ii
   \advance\secondx@\tYdimen@i
   \slope@
   \secondy@\secondx@\advance\secondy@-\firstx@
   \ifNESW@\else\multiply\secondy@\m@ne\fi
   \multiply\secondy@\tan@i\divide\secondy@\tan@ii\advance\secondy@\firsty@
  \else
   \ifa@
    \slope@
    \ifNESW@\global\advance\angcount@\exacount@\else
     \global\advance\angcount@-\exacount@\fi
    \ifnum\angcount@>23 \global\angcount@23 \fi
    \ifnum\angcount@<\@ne\global\angcount@\@ne\fi
    \slope@a\angcount@
    \ifY@
     \advance\secondy@\Ydimen@
    \else
     \ifX@
      \advance\secondx@\Xdimen@
      \dimen@\secondx@\advance\dimen@-\firstx@
      \ifNESW@\else\multiply\dimen@\m@ne\fi
      \multiply\dimen@\tan@i\divide\dimen@\tan@ii
      \advance\dimen@\firsty@\secondy@\dimen@
     \fi
    \fi
   \else
    \ifH@\else\ifV@\else\slope@\fi\fi
   \fi
  \fi
 \fi
 \ifH@\else\ifV@\else\ifsvertex@\else
  \dimen@6\p@\multiply\dimen@\tan@ii
  \count@\tan@i\advance\count@\tan@ii\divide\dimen@\count@
  \ifE@\advance\firstx@\dimen@\else\advance\firstx@-\dimen@\fi
  \multiply\dimen@\tan@i\divide\dimen@\tan@ii
  \ifN@\advance\firsty@\dimen@\else\advance\firsty@-\dimen@\fi
 \fi\fi\fi
 \ifp@
  \ifH@\else\ifV@\else
   \getcos@\pdimen@\advance\firsty@\dimen@\advance\secondy@\dimen@
   \ifNESW@\advance\firstx@-\dimen@ii\else\advance\firstx@\dimen@ii\fi
  \fi\fi
 \fi
 \ifH@\else\ifV@\else
  \ifnum\tan@i>\tan@ii
   \charht@\ten@\p@\charwd@\ten@\p@
   \multiply\charwd@\tan@ii\divide\charwd@\tan@i
  \else
   \charwd@\ten@\p@\charht@\ten@\p@
   \divide\charht@\tan@ii\multiply\charht@\tan@i
  \fi
  \ifnum\tcount@=\thr@@
   \ifN@\advance\secondy@-.3\charht@\else\advance\secondy@.3\charht@\fi
  \fi
  \ifnum\scount@=\tw@
   \ifE@\advance\firstx@.3\charht@\else\advance\firstx@-.3\charht@\fi
  \fi
  \ifnum\tcount@=12
   \ifN@\advance\secondy@-\charht@\else\advance\secondy@\charht@\fi
  \fi
  \iftY@
  \else
   \ifa@
    \ifX@
    \else
     \secondx@\secondy@\advance\secondx@-\firsty@
     \ifNESW@\else\multiply\secondx@\m@ne\fi
     \multiply\secondx@\tan@ii\divide\secondx@\tan@i
     \advance\secondx@\firstx@
    \fi
   \fi
  \fi
 \fi\fi
 \ifH@\harrow@\else\ifV@\varrow@\else\arrow@@\fi\fi}
\newdimen\mathaxis@
\mathaxis@90\p@\divide\mathaxis@36
\def\harrow@b{\ifE@\hskip\tocenter@\hskip\firstx@\fi}
\def\harrow@bb{\ifE@\hskip\xdimen@\else\hskip\Xdimen@\fi}
\def\harrow@e{\ifE@\else\hskip-\firstx@\hskip-\tocenter@\fi}
\def\harrow@ee{\ifE@\hskip-\Xdimen@\else\hskip-\xdimen@\fi}
\def\harrow@{\dimen@\secondx@\advance\dimen@-\firstx@
 \ifE@\let\next@\rlap\else\multiply\dimen@\m@ne\let\next@\llap\fi
 \next@{%
  \harrow@b
  \smash{\raise\pdimen@\hbox to\dimen@
   {\harrow@bb\arrow@ii
    \ifnum\arrcount@=\m@ne\else\ifnum\arrcount@=\thr@@\else
     \ifE@
      \ifnum\scount@=\m@ne
      \else
       \ifcase\scount@\or\or\char118 \or\char117 \or\or\or\char119 \or
       \char120 \or\char121 \or\char122 \or\or\or\arrow@i\char125 \or
       \char117 \hskip\thr@@\p@\char117 \hskip-\thr@@\p@\fi
      \fi
     \else
      \ifnum\tcount@=\m@ne
      \else
       \ifcase\tcount@\char117 \or\or\char117 \or\char118 \or\char119 \or
       \char120 \or\or\or\or\or\char121 \or\char122 \or\arrow@i\char125
       \or\char117 \hskip\thr@@\p@\char117 \hskip-\thr@@\p@\fi
      \fi
     \fi
    \fi\fi
    \dimen@\mathaxis@\advance\dimen@.2\p@
    \dimen@ii\mathaxis@\advance\dimen@ii-.2\p@
    \ifnum\arrcount@=\m@ne
     \let\leads@\null
    \else
     \ifcase\arrcount@
      \def\leads@{\hrule\height\dimen@\depth-\dimen@ii}\or
      \def\leads@{\hrule\height\dimen@\depth-\dimen@ii}\or
      \def\leads@{\hbox to\ten@\p@{%
       \leaders\hrule\height\dimen@\depth-\dimen@ii\hfil
       \hfil
      \leaders\hrule\height\dimen@\depth-\dimen@ii\hskip\z@ plus2fil\relax
       \hfil
       \leaders\hrule\height\dimen@\depth-\dimen@ii\hfil}}\or
     \def\leads@{\hbox{\hbox to\ten@\p@{\dimen@\mathaxis@\advance\dimen@1.2\p@
       \dimen@ii\dimen@\advance\dimen@ii-.4\p@
       \leaders\hrule\height\dimen@\depth-\dimen@ii\hfil}%
       \kern-\ten@\p@
       \hbox to\ten@\p@{\dimen@\mathaxis@\advance\dimen@-1.2\p@
       \dimen@ii\dimen@\advance\dimen@ii-.4\p@
       \leaders\hrule\height\dimen@\depth-\dimen@ii\hfil}}}\fi
    \fi
    \cleaders\leads@\hfil
    \ifnum\arrcount@=\m@ne\else\ifnum\arrcount@=\thr@@\else
     \arrow@i
     \ifE@
      \ifnum\tcount@=\m@ne
      \else
       \ifcase\tcount@\char119 \or\or\char119 \or\char120 \or\char121 \or
       \char122 \or \or\or\or\or\char123 \or\char124 \or
       \char125 \or\char119 \hskip-\thr@@\p@\char119 \hskip\thr@@\p@\fi
      \fi
     \else
      \ifcase\scount@\or\or\char120 \or\char119 \or\or\or\char121 \or\char122
      \or\char123 \or\char124 \or\or\or\char125 \or
      \char119 \hskip-\thr@@\p@\char119 \hskip\thr@@\p@\fi
     \fi
    \fi\fi
    \harrow@ee}}%
  \harrow@e}%
 \iflabel@i
  \dimen@ii\z@\setbox\ZER@\hbox{$\m@th\tsize@@\label@i$}%
  \ifnum\arrcount@=\m@ne
  \else
   \advance\dimen@ii\mathaxis@
   \advance\dimen@ii\dp\ZER@\advance\dimen@ii\tw@\p@
   \ifnum\arrcount@=\thr@@\advance\dimen@ii\tw@\p@\fi
  \fi
  \advance\dimen@ii\pdimen@
  \next@{\harrow@b\smash{\raise\dimen@ii\hbox to\dimen@
   {\harrow@bb\hskip\tw@\ldimen@i\hfil\box\ZER@\hfil\harrow@ee}}\harrow@e}%
 \fi
 \iflabel@ii
  \ifnum\arrcount@=\m@ne
  \else
   \setbox\ZER@\hbox{$\m@th\tsize@\label@ii$}%
   \dimen@ii-\ht\ZER@\advance\dimen@ii-\tw@\p@
   \ifnum\arrcount@=\thr@@\advance\dimen@ii-\tw@\p@\fi
   \advance\dimen@ii\mathaxis@\advance\dimen@ii\pdimen@
   \next@{\harrow@b\smash{\raise\dimen@ii\hbox to\dimen@
    {\harrow@bb\hskip\tw@\ldimen@ii\hfil\box\ZER@\hfil\harrow@ee}}\harrow@e}%
  \fi
 \fi}
\let\tsize@\tsize
\def\tsizeCDlabels{\let\tsize@\tsize}
\def\ssizeCDlabels{\let\tsize@\ssize}
\def\tsize@@{\ifnum\arrcount@=\m@ne\else\tsize@\fi}
\def\varrow@{\dimen@\secondy@\advance\dimen@-\firsty@
 \ifN@\else\multiply\dimen@\m@ne\fi
 \setbox\ZER@\vbox to\dimen@
  {\ifN@\vskip-\Ydimen@\else\vskip\ydimen@\fi
   \ifnum\arrcount@=\m@ne\else\ifnum\arrcount@=\thr@@\else
    \hbox{\arrow@iii
     \ifN@
      \ifnum\tcount@=\m@ne
      \else
       \ifcase\tcount@\char117 \or\or\char117 \or\char118 \or\char119 \or
       \char120 \or\or\or\or\or\char121 \or\char122 \or\char123 \or
       \vbox{\hbox{\char117}\nointerlineskip\vskip\thr@@\p@
       \hbox{\char117}\vskip-\thr@@\p@}\fi
      \fi
     \else
      \ifcase\scount@\or\or\char118 \or\char117 \or\or\or\char119 \or
      \char120 \or\char121 \or\char122 \or\or\or\char123 \or
      \vbox{\hbox{\char117}\nointerlineskip\vskip\thr@@\p@
      \hbox{\char117}\vskip-\thr@@\p@}\fi
     \fi}%
    \nointerlineskip
   \fi\fi
   \ifnum\arrcount@=\m@ne
    \let\leads@\null
   \else
    \ifcase\arrcount@\let\leads@\vrule\or\let\leads@\vrule\or
    \def\leads@{\vbox to\ten@\p@{%
     \hrule\height1.67\p@\depth\z@\width.4\p@
     \vfil
     \hrule\height3.33\p@\depth\z@\width.4\p@
     \vfil
     \hrule\height1.67\p@\depth\z@\width.4\p@}}\or
    \def\leads@{\hbox{\vrule\height\p@\hskip\tw@\p@\vrule}}\fi
   \fi
  \cleaders\leads@\vfill\nointerlineskip
   \ifnum\arrcount@=\m@ne\else\ifnum\arrcount@=\thr@@\else
    \hbox{\arrow@iv
     \ifN@
      \ifcase\scount@\or\or\char118 \or\char117 \or\or\or\char119 \or
      \char120 \or\char121 \or\char122 \or\or\or\arrow@iii\char123 \or
      \vbox{\hbox{\char117}\nointerlineskip\vskip-\thr@@\p@
      \hbox{\char117}\vskip\thr@@\p@}\fi
     \else
      \ifnum\tcount@=\m@ne
      \else
       \ifcase\tcount@\char117 \or\or\char117 \or\char118 \or\char119 \or
       \char120 \or\or\or\or\or\char121 \or\char122 \or\arrow@iii\char123 \or
       \vbox{\hbox{\char117}\nointerlineskip\vskip-\thr@@\p@
       \hbox{\char117}\vskip\thr@@\p@}\fi
      \fi
     \fi}%
   \fi\fi
   \ifN@\vskip\ydimen@\else\vskip-\Ydimen@\fi}%
 \ifN@
  \dimen@ii\firsty@
 \else
  \dimen@ii-\firsty@\advance\dimen@ii\ht\ZER@\multiply\dimen@ii\m@ne
 \fi
 \rlap{\smash{\hskip\tocenter@\hskip\pdimen@\raise\dimen@ii\box\ZER@}}%
 \iflabel@i
  \setbox\ZER@\vbox to\dimen@{\vfil
   \hbox{$\m@th\tsize@@\label@i$}\vskip\tw@\ldimen@i\vfil}%
  \rlap{\smash{\hskip\tocenter@\hskip\pdimen@
  \ifnum\arrcount@=\m@ne\let\next@\relax\else\let\next@\llap\fi
  \next@{\raise\dimen@ii\hbox{\ifnum\arrcount@=\m@ne\hskip-.5\wd\ZER@\fi
   \box\ZER@\ifnum\arrcount@=\m@ne\else\hskip\tw@\p@\fi}}}}%
 \fi
 \iflabel@ii
  \ifnum\arrcount@=\m@ne
  \else
   \setbox\ZER@\vbox to\dimen@{\vfil
    \hbox{$\m@th\tsize@\label@ii$}\vskip\tw@\ldimen@ii\vfil}%
   \rlap{\smash{\hskip\tocenter@\hskip\pdimen@
   \rlap{\raise\dimen@ii\hbox{\ifnum\arrcount@=\thr@@\hskip4.5\p@\else
    \hskip2.5\p@\fi\box\ZER@}}}}%
  \fi
 \fi
}
\newdimen\goal@
\newdimen\shifted@
\newcount\Tcount@
\newcount\Scount@
\newbox\shaft@
\newcount\slcount@
\def\getcos@#1{%
 \ifnum\tan@i<\tan@ii
  \dimen@#1%
  \ifnum\slcount@<8 \count@9 \else \ifnum\slcount@<12 \count@8 \else
   \count@7 \fi\fi
  \multiply\dimen@\count@\divide\dimen@\ten@
  \dimen@ii\dimen@\multiply\dimen@ii\tan@i\divide\dimen@ii\tan@ii
 \else
  \dimen@ii#1%
  \count@-\slcount@\advance\count@24
  \ifnum\count@<8 \count@9 \else \ifnum\count@<12 \count@8
   \else\count@7 \fi\fi
  \multiply\dimen@ii\count@\divide\dimen@ii\ten@
  \dimen@\dimen@ii\multiply\dimen@\tan@ii\divide\dimen@\tan@i
 \fi}
\newdimen\adjust@
\def\Nnext@{\ifN@\let\next@\raise\else\let\next@\lower\fi}
\def\arrow@@{\slcount@\angcount@
 \ifNESW@
  \ifnum\angcount@<\ten@
   \let\arrowfont@\arrow@i\global\advance\angcount@\m@ne
   \global\multiply\angcount@13
  \else
   \ifnum\angcount@<19
    \let\arrowfont@\arrow@ii\global\advance\angcount@-\ten@
    \global\multiply\angcount@13
   \else
    \let\arrowfont@\arrow@iii\global\advance\angcount@-19
    \global\multiply\angcount@13
  \fi\fi
  \Tcount@\angcount@
 \else
  \ifnum\angcount@<5
   \let\arrowfont@\arrow@iii\global\advance\angcount@\m@ne
   \global\multiply\angcount@13 \global\advance\angcount@65
  \else
   \ifnum\angcount@<14
    \let\arrowfont@\arrow@iv\global\advance\angcount@-5
    \global\multiply\angcount@13
   \else
    \ifnum\angcount@<23
     \let\arrowfont@\arrow@v\global\advance\angcount@-14
     \global\multiply\angcount@13
    \else
     \let\arrowfont@\arrow@i\global\angcount@117
  \fi\fi\fi
  \ifnum\angcount@=117 \Tcount@115 \else\Tcount@\angcount@\fi
 \fi
 \Scount@\Tcount@
 \ifE@
  \ifnum\tcount@=\z@\advance\Tcount@\tw@\else\ifnum\tcount@=13
   \advance\Tcount@\tw@\else\advance\Tcount@\tcount@\fi\fi
  \ifnum\scount@=\z@\else\ifnum\scount@=13 \advance\Scount@\thr@@\else
   \advance\Scount@\scount@\fi\fi
 \else
  \ifcase\tcount@\advance\Tcount@\thr@@\or\or\advance\Tcount@\thr@@\or
  \advance\Tcount@\tw@\or\advance\Tcount@6 \or\advance\Tcount@7
  \or\or\or\or\or\advance\Tcount@8 \or\advance\Tcount@9 \or
  \advance\Tcount@12 \or\advance\Tcount@\thr@@\fi
  \ifcase\scount@\or\or\advance\Scount@\thr@@\or\advance\Scount@\tw@\or
  \or\or\advance\Scount@4 \or\advance\Scount@5 \or\advance\Scount@\ten@
  \or\advance\Scount@11 \or\or\or\advance\Scount@12 \or\advance
  \Scount@\tw@\fi
 \fi
 \ifcase\arrcount@\or\or\global\advance\angcount@\@ne\else\fi
 \ifN@\shifted@\firsty@\else\shifted@-\firsty@\fi
 \ifE@\else\advance\shifted@\charht@\fi
 \goal@\secondy@\advance\goal@-\firsty@
 \ifN@\else\multiply\goal@\m@ne\fi
 \setbox\shaft@\hbox{\arrowfont@\char\angcount@}%
 \ifnum\arrcount@=\thr@@
  \getcos@{1.5\p@}%
  \setbox\shaft@\hbox to\wd\shaft@{\arrowfont@
   \rlap{\hskip\dimen@ii
    \smash{\ifNESW@\let\next@\lower\else\let\next@\raise\fi
     \next@\dimen@\hbox{\arrowfont@\char\angcount@}}}%
   \rlap{\hskip-\dimen@ii
    \smash{\ifNESW@\let\next@\raise\else\let\next@\lower\fi
      \next@\dimen@\hbox{\arrowfont@\char\angcount@}}}\hfil}%
 \fi
 \rlap{\smash{\hskip\tocenter@\hskip\firstx@
  \ifnum\arrcount@=\m@ne
  \else
   \ifnum\arrcount@=\thr@@
   \else
    \ifnum\scount@=\m@ne
    \else
     \ifnum\scount@=\z@
     \else
      \setbox\ZER@\hbox{\ifnum\angcount@=117 \arrow@v\else\arrowfont@\fi
       \char\Scount@}%
      \ifNESW@
       \ifnum\scount@=\tw@
        \dimen@\shifted@\advance\dimen@-\charht@
        \ifN@\hskip-\wd\ZER@\fi
        \Nnext@
        \next@\dimen@\copy\ZER@
        \ifN@\else\hskip-\wd\ZER@\fi
       \else
        \Nnext@
        \ifN@\else\hskip-\wd\ZER@\fi
        \next@\shifted@\copy\ZER@
        \ifN@\hskip-\wd\ZER@\fi
       \fi
       \ifnum\scount@=12
        \advance\shifted@\charht@\advance\goal@-\charht@
        \ifN@\hskip\wd\ZER@\else\hskip-\wd\ZER@\fi
       \fi
       \ifnum\scount@=13
        \getcos@{\thr@@\p@}%
        \ifN@\hskip\dimen@\else\hskip-\wd\ZER@\hskip-\dimen@\fi
        \adjust@\shifted@\advance\adjust@\dimen@ii
        \Nnext@
        \next@\adjust@\copy\ZER@
        \ifN@\hskip-\dimen@\hskip-\wd\ZER@\else\hskip\dimen@\fi
       \fi
      \else
       \ifN@\hskip-\wd\ZER@\fi
       \ifnum\scount@=\tw@
        \ifN@\hskip\wd\ZER@\else\hskip-\wd\ZER@\fi
        \dimen@\shifted@\advance\dimen@-\charht@
        \Nnext@
        \next@\dimen@\copy\ZER@
        \ifN@\hskip-\wd\ZER@\fi
       \else
        \Nnext@
        \next@\shifted@\copy\ZER@
        \ifN@\else\hskip-\wd\ZER@\fi
       \fi
       \ifnum\scount@=12
        \advance\shifted@\charht@\advance\goal@-\charht@
        \ifN@\hskip-\wd\ZER@\else\hskip\wd\ZER@\fi
       \fi
       \ifnum\scount@=13
        \getcos@{\thr@@\p@}%
        \ifN@\hskip-\wd\ZER@\hskip-\dimen@\else\hskip\dimen@\fi
        \adjust@\shifted@\advance\adjust@\dimen@ii
        \Nnext@
        \next@\adjust@\copy\ZER@
        \ifN@\hskip\dimen@\else\hskip-\dimen@\hskip-\wd\ZER@\fi
       \fi	
      \fi
  \fi\fi\fi\fi
  \ifnum\arrcount@=\m@ne
  \else
   \loop
    \ifdim\goal@>\charht@
    \ifE@\else\hskip-\charwd@\fi
    \Nnext@
    \next@\shifted@\copy\shaft@
    \ifE@\else\hskip-\charwd@\fi
    \advance\shifted@\charht@\advance\goal@-\charht@
    \repeat
   \ifdim\goal@>\z@
    \dimen@\charht@\advance\dimen@-\goal@
    \divide\dimen@\tan@i\multiply\dimen@\tan@ii
    \ifE@\hskip-\dimen@\else\hskip-\charwd@\hskip\dimen@\fi
    \adjust@\shifted@\advance\adjust@-\charht@\advance\adjust@\goal@
    \Nnext@
    \next@\adjust@\copy\shaft@
    \ifE@\else\hskip-\charwd@\fi
   \else
    \adjust@\shifted@\advance\adjust@-\charht@
   \fi
  \fi
  \ifnum\arrcount@=\m@ne
  \else
   \ifnum\arrcount@=\thr@@
   \else
    \ifnum\tcount@=\m@ne
    \else
     \setbox\ZER@
      \hbox{\ifnum\angcount@=117 \arrow@v\else\arrowfont@\fi\char\Tcount@}%
     \ifnum\tcount@=\thr@@
      \advance\adjust@\charht@
      \ifE@\else\ifN@\hskip-\charwd@\else\hskip-\wd\ZER@\fi\fi
     \else
      \ifnum\tcount@=12
       \advance\adjust@\charht@
       \ifE@\else\ifN@\hskip-\charwd@\else\hskip-\wd\ZER@\fi\fi
      \else
       \ifE@\hskip-\wd\ZER@\fi
     \fi\fi
     \Nnext@
     \next@\adjust@\copy\ZER@
     \ifnum\tcount@=13
      \hskip-\wd\ZER@
      \getcos@{\thr@@\p@}%
      \ifE@\hskip-\dimen@\else\hskip\dimen@\fi
      \advance\adjust@-\dimen@ii
      \Nnext@
      \next@\adjust@\box\ZER@
     \fi
  \fi\fi\fi}}%
 \iflabel@i
  \rlap{\hskip\tocenter@
  \dimen@\firstx@\advance\dimen@\secondx@\divide\dimen@\tw@
  \advance\dimen@\ldimen@i
  \dimen@ii\firsty@\advance\dimen@ii\secondy@\divide\dimen@ii\tw@
  \global\multiply\ldimen@i\tan@i\global\divide\ldimen@i\tan@ii
  \ifNESW@\advance\dimen@ii\ldimen@i\else\advance\dimen@ii-\ldimen@i\fi
  \setbox\ZER@\hbox{\ifNESW@\else\ifnum\arrcount@=\thr@@\hskip4\p@\else
   \hskip\tw@\p@\fi\fi
   $\m@th\tsize@@\label@i$\ifNESW@\ifnum\arrcount@=\thr@@\hskip4\p@\else
   \hskip\tw@\p@\fi\fi}%
  \ifnum\arrcount@=\m@ne
   \ifNESW@\advance\dimen@.5\wd\ZER@\advance\dimen@\p@\else
    \advance\dimen@-.5\wd\ZER@\advance\dimen@-\p@\fi
   \advance\dimen@ii-.5\ht\ZER@
  \else
   \advance\dimen@ii\dp\ZER@
   \ifnum\slcount@<6 \advance\dimen@ii\tw@\p@\fi
  \fi
  \hskip\dimen@
  \ifNESW@\let\next@\llap\else\let\next@\rlap\fi
  \next@{\smash{\raise\dimen@ii\box\ZER@}}}%
 \fi
 \iflabel@ii
  \ifnum\arrcount@=\m@ne
  \else
   \rlap{\hskip\tocenter@
   \dimen@\firstx@\advance\dimen@\secondx@\divide\dimen@\tw@
   \ifNESW@\advance\dimen@\ldimen@ii\else\advance\dimen@-\ldimen@ii\fi
   \dimen@ii\firsty@\advance\dimen@ii\secondy@\divide\dimen@ii\tw@
   \global\multiply\ldimen@ii\tan@i\global\divide\ldimen@ii\tan@ii
   \advance\dimen@ii\ldimen@ii
   \setbox\ZER@\hbox{\ifNESW@\ifnum\arrcount@=\thr@@\hskip4\p@\else
    \hskip\tw@\p@\fi\fi
    $\m@th\tsize@\label@ii$\ifNESW@\else\ifnum\arrcount@=\thr@@\hskip4\p@
    \else\hskip\tw@\p@\fi\fi}%
   \advance\dimen@ii-\ht\ZER@
   \ifnum\slcount@<9 \advance\dimen@ii-\thr@@\p@\fi
   \ifNESW@\let\next@\rlap\else\let\next@\llap\fi
   \hskip\dimen@\next@{\smash{\raise\dimen@ii\box\ZER@}}}%
  \fi
 \fi
}
\def\outCD@#1{\def#1{\Err@{\noexpand#1must not be used within \string\CD}}}
\newskip\preCDskip@
\newskip\postCDskip@
\preCDskip@\z@
\postCDskip@\z@
\def\preCDspace#1{\RIfMIfI@
 \onlydmatherr@\preCDspace\else\advance\preCDskip@#1\relax\fi\else
 \onlydmatherr@\preCDspace\fi}
\def\postCDspace#1{\RIfMIfI@
 \onlydmatherr@\postCDspace\else\advance\postCDskip@#1\relax\fi\else
 \onlydmatherr@\postCDspace\fi}
\def\predisplayspace#1{\RIfMIfI@
 \onlydmatherr@\predisplayspace\else
 \advance\abovedisplayskip#1\relax
 \advance\abovedisplayshortskip#1\relax\fi
 \else\onlydmatherr@\preCDspace\fi}
\def\postdisplayspace#1{\RIfMIfI@
 \onlydmatherr@\postdisplayspace\else
 \advance\belowdisplayskip#1\relax
 \advance\belowdisplayshortskip#1\relax\fi
 \else\onlydmatherr@\postdisplayspace\fi}
\def\PreCDSpace#1{\global\preCDskip@#1\relax}
\def\PostCDSpace#1{\global\postCDskip@#1\relax}
\def\CD#1\endCD{%
 \outCD@\cgaps\outCD@\rgaps\outCD@\Cgaps\outCD@\Rgaps
 \preCD@#1\endCD
 \advance\abovedisplayskip\preCDskip@
 \advance\abovedisplayshortskip\preCDskip@
 \advance\belowdisplayskip\postCDskip@
 \advance\belowdisplayshortskip\postCDskip@
 \vcenter{\offinterlineskip
  \vskip\preCDskip@\Let@\global\colcount@\@ne\global\rowcount@\z@
  \everycr{%
   \noalign{%
    \ifnum\rowcount@=\Rowcount@
    \else
     \getrgap@\rowcount@\vskip\getdim@
     \global\advance\rowcount@\@ne\global\colcount@\@ne
    \fi}}%
  \tabskip\z@
  \halign{&\global\xoff@\z@\global\yoff@\z@
   \getcgap@\colcount@\hskip\getdim@
   \hfil\vrule\height\ten@\p@\width\z@\depth\z@
   $\m@th\displaystyle{##}$\hfil
   \global\advance\colcount@\@ne\cr
   #1\crcr}\vskip\postCDskip@}%
 \preCDskip@\z@\postCDskip@\z@
 \def\getcgap@##1{\ifcase##1\or\getdim@\z@\else\getdim@\standardcgap\fi}%
 \def\getrgap@##1{\ifcase##1\getdim@\z@\else\getdim@\standardrgap\fi}%
 \let\Width@\relax\let\Height@\relax\let\Depth@\relax\let\Rowheight@\relax
 \let\Rowdepth@\relax\let\Colwidth@\relax
}
\let\enddocument\bye
\def\alloc@#1#2#3#4#5{\global\advance\count1#1by\@ne
  \ch@ck#1#4#2%
  \allocationnumber=\count1#1%
  \global#3#5=\allocationnumber
  \wlog{\string#5=\string#2\the\allocationnumber}}
\catcode`\@=\active

\input epsf

\loadmsam
\loadmsbm

\magnification1200

\NoBlackBoxes
\RefWarnings

\font\grosso=cmbx12

\def\mm{\overline{\Cal M}}
\def\m{\Cal M}
\def\pu{\Cal P}
\def\del{\partial}
\def\cic{\Bbb C}
\def\pip{\Bbb P}
\def\quq{\Bbb Q}

\def\sez#1#2{\nopunct\hl1"#1."{#2}
\smallskip
\Reset\tag1
\Reset\claim1
\newpre\tag{#1.}
\newpre\claim{(#1.}}

\newpost\claim{)}
\def\clail#1#2{\claim{#1}{\label{#2}}\Offset\tag2}
\def\clailm#1#2#3{\claim{#1}{\label{#2}#3}\Offset\tag2}
\def\tal#1{\tag\label{#1}\Offset\claim2}

\document

\centerline{\grosso Calculating cohomology groups of moduli spaces of curves}
\centerline{\grosso via algebraic geometry}
\vskip15pt
\centerline{Enrico Arbarello and Maurizio Cornalba\footnote{\sevenrm Both
authors
supported in part by the Human Capital and Mobility Programme of the European
Union
under contract ERBCHRXCT940557 and by the MURST national project ``Geometria
Algebrica''}}
\vskip20pt

In this paper we compute the first, second, third, and fifth rational
cohomology groups
of  $\mm_{g,n}$, the moduli space of stable $n$-pointed genus $g$ curves. It
turns out
that $H^1(\mm_{g,n},\quq)$, $H^3(\mm_{g,n},\quq)$, and $H^5(\mm_{g,n},\quq)$
are
zero for all values of $g$ and $n$, while $H^2(\mm_{g,n},\quq)$ is generated by
tautological classes, modulo relations that can be written down explicitly; the
precise
statements are given by Theorems \Ref{H1H3} and \Ref{H2}. We are convinced that
the computation of the fourth cohomology of all moduli spaces
$\mm_{g,n}$ should also be accessible to our methods.

It must be observed that some of these results are
not new. In fact, it is known that $\mm_{g,n}$ is simply connected (cf.
[\Ref{BoPik}], for
instance), while Harer has determined $H^2(\m_{g,n},\quq)$
[\Ref{HarerH2}]; once this is known, it is not hard to compute the
corresponding group
for $\mm_{g,n}$. Harer [\Ref{HarerH3}] has also shown that $H^3(\m_{g,n},\quq)$
vanishes, at least for large enough genus. What is really new here is the
method of proof,
which is mostly based on standard algebro-geometric techniques, rather than
geometric
topology. Especially for odd cohomology, this provides proofs that are quite
short and,
we hope, rather transparent. It should also be noticed that the odd cohomology
of
$\mm_{g,n}$, at least in the range we can deal with, seems to be somewhat
better
behaved than the one of $\m_{g,n}$, for it is certainly not the case that the
first and
third cohomology groups of $\m_{g,n}$ are always zero.

Roughly speaking, the idea of the proof is as follows. If one could apply the
Lefschetz
hyperplane theorem, one might reduce the computation of $H^k(\mm_{g,n},\quq)$
to the
one of $H^k(\del\m_{g,n},\quq)$, for low $k$. Although the standard Lefschetz
theorem
cannot be used, since $\del\m_{g,n}$ is almost never ample, a foundational
result of
Harer, which is a direct consequence of the construction of a cellular
decomposition of
$\m_{g,n}$ by means of Strebel differentials, provides a suitable substitute.
Once this is
established, a little Hodge theory shows that, always for low enough $k$,
$H^k(\mm_{g,n},\quq)$ injects not only in
$H^k(\del\m_{g,n},\quq)$, but also in the $k$-th cohomology group of the
normalization
$N$ of $\del\m_{g,n}$. Put otherwise, the combinatorics of the boundary does
not
contribute to $H^k(\mm_{g,n},\quq)$. Since the components of $N$ are
essentially products of moduli spaces $\mm_{\gamma,\nu}$ such that either
$\gamma<g$ or $\gamma=g$ and $\nu<n$, one may try to compute
$H^k(\mm_{g,n},\quq)$ by double induction on $g$ and $n$, starting from a few
seed
case to be handled directly. This turns out to be possible, and reduces to
elementary
linear algebra. An interesting, and somewhat unexpected, byproduct of the proof
is that,
for any $k$, the $k$-th cohomology group of $\mm_{g,n}$ injects into the $k$-th
cohomology of the normalization of the component of the boundary parametrizing
irreducible singular curves, provided $g$ is large enough.

There are other cases, in addition to the ones mentioned above, in which the
cohomology of moduli spaces of curves has been computed. First of all, Harer
[\Ref{HarerH4}] has computed the fourth cohomology of $\m_{g,n}$ for large
enough $g$.
In a different direction, the entire cohomology ring
of $\mm_{0,n}$ has been described for any $n$ by Keel [\Ref{Keel}] in terms of
generators and relations. Likewise, Getzler [\Ref{Ezra}] has announced the
computation
of the {\it even} cohomology ring $H^{even}(\mm_{1,n},\quq)$ for any $n$.
Mumford
[\Ref{Mumfenum}] and Faber [\Ref{Faber}] have computed the cohomology rings of
$\mm_{2}$ and
$\mm_{2,1}$, while Getzler [\Ref{EzraM22}] has computed the cohomology ring of
$\mm_{2,2}$ and announced the determination of the cohomology of $\mm_{2,3}$.
Finally, Looijenga [\Ref{LooijM3}] has computed the cohomology of $\m_{3}$ and
$\m_{3,1}$.

We will assume Keel's result, which is derived entirely via algebraic geometry,
although
the part of it that we use could be proved without much effort by our methods.
Some of
Getzler's and Looijenga's results will be needed, while computing $H^5$, to
deal
with some of the initial cases of the induction; except for this, our treatment
will be
self-contained and for the most part quite elementary.

We are grateful to Eduard Looijenga for indicating to us that Proposition
\Ref{injectinbdry} could be proved via Hodge theory; our original proof was
based on a
fairly involved combinatorial argument and worked with certainty only for $k\le
2$.

\sez{1}{Tautological classes}
The purpose of this section is to fix notation and to collect a few results
about
tautological cohomology classes that are well known to specialists but for
which a
comprehensive reference seems to be lacking. All the varieties we shall
consider will be
over $\cic$. Only rational cohomology will be used; when we omit mention of the
coefficient group, we always implicitly assume rational coefficients.

Let $g$ and $n$ be non-negative integers such that $2g-2+n>0$. We
denote by $\mm_{g,n}$ the moduli space of stable $n$-pointed genus $g$ curves
and by
$\m_{g,n}$ its subspace parametrizing smooth curves. More generally, if $P$ is
a set with $n$ elements, it will be
technically convenient to consider also stable $P$-pointed curves. These are
simply
stable curves whose marked points are indexed by
$P$, and not by $\{1,\dots,n\}$. We shall denote by $\mm_{g,P}$ and $\m_{g,P}$
the
corresponding moduli spaces. The {\it boundary} of $\m_{g,P}$ is
$\del\m_{g,P}=\mm_{g,P}\setminus\m_{g,P}$.

By a {\it graph} we shall mean the datum $G$ of:
\list
\item"-"a non-empty finite set $V=V(G)$ (the {\it vertices} of $G$),
\item"-"a non-negative integer $g_v$ for every $v\in V$,
\item"-"a finite set $L=L(G)$ (the {\it half-edges} of $G$),
\item"-"a partition $\pu$ of $L$ in subsets with one or two elements,
\item"-"a subset $L_v\subset L$ for every $v\in V$,
\endlist
with the property that
$$
L=\coprod_{v\in V}L_v\,.
$$
Whe shall call the elements of $\pu$ with one element {\it legs} of $G$ and
those with
two elements {\it edges}; the set of all the latter will be denoted $E=E(G)$.
We also
set $l_v=\vert L_v\vert$.
In what follows we shall implicitly consider only graphs which are {\it
connected}, in an
obvious sense. If $P$ is a finite set, by a $P$-{\it labelled graph} we shall
mean the
datum of a graph $G$ plus a bijection between the set of its legs and $P$.

To every stable $P$-pointed genus $g$ curve $(C;\{q_p\}_{p\in P})$ we may
associate a
$P$-labelled graph $G$ as follows. Let $\pi:N\to C$ be the normalization of
$C$. We let
$V(G)$ be the set of all components of $N$, and $L(G)$ the set of all points of
$N$ which
map either to nodes or to marked points of $C$; two of these constitute an edge
if they
map to the same point (a node) of $C$, while the remaining ones are legs. The
indexing
of the legs by $P$ is the obvious one. We also set $g_v=\text{genus of }v$, and
let $L_v$
be the set of all elements of $L$ belonging to $v$. Notice that
$$
g=\sum_{v\in V(G)}g_v+1-\vert V(G)\vert+\vert E(G)\vert\,.
$$
Conversely, this formula can be used to define the {\it genus} of any connected
$P$-labelled graph. The graph associated to a stable $P$-pointed genus $g$
curve is {\it
stable} in the sense that $2g_v-2+l_v>0$ for every vertex $v$; more exactly, to
say that
a curve is stable is equivalent to saying that its graph is.

Let $G$ be a connected stable $P$-labelled graph of genus $g$. We denote by
$\m(G)$
the moduli space of all $P$-pointed genus $g$ stable curves whose associated
graph is
$G$; it a locally closed subspace of $\mm_{g,P}$ of codimension $\vert
E(G)\vert$. We
also denote by $\delta_G$ the orbifold fundamental class of $\overline{\m(G)}$,
that is,
the crude fundamental class divided by the order of the automorphism group of a
general element of $\m(G)$. The degree two classes correspond to graphs with
one
edge. These come in two kinds; there is the graph $G_{irr}$, with one vertex of
genus
$g-1$, and there are the graphs $G_{a,A}$, which have two vertices, one of
genus $a$, with attached the legs indexed by $A$, and one of genus $g-a$, with
attached
the legs indexed by $A^c=P\setminus A$. We shall normally write
$\delta_{irr}$ and $\delta_{a,A}$ instead of $\delta_{G_{irr}}$ and
$\delta_{G_{a,A}}$.
We shall also write $\Delta_{irr}$ and $\Delta_{a,A}$ to denote
$\overline{\m(G_{irr})}$
and $\overline{\m(G_{a,A})}$, respectively. It is clear that
$$
\delta_{a,A}=\delta_{g-a,A^c}\,,\tal{complrel}
$$
and also that $\delta_{a,A}$ does not make sense as defined unless $G_{a,A}$ is
a stable
graph, i.e., unless $2a-2+\vert A\vert\ge 0$ and $2(g-a)-2+\vert A^c\vert\ge
0$. In
practice, this means that $\delta_{a,A}$ is still undefined if
$a=0$ and $\vert A\vert<2$ or
$a=g$ and $\vert A\vert>\vert P\vert-2$. We will find it convenient to set
$\delta_{a,A}$ to zero if $a<0$, $a>g$, $2a-2+\vert A\vert<0$, or
$2(g-a)-2+\vert A^c\vert<0$.
 The class
$\delta_{irr}$, as defined above, also does not make sense in genus zero; we
will set it to
zero in this case.

The basic maps between moduli spaces are
$$
\aligned
&\pi:\mm_{g,P\cup\{q\}}\to\mm_{g,P}\,,\\
&\xi:\mm_{g-1,P\cup\{q,r\}}\to\mm_{g,P}\,,\\
&\eta:\mm_{a,A\cup\{q\}}\times\mm_{g-a,A^c\cup\{r\}}\to\mm_{g,P}\,,
\endaligned
$$
which are defined as follows. The image under $\pi$ of a $P\cup\{q\}$-pointed
curve is
obtained by ``forgetting'' the point labelled by $q$ and passing to the stable
model. The
image under $\xi$ of a $P\cup\{q,r\}$-pointed genus $g-1$ curve is obtained by
identifying the points labelled by $q$ and $r$; likewise, the image under
$\eta$ of a pair
consisting of an $A\cup\{q\}$-pointed genus $a$ curve and an
$A^c\cup\{r\}$-pointed
genus $g-a$ curve is the $P$-pointed curve of genus $g$ obtained by identifying
the
points labelled by $q$ and $r$.

The map $\pi:\mm_{g,P\cup\{q\}}\to\mm_{g,P}$ is also called the {\it universal
curve}
over $\mm_{g,P}$. It has $\vert P\vert$ sections, indexed by $P$; the section
$\sigma_p$ attaches to any $P$-pointed curve $(C;\{x_i\}_{i\in P})$ the
$P\cup\{q\}$-pointed curve obtained by attaching to $C$ a copy of $\pip^1$ by
identifying $x_p$ and $0\in \pip^1$, and labelling the points $1$ and $\infty$
by $p$
and $q$. One may use the universal curve to define further cohomology classes
on
$\mm_{g,P}$ as follows. We denote by $\omega_\pi$ the relative dualizing sheaf
and by
$D_p$ the image of $\sigma_p$. One then sets
$$
\aligned
\psi_p&=\sigma_p^*(c_1(\omega_\pi))\,,\qquad p\in P\,,\\
\kappa_i&=\pi_*(c_1(\omega_\pi(\textstyle{\sum}D_p))^{i+1})\,,\qquad i\ge 0\,.
\endaligned
$$
The classes $\psi_p$ have degree 2, while $\kappa_i$ has degree $2i$. In the
rest of
this paper, whenever we speak of {\it tautological} or {\it natural} classes
(of degree 2)
on $\mm_{g,P}$, we refer to $\kappa_1$, the $\psi_p$, $\delta_{irr}$, and the
$\delta_{a,A}$. The classes $\delta_{irr}$ and $\delta_{a,A}$ will be called
{\it
boundary} classes.

Our next task is to describe how the natural classes pull back under $\pi$,
$\xi$, and
$\eta$.

\clail{Lemma}{pipullbacks}i) $\pi^*(\kappa_1)=\kappa_1-\psi_q$;
\item{ii)}$\pi^*(\psi_p)=\psi_p-\delta_{0,\{p,q\}}$ for any $p\in P$;
\item{iii)}$\pi^*(\delta_{irr})=\delta_{irr}$;
\item{iv)}$\pi^*(\delta_{a,A})=\delta_{a,A}+\delta_{a,A\cup\{q\}}$.
\endclaim
Part {\sl i)} of the lemma is proved in [\Ref{ArbCo}], while {\sl iii)} and
{\sl iv)} are
clear. To prove {\sl ii)} we reason as follows. Consider the diagram
$$
\cgaps{1;.07}
\CD
\mm_{g,P\cup\{q,r\}}@()\L{\mu}@(1,0)@()\L{\varphi'}@(0,-1)&
\mm_{g,P\cup\{r\}}@()\L{\varphi}@(0,-1)&\cr
\mm_{g,P\cup\{q\}}@()\L{\pi}@(1,0)&\mm_{g,P}&,
\endCD
$$
where $\varphi'$ and $\varphi$ are defined as ``forgetting the point labelled
by $r$'' and $\mu$ as ``forgetting the point labelled by $q$''. It is known
(cf.
[\Ref{ArbCo}], for instance) that
$$
\mu^*(\omega_{\varphi})=\omega_{\varphi'}(-\sum_{x\in P}
\Delta_{0,\{x,q,r\}})\,.
$$
Thus, if $\tau_x$, $x\in P$ (resp., $\tau'_x$, $x\in P\cup\{q\}$), are the
canonical
sections of $\varphi$ (resp., $\varphi'$), then
$$
\pi^*(\tau_p^*(\omega_{\varphi}))={\tau'_p}^*(\mu^*(\omega_{\varphi}))
={\tau'_p}^*(\omega_{\varphi'}(-\sum_{x\in P} \Delta_{0,\{x,q,r\}}))
$$
for any $p\in P$. This translates into {\sl ii)}, finishing the proof of the
lemma.

\clail{Lemma}{xipullbacks}i) $\xi^*(\kappa_1)=\kappa_1$;
\item{ii)} $\xi^*(\psi_p)=\psi_p$ for any $p\in P$;
\item{iii)} $\xi^*(\delta_{irr})=\delta_{irr}-\psi_q-\psi_r
+\displaystyle{\sum_{q\in B, r\not\in B}}\delta_{b,B}$;
\item{iv)} $\xi^*(\delta_{a,A})=\left\{\matrix\delta_{a,A}\hfill&&\text{if }
g=2a,\ A=P=\emptyset,\hfill\\
\delta_{a,A}+\delta_{a-1,A\cup\{q,r\}}\hfill&&\text{otherwise}.\hfill
\endmatrix\right.$
\endclaim
Part {\sl i)} is proved in [\Ref{ArbCo}], while the other parts of the lemma
are
straightforward. We now turn to $\eta$. We shall not actually compute the
pullbacks of
the natural classes under $\eta$, but only under the map
$$
\vartheta:\mm_{a,A\cup\{q\}}\to\mm_{g,P}
$$
which associates to any
$A\cup\{q\}$-pointed genus $a$ curve the $P$-pointed genus $g$ curve
obtained by glueing to it a {\it fixed} $A^c\cup\{r\}$-pointed genus $g-a$
curve $C$ via
identification of $q$ and $r$. On the other hand, as we announced, the first
cohomology
of $\mm_{\gamma,\nu}$ always vanishes, by the K\"unneth formula, so the second
cohomology of $\mm_{a,A}\times \mm_{g-a,A^c}$ is the direct sum of
$H^2(\mm_{a,A})$ and $H^2(\mm_{g-a,A^c})$. Thus knowing how the natural classes
pull back under $\vartheta$ actually tells us how they pull back under $\eta$
(once the
vanishing of the first cohomology has been proved). It is important to stress
that,
although of course $\vartheta$ depends on the choice of $C$, any two choices
give rise
to homotopic maps so that, in cohomology, the pullback map
$\vartheta^*$ is independent of the choice of $C$.

\clail{Lemma}{thetapullbacks}i) $\vartheta^*(\kappa_1)=\kappa_1$;
\item{ii)} $\vartheta^*(\psi_p)=\left\{\matrix\psi_p\hfill&&\text{if }p\in
A,\hfill\\
0\hfill&&\text{if }p\in A^c;\hfill
\endmatrix\right.$
\item{iii)} $\vartheta^*(\delta_{irr})=\delta_{irr}$.

\noindent Suppose $A=P$. Then
\item{iv)} $\vartheta^*(\delta_{b,B})=\left\{
\matrix
\delta_{2a-g,P\cup\{q\}}-\psi_q\hfill&&\text{if }(b,B)=(a,P)\text{ or }
(b,B)=(g-a,\emptyset),\hfill\\
\delta_{b,B}+\delta_{b+a-g,B\cup\{q\}}\hfill&&\text{otherwise}.\hfill
\endmatrix
\right.$

\noindent Suppose $A\neq P$. Then
\item{iv')} $\vartheta^*(\delta_{b,B})=\left\{
\matrix
-\psi_q\hfill&&\text{if }(b,B)=(a,A)\text{ or }
(b,B)=(g-a,A^c),\\
\delta_{b,B}\hfill&&\text{if }B\subset A\text{ and }(b,B)\neq (a,A),\hfill\\
\delta_{b+a-g,(B\setminus A^c)\cup\{q\}}\hfill&&\text{if }B\supset A^c\text{
and
}(b,B)\neq (g-a,A^c),\hfill\\ 0\hfill&&\text{otherwise}.\hfill
\endmatrix
\right.$
\endclaim
Again, the only part that needs justification is {\sl i)}, which is proved in
[\Ref{ArbCo}].

\bigskip
We may now determine all relations among tautological classes in degree two.
The
answer depends on the genus. We begin with genus zero. In this case it has been
observed by Keel [\Ref{Keel}] that, for any four distinct elements $p,q,r,s$ of
$P$, the
following relations hold among the classes $\delta_{0,A}$ such that $\vert
A\vert\ge 2$
and
$\vert A^c\vert\ge 2$:
$$
\sum_{A\ni p,q\atop A\not\ni r,s}\delta_{0,A}
=\sum_{A\ni p,r\atop A\not\ni q,s}\delta_{0,A}
=\sum_{A\ni p,s\atop A\not\ni q,r}\delta_{0,A}\,.\tal{Keelrel}
$$
What is more important, Keel proves that $H^2(\mm_{0,P})$ is the quotient of
the vector
space generated by the $\delta_{0,A}$ such that $\vert A\vert\ge 2$
and $\vert A^c\vert\ge 2$ modulo the trivial relations \Ref{complrel} and the
relations
\Ref{Keelrel} for all possible choices of $p,q,r,s$.

\clail{Proposition}{tautog0}For any choice of distinct elements $x,y,z\in P$,
the following
relations hold in $H^2(\mm_{0,P})$:
$$
\align
\psi_z&=\sum_{A\ni z\atop A\not\ni x,y}\delta_{0,A}\,,\tal{psig0}\\
\kappa_1&=\sum_{A\not\ni x,y}(\vert A\vert-1)\delta_{0,A}\,.\tal{kappag0}
\endalign
$$
These, together with the relation $\delta_{irr}=0$ and relations \Ref{complrel}
and
\Ref{Keelrel}, generate all relations in $H^2(\mm_{0,P})$ among the natural
classes
$\kappa_1$, $\psi_i$, $\delta_{irr}$, and $\delta_{0,A}$ with $\vert
A\vert\ge 2$ and $\vert A^c\vert\ge 2$.
\endclaim

In view of Keel's result, all that needs to be shown is that \Ref{psig0} and
\Ref{kappag0}
hold. The proof is by induction on $\vert P\vert$, starting from the obvious
remark
that $0=\kappa_1=\psi_1=\psi_2=\psi_3$ on $\mm_{0,P}$ when $\vert P\vert=3$.
The
induction step is based on Lemma \Ref{pipullbacks}. Suppose that \Ref{psig0}
and \Ref{kappag0} hold in $\mm_{0,P}$. By simmetry, it suffice to prove their
analogues
in $\mm_{0,P\cup\{q\}}$ for $x,y,z\in P$. Pulling back \Ref{psig0} via $\pi$
gives that
$$
\psi_z-\delta_{0,\{z,q\}}=\sum_{z\in A\subset P\atop x,y\not\in
A}\delta_{0,A}+\delta_{0,A\cup\{q\}}\,,
$$
which is nothing but the analogue of \Ref{psig0} for $\mm_{0,P\cup\{q\}}$.
Similarly,
pulling back \Ref{kappag0} yields
$$
\kappa_1-\psi_q=\sum_{x,y\not\in
A\subset P}(\vert A\vert-1)(\delta_{0,A}+\delta_{0,A\cup\{q\}})\,,
$$
that is, using \Ref{psig0} to express $\psi_q$ in terms of boundary classes,
$$
\kappa_1=\sum_{x,y\not\in
A\subset P}\delta_{0,A\cup\{q\}}+\sum_{x,y\not\in
A\subset P}(\vert A\vert-1)(\delta_{0,A}+\delta_{0,A\cup\{q\}})\,,
$$
which is exactly what had to be shown.

\bigskip
To state and prove the analogues of \Ref{tautog0} in higher genus, it will be
convenient
to adopt the following notational conventions. We set $\psi=\sum \psi_i$ and,
for any
integer $a$, we let $\delta_a$ be the sum of all classes $\delta_{a,A}$; notice
that, in
case $g=2a$, the summand $\delta_{a,A}=\delta_{a,A^c}$ occurs only once, and
not
twice, in this sum. We also let $\delta$ denote the sum of $\delta_{irr}$ and
of all the
$\delta_a$ with $2a\le g$. Finally, as is customary, we denote by $\lambda$ the
Hodge
class, that is, $c_1(\pi_*(\omega_\pi))$.

\clail{Proposition}{tautog>0}i) The following relations hold in
$H^2(\mm_{1,P})$, for any
$p\in P$:
$$
\align
\kappa_1&=\psi-\delta_0\,,\tal{kappag1}\\
12\psi_p&=\delta_{irr}+12\sum_{{S\ni p}\atop{\vert S\vert\ge
2}}\delta_{0,S}\,.\tal{psig1}
\endalign
$$
These, together with the \Ref{complrel}, generate all relations in
$H^2(\mm_{1,P})$ among the natural classes
$\kappa_1$, $\psi_i$, $\delta_{irr}$, and $\delta_{a,A}$ with $0\le a\le
1$ and $2\le\vert A\vert\le \vert P\vert-2$ if $a=0$.
\medskip\noindent
ii) The following relation holds in $H^2(\mm_{2,P})$:
$$
5\kappa_1=5\psi+\delta_{irr}-5\delta_0+7\delta_1\,.\tal{kappag2}
$$
This relation and the \Ref{complrel} generate all relations in
$H^2(\mm_{2,P})$ among the natural classes
$\kappa_1$, $\psi_i$, $\delta_{irr}$, and $\delta_{a,A}$ with $0\le a\le
2$ and $2\le\vert A\vert\le \vert P\vert-2$ if $a=0$.
\medskip\noindent
iii) If $g\ge 3$, the \Ref{complrel} generate all relations in $H^2(\mm_{g,P})$
among the
classes $\kappa_1$, $\psi_i$, $\delta_{irr}$, and $\delta_{a,A}$ with $0\le
a\le g$ and
$2\le\vert A\vert\le \vert P\vert-2$ if $a=0$.
\endclaim

The proofs of \Ref{kappag1}, \Ref{psig1}, and \Ref{kappag2} are the exact
analogues of
those of \Ref{kappag0} and \Ref{psig0}. The initial cases of the induction are
as follows.
First of all, for any $g$ and any $P$, one has Mumford's relation
[\Ref{Mumfenseign}][\Ref{MDT}]
$$
\kappa_1=12\lambda-\delta+\psi\,.\tal{Mumfrel}
$$
For $g=\vert P\vert=1$, one knows that $\psi=\lambda$, and that
$12\lambda-\delta=0$ [\Ref{HarrisMumf}], so that $\kappa_1=\psi$ and
$\delta=12\psi$. Since in this case $\delta=\delta_{irr}$ and $\delta_0=0$,
these
identities are just the relations \Ref{kappag1} and \Ref{psig1}. For $g=2$,
$P=\emptyset$, Mumford [\Ref{Mumfenum}] has shown that
$10\lambda=\delta_{irr}+2\delta_1$. Coupled with \Ref{Mumfrel}, this says that
$$
\aligned
5\kappa_1&=60\lambda-5\delta+5\psi=6\delta_{irr}+12\delta_1-5\delta_{irr}
-5\delta_0-5\delta_1+5\psi\\
&=5\psi+\delta_{irr}-5\delta_0+7\delta_1\,,
\endaligned
$$
as desired.

What remains to be shown is that there are no relations in addition to the ones
listed above. We begin with the case of genus 1. We need the following simple
remark.

\clail{Lemma}{deltairrtozero}The homomorphism $\xi^*:H^2(\mm_{1,P})\to
H^2(\mm_{0,P\cup\{q,r\}})$ maps $\delta_{irr}$ to zero.
\endclaim

Let $p$ be an element of $P$, and let $\rho:\mm_{1,P}\to\mm_{1,\{p\}}$ be the
morphism defined by forgetting the points labelled by elements of $P$ other
than $p$.
The cycle $\Delta_{irr}$ is the inverse image of the point at infinity
$x_0\in\mm_{1,\{p\}}$. Therefore the cohomology class $\delta_{irr}$ is the
fundamental class of $\rho^{-1}(x)$, where $x$ is any point of $\mm_{1,\{p\}}$.
Since
$\rho^{-1}(x)$ does not touch $\Delta_{irr}=\xi(\mm_{0,P\cup\{q,r\}})$ if
$x\neq x_0$, it
follows that $\xi^*(\delta_{irr})=0$, as desired.

\bigskip
What we need to do to finish the genus 1 case is to show that $\delta_{irr}$
and the
classes $\delta_{1,S}$ are independent.
First consider the inclusion
$\vartheta:\mm_{1,1}\hookrightarrow\mm_{1,P}$ obtained by sending any 1-pointed
genus 1 curve $C$ to the union of $C$ with a fixed
$P\cup\{q\}$-pointed smooth rational curve E, with the marked point of $C$
identified
with the point of $E$ labelled by $q$. Notice that, by \Ref{thetapullbacks},
$\vartheta^*(\delta_{irr})=\delta_{irr}$, so
$\delta_{irr}\in H^2(\mm_{1,P})$ is not zero since its pullback to
$H^2(\mm_{1,1})$ does not vanish.

Now look at the pullbacks of the boundary classes via the morphism
$\xi:\mm_{0,P\cup\{q,r\}}\to\mm_{1,P}$. We have seen that $\delta_{irr}$ pulls
back to
zero. On the other hand it is clear that
$\xi^*(\delta_{1,S})=\delta_{0,S\cup\{q,r\}}$.
The independence of the classes $\delta_{irr}$ and $\delta_{1,S}$ will then
follow from
the remark that $\delta_{irr}\neq 0$ and from the following result.

\clail{Lemma}{T3}The classes $\delta_{S\cup\{q,r\}}$, where $S$ runs
through all subsets of $P$ with at most $\vert P\vert-2$ elements, are
independent in
$H^2(\mm_{0,P\cup\{q,r\}})$.
\endclaim

This would follow, after some work, from Keel's description of the cohomology
ring of
$\mm_{0,P\cup\{q,r\}}$, but we choose a different approach. To begin with, it
suffices
to do the proof when $P=\{1,\dots,n\}$, $q=n+1$, $r=n+2$. The lemma is clearly
true for
$n\le 2$, so we can assume that $n>2$. For each
$S$, we shall produce a 1-parameter family $F_S$ of stable $(n+2)$-pointed
genus zero
curves and will show that the matrix whose entries are the degrees of the
classes
$\delta_{T\cup\{n+1,n+2\}}$ on these families is non-degenerate.

Let $f:X=\Bbb P^1\times \Bbb P^1\to \Bbb P^1$ be the projection onto the first
factor, and let $D'_i$, $i\in \{1,\dots,n+3\}$ be disjoint constant sections of
$f$. Suppose
first that $S\neq\emptyset$. We then let $F_S'$ be the family of pointed curves
over
$\Bbb P^1$ consisting of the blowup of $f:X\to \Bbb P^1$ at the points of
intersection of
the diagonal with the sections $D_i'$, $i\in S\cup\{n+1,n+2\}$, together with
the sections
$D_i$, $i\in S\cup\{n+1,n+2,n+4\}$, where $D_i$ is the proper transform of
$D_i'$
for$i\neq n+4$, and $D_{n+4}$ is the proper transform of the diagonal. We also
let
$F_S''$ be the family of pointed curves over $\Bbb P^1$ consisting of $f:X\to
\Bbb P^1$
together with the sections $D_i=D_i'$, $i\in\{1,\dots,n,n+3\}-S$. When
$S=\emptyset$, the definitions of $F_S'$ and $F_S''$ are different. The family
$F_S'$ is just $f:X\to \Bbb P^1$, together with the sections $D_i=D_i'$,
$i=n+1,n+2,n+3$.
As for $F_S''$, blow up $f:X\to \Bbb P^1$ at the points where the diagonal
meets the
sections $D_i'$, $i=1,\dots,n$, and take as sections $D_i$, $i=1,\dots,n,n+4$,
the proper
transforms of the $D_i'$, $i\le n$, and of the diagonal. In any case, $F_S$ is
defined to be
the family of stable $(n+2)$-pointed genus zero curves over $\Bbb P^1$ obtained
from
$F_S'$ and $F_S''$ by identifying the sections $D_{n+3}$ and $D_{n+4}$.

The degrees of the classes $\delta_{T\cup\{n+1,n+2\}}$ are readily calculated.
They are
as follows.
\medskip
\item{-} $S\neq\emptyset$
$$
\alignedat2
\deg_{F_S}\delta_{S\cup\{n+1,n+2\}}&=-\vert S\vert\,,\\
\deg_{F_S}\delta_{S-\{s\}\cup\{n+1,n+2\}}&=1&\text{if } s\in S\,,\\
\deg_{F_S}\delta_{T\cup\{n+1,n+2\}}&=0&\text{otherwise.}
\endalignedat
$$
\medskip
\item{-} $S=\emptyset$
$$
\alignedat2
\deg_{F_S}\delta_{\{n+1,n+2\}}&=2-n\,,\\
\deg_{F_S}\delta_{\{s,n+1,n+2\}}&=1 &\text{if } s=1,\dots,n\,,\\
\deg_{F_S}\delta_{T\cup\{n+1,n+2\}}&=0 &\text{otherwise.}
\endalignedat
$$

\noindent Now suppose there is a relation $\sum_T
a_T\delta_{T\cup\{n+1,n+2\}}=0$.
Evaluating on $F_\emptyset$ yields
$$
(2-n)a_\emptyset+\sum_{s=1}^n a_{\{s\}}=0\,,
$$
while evaluating on $F_{\{s\}}$ gives
$$
-a_{\{s\}}+a_\emptyset=0\,,\qquad s=1,\dots,n\,.
$$
Combining these relations gives $a_\emptyset=0$ and $a_{\{s\}}=0$ for
$s=1,\dots,n$.
To show that, in fact, $a_S=0$ for any $S$, one proceeds by induction on $\vert
S\vert$.
We may assume that $S\neq\emptyset$. Suppose we know that $a_T=0$ whenever
$\vert T\vert<\vert S\vert$. Evaluating on
$F_S$ yields
$$
\vert S\vert a_S=\sum_{s\in S}a_{S-\{s\}}=0\,,
$$
so $a_S=0$. This finishes the proof of part {\sl i)} of \Ref{tautog>0}.

To complete the proof of part {\sl ii)} it remains to show that the boundary
classes and
the classes $\psi_p$, $p\in P$, are independent in $H^2(\mm_{2,P})$ modulo the
trivial
relations \Ref{complrel}. We distinguish several cases.
\medskip
\noindent {\bf Case 1: $P=\emptyset$}. Let $F'$ be a non-isotrivial family of
stable
1-pointed genus 1 curves over a smooth complete curve; by attaching a fixed
elliptic tail to the marked point of each fiber we get a family $F$ of stable
genus 2
curves. We then have that
$$
\aligned
\deg_F\delta_{irr}&=\deg_{F'}\delta_{irr}=12\deg_{F'}\lambda\,,\\
\deg_F\delta_{1,\emptyset}&=-\deg_{F'}\psi_1=-\deg_{F'}\lambda\,.
\endaligned\tal{T14}
$$
Next, let $C$ be a smooth elliptic curve, and let $x$ be a fixed point on $C$.
Let $X$ be
the surface obtained by blowing up $C\times C$ at the intersection point
between the
diagonal and $\{x\}\times C$, and let $X\to C$ be the composition of the
blow-down
map and of the projection from $C\times C$ onto the second factor. We get a
family $E$
of stable genus 2 curves over $C$ by identifying the proper transforms of the
diagonal
of $C\times C$ and of $\{x\}\times C$. For this family we have
$$
\deg_E\delta_{irr}=-2\,,\qquad
\deg_E\delta_{1,\emptyset}=1\,.
$$
Putting this together with \Ref{T14} shows that $\delta_{irr}$ and
$\delta_{1,\emptyset}$ are independent; in fact, the non-isotriviality of $F'$
implies
that $\deg_{F'}\lambda\neq 0$.
\medskip
\noindent {\bf Case 2: $\vert P\vert=1$}. Denote by $p$ the only point of $P$.
We
construct two families of stable 1-pointed genus 2 curves by modifying slightly
those
constructed in case 1. Let $F_1$ be as
$F$, but with the addition of the section traced out by a fixed point on the
constant
elliptic tail. To construct $E_1$, choose a point $y\in C$ distinct from $x$,
and blow up
$E\times E$ at the intersections between the diagonal, $\{x\}\times C$, and
$\{y\}\times
C$. Then identify the proper transforms of the diagonal and $\{x\}\times C$;
the
canonical section is the proper transform of $\{y\}\times C$. A third family
$G_1$ is
$\Gamma\times \Gamma\to \Gamma$, where $\Gamma$ is a smooth genus 2 curve,
with the diagonal as section. The invariants for these families are readily
calculated.
They are
$$
\aligned
\deg_{F_1}\delta_{irr}&=12\deg_{F'}\lambda\,,\qquad
\deg_{F_1}\delta_{1,\emptyset}=-\deg_{F'}\lambda\,,\qquad
\deg_{F_1}\psi_p=0\,,\\
\deg_{E_1}\delta_{irr}&=-2\,,\qquad
\deg_{E_1}\delta_{1,\emptyset}=1\,,\qquad
\deg_{E_1}\psi_p=1\,,\\
\deg_{G_1}\delta_{irr}&=0\,,\qquad
\deg_{G_1}\delta_{1,\emptyset}=0\,.\qquad
\deg_{G_1}\psi_p=2\,,\\
\endaligned
$$
Thus $\delta_{irr}$, $\delta_{1,\emptyset}$, and $\psi_p$ are independent.
\medskip
\noindent {\bf Case 3: $\vert P\vert>1$}. We proceed by induction on $\vert
P\vert$.
Let $p,q$ be distinct elements of $P$; set $Q=P-\{p,q\}$ and let $r,r'$ be
distinct and not
belonging to $P$. Let
$$
\vartheta:\mm_{2,Q\cup\{r\}}\simeq\mm_{2,Q\cup\{r\}}\times\mm_{0,\{r',p.q\}}
\to \mm_{2,P}
$$
be the map obtained by identifying the points labelled by $r$ and $r'$. It
follows from
\Ref{thetapullbacks} that
$$
\aligned
\vartheta^*(\psi_i)&=\psi_i\,,\qquad i\in Q\,,\\
\vartheta^*(\psi_p)&=\vartheta^*(\psi_q)=0\,,\\
\vartheta^*(\delta_{irr})&=\delta_{irr}\,,\\
\vartheta^*(\delta_{0,\{p,q\}})&=-\psi_r\,,\\
\vartheta^*(\delta_{a,A})&=\delta_{a,A}\,,\qquad a\le 1,\ A\subset Q\,,\\
\vartheta^*(\delta_{a,A\cup\{p,q\}})&=\delta_{a,A\cup\{r\}}\,,\qquad a\le 1,\
A\subset Q,\
A\neq\emptyset\text{ if }a=0\,,\\
\vartheta^*(\delta_{a,A})&=0\,,\qquad \{p,q\}\not\subset A,\ \{p,q\}\not\subset
P-A\,.
\endaligned\tal{T15}
$$
Suppose there is a relation
$$
\sum_{i\in P} a_i\psi_i+b\delta_{irr}+\sum_{\vert A\vert\ge 2}c_A\delta_{0,A}+
{\sum_{\vert A\vert\le \vert P\vert/2}}^{\hskip-10pt
'}\hskip7pt d_A\delta_{1,A}=0\,,\tal{T16}
$$
where the quotation mark affixed to the last summation symbol means that, for
each
$A\subset P$ such that $\vert A\vert=\vert P\vert/2$, {\it only one} of the
summands
$d_A\delta_{1,A}$, $d_{P-A}\delta_{1,P-A}$ occurs. Then the induction
hypothesis and
\Ref{T15} imply that $a_i=0$ unless $i=p$ or
$i=q$, that $b=0$, and that $c_A=d_A=0$ unless $\{p,q\}\not\subset A$,
$\{p,q\}\not\subset P-A$. Since $p$ and $q$ can be chosen at will, for $\vert
P\vert\ge
3$ this implies that all coefficients of \Ref{T16} vanish. If $P=\{p,q\}$, all
we can say is
that all the coefficients of \Ref{T16} are zero, except possibly for
$a_p$, $a_q$, and $d_{\{p\}}$. To see that these are zero as well, let $C$ be a
smooth
genus 2 curve, and let $x$ be a point on it. Blowing up $C\times C$ at the
intersection
point of $\{x\}\times C$ and the diagonal yields a family of genus 2 curves
with
two sections, namely the proper transforms of $\{x\}\times C$ and of the
diagonal.
Labelling the first of these with $p$ and the second with $q$, for the
resulting family of
stable $\{p,q\}$-pointed curves we have that
$$
\deg\psi_p=-1\,,\qquad\deg\psi_q=-3\,,\qquad\deg\delta_{1,\{p\}}=0\,.
$$
Since the roles of $p$ and $q$ can be interchanged, this implies that
$a_p=a_q=0$. To
conclude it suffices to produce a family of stable $\{p,q\}$-pointed curves of
genus 2
for which $\deg\delta_{1,\{p\}}$ does not vanish. One such is obtained by
attaching a
fixed elliptic tail with two extra fixed marked points to the points of the
canonical
section of the family $F'$ of elliptic curves considered in the analysis of
Case 1 above.
The proof of part {\sl ii)} of \Ref{tautog>0} is now complete.

To do part {\sl iii)} we proceed by induction on $g$. When
$P\neq\emptyset$, fix an element
$p\in P$; we have to show that $\kappa_1$, the $\psi_i$, $\delta_{irr}$ and the
$\delta_{a,A}$ such that $p\in A$ are independent. When $P=\emptyset$, instead,
we
have to show independence of $\kappa_1$, the $\psi_i$, $\delta_{irr}$ and the
$\delta_a$ with $2a\le g$. For $g>3$, the formulas in Lemma \Ref{xipullbacks},
plus
the induction hypothesis, show directly that the pullbacks of these classes via
$\xi:\mm_{g,P}\to\mm_{g-1,P\cup\{q,r\}}$ are already independent. For $g=3$, we
argue as follows. Suppose there is a relation
$$
0=a\kappa_1+\sum b_i\psi_i+c\delta_{irr}+\cdots
$$
among them. Pulling back via $\xi$, and using \Ref{xipullbacks} and
\Ref{kappag2}, we
find a relation on $\mm_{2,P\cup\{q,r\}}$ or the form
$$
0=(a-c)\psi_q+\dots+(c+a/5)\delta_{irr}+\cdots
$$
By {\sl ii)}, all coefficients in this relation must vanish, so $a=0$. At this
point we may
proceed as for $g>3$. The proof of \Ref{tautog>0} is now complete.

\sez{2}{The main results}
In this section we state our main results and present the core of their
proofs. The first theorem describes the first, third and fifth cohomology
groups of the
moduli spaces of stable curves.

\clail{Theorem}{H1H3}$H^k(\mm_{g,n})=0$ for $k=1,3,5$ and all $g$ and $n$ such
that
$2g-2+n>0$.
\endclaim

The next result describes the second cohomology group of $\mm_{g,n}$ in terms
of
generators and relations. It turns out that this group is always generated by
the natural
classes. The relations among these have already been determined in section 1.

\clail{Theorem}{H2}For any $g$ and $n$ such that
$2g-2+n>0$, $H^2(\mm_{g,n})$ is generated by $\kappa_1$,
the classes $\psi_i$, $\delta_{irr}$, and the classes $\delta_{a,A}$ such that
$0\le a\le
g$, $2a-2+\vert A\vert\ge 0$ and $2(g-a)-2+\vert A^c\vert\ge 0$. The relations
among
these classes are as follows.
\item{a)}If $g>2$ all relations are generated by those of the form
$$
\delta_{a,A}=\delta_{g-a,A^c}\,.\tal{complrelbis}
$$
\item{b)}If $g=2$ all relations are generated by the \Ref{complrelbis} plus the
following
one
$$
5\kappa_1=5\psi+\delta_{irr}-5\delta_0+7\delta_1\,.
$$
\item{c)}If $g=1$ all relations are generated by the \Ref{complrelbis} plus the
following
ones
$$
\align
\kappa_1&=\psi-\delta_0\,,\\
12\psi_p&=\delta_{irr}+12\sum_{{S\ni p}\atop{\vert S\vert\ge
2}}\delta_{0,S}\,.
\endalign
$$
\item{d)}If $g=0$ all relations are generated by the \Ref{complrelbis} plus the
following
ones$$
\align
\kappa_1&=\sum_{A\not\ni x,y}(\vert A\vert-1)\delta_{0,A}\,\\
\psi_z&=\sum_{A\ni z\atop A\not\ni x,y}\delta_{0,A}\,,\\
\delta_{irr}&=0\,.
\endalign
$$
\endclaim

Observe, first of all, that the moduli spaces $\m_{g,n}$ and $\mm_{g,n}$,
although in
general not smooth, are orbifolds; in particular, Poincar\'e duality holds for
them in
rational cohomology.

The proof of \Ref{H1H3} and \Ref{H2} begins with a simple remark. Look first at
$\m_{0,n}$; it can be viewed as the space of all $n$-tuples
$(0,1,\infty,z_4,\dots,z_n)$ of
distinct points of
$\pip^1$ or, which is the same, as the space of all $(n-3)$-tuples
$(z_4,\dots,z_n)$ in
$\cic^{n-3}$ such that
$z_i\neq 0,1$ for all $i$ and $z_i\neq z_j$ for all $i\neq j$. In other words,
$\m_{0,n}$
is nothing but $\cic^{n-3}$ minus a bunch of hyperplanes, so in particular it
is an
$(n-3)$-dimensional affine variety. It follows that $H_k(\m_{0,n})=0$ for
$k>n-3$. Things are similar in higher genus. In fact, when $g>0$, $n>0$, Harer
[\Ref{HarerVirtDim}] (see also [\Ref{Looij}]) constructs a
$(4g-4+n)$-dimensional spine
for $\m_{g,n}$; thus $H_k(\m_{g,n})$ vanishes for $k>4g-4+n$. The spine in
question is
constructed starting from the cellular decomposition of $\m_{g,n}$ defined in
terms of
Strebel differentials. Harer [\Ref{HarerVirtDim}] also shows, by a
spectral sequence argument, that $H_k(\m_g)=0$ for $k>4g-5$.
By Poincar\'e duality all this is equivalent to saying that the cohomology with
compact
support $H^k_c(\m_{g,n})$ vanishes for $k\le d(g,n)$, where
$$
d(g,n)=\left\{\matrix
n-4\hfill&&&\text{if }g=0,\hfill\\
2g-2\hfill&&&\text{if }n=0,\hfill\\
2g-3+n\hfill&&&\text{if }g>0,n>0.
\endmatrix
\right.\tal{dgn}
$$
Looking at the exact sequence of cohomology with compact supports
$$
\cdots\to H^k_c(\m_{g,n})\to H^k(\mm_{g,n})\to H^k(\del\m_{g,n})\to
H^{k+1}_c(\m_{g,n})\to\cdots
$$
then proves
\clail{Lemma}{1.10}The homomorphism $H^k(\mm_{g,n})\to H^k(\del\m_{g,n})$ is an
isomorphism for $k<d(g,n)$ and is injective for $k=d(g,n)$.
\endclaim

Let $g$ and $n$ be non-negative integers such that $2g-2+n>0$, and let $P$ be a
set
with $n$ elements. Denote by $D_1$, $D_2$,\dots the different irreducible
components of
$\del\m_{g,P}$. Each of these is the image of a map $\mu_i:X_i\to \mm_{g,P}$,
where
$X_i$ can be of two different kinds. Either $X_i=\mm_{g-1,P\cup\{q,r\}}$, or
else
$X_i=\mm_{a,A\cup\{q\}}\times \mm_{b,B\cup\{r\}}$, where $a+b=g$, $A\coprod
B=P$, and both $2a-2+\vert A\vert$ and $2b-2+\vert B\vert$ are non-negative. In
any
case $q$ and $r$ are distinct points not belonging to $P$, and the map $\mu_i$
is gotten
by identifying $q$ and $r$.

\clail{Lemma}{preinjectinbdry}The map $H^k(\mm_{g,P})\to \bigoplus_i
H^k(X_i)$ is injective whenever $H^k(\mm_{g,P})\to H^k(\del\m_{g,P})$ is.
\endclaim

The proof uses a bit of Hodge theory, in the form of the following result of
Deligne.

\clailm{Proposition}{deligne}{([\Ref{HodgeIII}], Proposition
(8.2.5))}Let $Y$ be proper. If $u:X\to Y$ is a proper surjective morphism, and
$X$
is smooth, then the weight $k$ quotient of $H^k(Y,\Bbb Q)$ is the image of
$H^k(Y,\Bbb Q)$ in $H^k(X,\Bbb Q)$.
\endclaim

In our application, $Y$ is $\del\m_{g,P}$, and $X$ is the disjoint union of the
$X_i$. Of
course, Deligne's result is stated for varieties and not for orbifolds, and the
$X_i$ are
smooth as orbifolds, but usually not as varieties. There are at least two ways
out. One is
to convince oneself that Deligne's proof works also in the orbifold context.
The other is to
appeal to the results of Looijenga [\Ref{Looij2}] and Boggi-Pikaart
[\Ref{BoPik}] which
imply that each of the $X_i$ is the quotient of a smooth variety $Z_i$ by the
action of a
finite group; one may then take as $X$ the disjoint union of the $Z_i$ and
prove
injectivity of the map $H^k(\mm_{g,P})\to H^k(X)=\bigoplus_i H^k(Z_i)$, which
obviously
implies the injectivity of $H^k(\mm_{g,P})\to \bigoplus_i H^k(X_i)$.

Whatever road we choose, the proof proceeds as follows. The homomorphism
$$
\rho:H^k(\mm_{g,P})\to H^k(\del\m_{g,P})
$$
is a morphism of mixed Hodge structures, and hence is strictly compatible
with the filtrations. Thus
$$
\rho(H^k(\mm_{g,P}))\cap
W_{k-1}(H^k(\del\m_{g,P}))=\rho(W_{k-1}(H^k(\mm_{g,P}))
=\rho(\{0\})=\{0\}\,,
$$
since $H^k(\mm_{g,P})$ is pure of weight $k$. As we are assuming that $\rho$ is
injective, this shows that $H^k(\mm_{g,P})$ injects
into $H^k(\del\m_{g,P})/W_{k-1}(H^k(\del\m_{g,P}))$.
On the other hand \Ref{deligne} says that
$H^k(\del\m_{g,P})/W_{k-1}(H^k(\del\m_{g,P}))$ injects into $H^k(X)$.

\bigskip
In view of \Ref{1.10}, an immediate corollary of Lemma \Ref{preinjectinbdry} is
the
following result.
\clail{Proposition}{injectinbdry}The map $H^k(\mm_{g,P})\to \bigoplus_i
H^k(X_i)$ is
injective whenever $k\le d(g,n)$, where $n=\vert P\vert$.
\endclaim

Proposition \Ref{injectinbdry} makes it possible to give a quick inductive
proof of
\Ref{H1H3}, based on the following intermediate result.

\clail{Lemma}{oddindscheme}Let $k$ be an odd integer, $h$ a non-negative
integer, and suppose $H^q(\mm_{g,n})=0$ for all odd $q\le k$, all $g\le h$, and
all $n$
such that $q>d(g,n)$. Then $H^q(\mm_{g,n})=0$ for all odd $q\le k$, all $g\le
h$, and all
$n$. If $H^q(\mm_{g,n})=0$ for all odd $q\le k$ and all $g$ and $n$ such that
$q>d(g,n)$, then $H^q(\mm_{g,n})=0$ for all odd $q\le k$ and all $g$ and $n$.
\endclaim
Clearly, it suffices to prove the first assertion. We argue by induction on
$k$. We may
assume, inductively, that $H^q(\mm_{g,n})=0$ for all odd $q< k$, all $g\le h$
and all
$n$. If $g\le h$ and $k\le d(g,n)$, Proposition \Ref{injectinbdry} says that
$H^k(\mm_{g,n})$ injects into a direct sum of vector spaces $H^k(X_i)$, where
$X_i$ is $\mm_{g-1,n+2}$ or a product of two moduli spaces
$\mm_{a,\alpha}$ and $\mm_{b,\beta}$ such that $a+b=g$,
$\alpha+\beta=n+2$; in the latter case either $a<g$ or $a=g$ and $\alpha<n$,
and
similarly for $b$ and $\beta$. By the K\"unneth formula, $H^k(\mm_{g,n})$
injects into
the direct sum of $H^k(\mm_{g-1,n+2})$ and of all the tensor products
$H^l(\mm_{a,\alpha})\otimes H^m(\mm_{b,\beta})$ with $l+m=k$. Since either $l$
or
$m$ must be odd, the induction hypothesis guarantees that all these tensor
products
vanish, except possibly those for which $l$ or $m$ is zero. This means that
$H^k(\mm_{g,n})$ injects into a direct sum of vector spaces
$H^k(\mm_{\gamma,\nu})$
such that $\gamma<g$ or $\gamma=g$ and $\nu<n$. But then the result follows by
double induction on $g$ and $n$.

\bigskip
Lemma \Ref{oddindscheme} reduces the proof of the vanishing of odd
cohomology (so long as it does vanish!) to checking it explicitly for finitely
many
values of $g$ and $n$ in each odd degree $k$, that is, those for which
$k>d(g,n)$. When
$k=1$ this means $g=0$, $n\le 4$ or $g=n=1$. Now, $\mm_{0,3}$ is a point, while
$\mm_{0,4}$ and $\mm_{1,1}$ are both isomorphic to the projective line, so the
first
cohomology groups of all three are zero. This concludes the proof of \Ref{H1H3}
in case
$k=1$. For $k=3,5$ the initial cases of the induction are slightly more
complicated and
will be dealt with in section 3. Lemma \Ref{oddindscheme} cannot be applied as
such if
$k\ge 11$, since it is known that $H^{11}(\mm_{1,11})$ is not zero. So far
as we know, the cases $k=7,9$ are still open.

\bigskip
We now turn to the proof of \Ref{H2}. Clearly, all that has to be
shown is that $H^2(\mm_{g,n})$ is always generated by tautological classes, the
relations among these having already been determined in section 1.
As was the case for \Ref{H1H3}, the proof is by double
induction on $g$ and $n$. Here we will describe the inductive step in genus 3
or more.
The cases of lower genus are a bit more involved, and will be treated in
section 4. The
initial cases of the induction, that is, those for which $2>d(g,n)$, will be
dealt with in
section 3.

Our strategy for the inductive step is quite simple. Suppose we want to show
that
$H^2(\mm_{g,n})$ is generated by tautological classes, assuming the same is
known to be
true in genus less then $g$, or in genus $g$ but with fewer than $n$ marked
points.
Proposition \Ref{injectinbdry} shows that $H^2(\mm_{g,n})$ injects into the
direct sum
of the second cohomology groups of the $X_i$. By induction hypothesis, these
are
generated by tautological classes, all relations among which are known. By
\Ref{xipullbacks} and \Ref{thetapullbacks}, we have complete control on the
effect of
each map $H^2(\mm_{g,n})\to H^2(X_i)$ on tautological classes, so that, at
least in
principle, we can decide which classes in $\bigoplus H^2(X_i)$ come from
tautological
classes on $H^2(\mm_{g,n})$. On the other hand, let $\alpha$ be any class in
$H^2(\mm_{g,n})$; if we denote by $\alpha_i$ its pullback to $X_i$, these
classes
satisfy obvious compatibility relations on the ``intersections'' of the $X_i$.
The subspace
of $\bigoplus H^2(X_i)$ defined by these compatibility relations can be
completely described using \Ref{xipullbacks} and \Ref{thetapullbacks}, at least
in
principle, because the spaces $H^2(X_i)$ are generated by tautological classes.
What we
will show, in essence, is that it coincides with the one generated by the
images of the
tautological classes of $H^2(\mm_{g,n})$. By the injectivity of
$H^2(\mm_{g,n})\to
\bigoplus H^2(X_i)$, this will conclude the proof.

A first step in the strategy outlined above is the following result, which is
also of independent interest.

\clail{Theorem}{injectindirr}Let $g\ge 1$ be an integer, let $P$ be a finite
set
such that $2g-2+\vert P\vert>0$, and let $q,r$ be distinct and not belonging to
$P$.
Then, if $\xi:\mm_{g-1,P\cup\{q,r\}}\to \mm_{g,P}$ is the morphism obtained by
identifying the points labelled by $q$ and $r$, the pullback map
$\xi^*:H^k(\mm_{g,P})\to H^k(\mm_{g-1,P\cup\{q,r\}})$ is injective for any
$k\le 2g-2$
if $g\le 7$, and for any $k\le g+5$ if $g\ge 7$.
\endclaim

This we will prove by triple induction on $k$, $g$ and $n=\vert P\vert$. The
statement is true when $k=0$, and also when $k=1$, since $H^1(\mm_{g,n})=0$ for
any
$g$ and $n$. Suppose then that $g\le 7$, that $k\le 2g-2$, and that the result
is
known to hold for all triples
$(k',g',n')$ such that either $k'<k$, or $k'=k$ and $g'<g$, or $k'=k$, $g'=g$,
and $n'<n$. In
view of
\Ref{injectinbdry}, what we have to show is that, if $x$ is any element of
$H^k(\mm_{g,P})$ such that $\xi^*(x)=0$, then $x$ pulls back to zero under any
one of
the maps
$$
\mm_{a,A\cup\{s\}}\times\mm_{b,B\cup\{t\}}\to\mm_{g,P}\,,
$$
where $g=a+b$ and $P=A\coprod B$. By
the K\"unneth formula, $H^k(\mm_{a,A\cup\{s\}}\times\mm_{b,B\cup\{t\}})$ breaks
up into a direct sum of summands $H^l(\mm_{a,A\cup\{s\}})\otimes
H^m(\mm_{b,B\cup\{t\}})$, where $k=l+m$. Thus we have to show that $x$ goes to
zero
under any one of the maps
$$
\rho:H^k(\mm_{g,P})\to
H^l(\mm_{a,A\cup\{s\}})\otimes H^m(\mm_{b,B\cup\{t\}})\,.
$$
Suppose $l\ge 2a-1$ and $m\ge 2b-1$; then $k=l+m\ge 2(a+b)-2=2g-2$. The only
possibility is that $l=2a-1$, $m=2b-1$; in particular, $l$ and $m$ are both
odd. Since
they add to $k\le 2g-2\le 12$, one of them must equal 5 or less. In view of
\Ref{H1H3},
this implies that $H^l(\mm_{a,A\cup\{s\}})\otimes H^m(\mm_{b,B\cup\{t\}})=0$,
so we
are done in this case. We may then suppose that either $l\le 2a-2$ or $m\le
2b-2$. Say
$l\le 2a-2$; this implies, in particular, that $a>0$. Then $\rho$ fits into a
commutative
diagram
$$
\CD
H^k(\mm_{g,P})@()\L{\rho}@(1,0)@()\L{\xi^*}@(0,-1)
&H^l(\mm_{a,A\cup\{s\}})\otimes H^m(\mm_{b,B\cup\{t\}})@()\L{\lambda}@(0,-1)\\
H^k(\mm_{g-1,P\cup\{q,r\}})@(1,0)
&H^l(\mm_{a-1,A\cup\{s,q,r\}})\otimes H^m(\mm_{b,B\cup\{t\}})
\endCD
$$
If $a=g$, that is, if $b=0$, then $\vert B\vert>1$, and hence $\vert
A\cup\{s\}\vert< n$. Thus $\lambda$ is always injective, by induction
hypothesis. Since
$\xi^*(x)=0$, and so $\lambda\rho(x)=0$, this implies that $\rho(x)=0$, as
desired.

When $g>7$ and $k\le g+5$, the argument is similar, but simpler. Set
$f(n)=2n-2$ if
$n\le 7$ and $f(n)=n+5$ if $n\ge 7$. To show that $\rho(x)=0$ we may argue as
in the
previous case, provided we can show that either $a>0$ and $l\le f(a)$ or $b>0$
and
$m\le f(b)$. If $a=0$, then $m\le k\le f(g)=f(b)$, and similarly if $b=0$; we
may
therefore assume that $a> 0$ and $b> 0$. Three cases are possible. Suppose
first that
$a\le 7$ and $b\le 7$. If $l>f(a)=2a-2$, $m>f(b)=2b-2$, then $k\ge 2g-2>g+5$,
against the
assumptions. Suppose next that $a\le 7$ and $b> 7$. If $l>f(a)=2a-2$,
$m>f(b)=b+5$, then
$k\ge a+g+5>g+5$, contrary to what we have assumed. Suppose finally that $a\ge
7$ and
$b\ge 7$. If $l>f(a)=a+5$, $m>f(b)=b+5$, then $k> g+10>g+5$, again against the
assumptions.

\clail{Remark}{injectindirrrem}{\rm The above argument actually shows that, if
$v$ is
an odd integer and we know that the odd cohomology of all the $\mm_{g,n}$
vanishes in
degree not exceeding
$v$, then $H^k(\mm_{g,P})\to H^k(\mm_{g-1,P\cup\{q,r\}})$ is injective
for $k\le 2g-2$ if $g\le v+2$, and for $k\le g+v$ if $g\ge v+2$. Thus
we could improve slightly on \Ref{injectindirr} if we could prove the
vanishing of the odd cohomology in degree greater than 5. In a different
direction, just knowing that the first cohomology vanishes, which is all that
has been fully proved up to now, suffices to show that $H^2(\mm_{g,P})$ injects
into
$H^2(\mm_{g-1,P\cup\{q,r\}})$ as soon as $g\ge 2$. This is the only
consequence of \Ref{injectindirr} that we will need.}
\endclaim

We are now in a position to describe the inductive step in the proof of
\Ref{H2}, in
genus 3 or more. As we announced, the cases of lower genus will be treated in
section 4.

Let then $g\ge 3$ be an integer, and let $P$ be a finite set.
If $P$ is not empty, let $p$ be a fixed element of $P$. Let $q,r$ be distinct
and not
belonging to $P$. Let $\xi:\mm_{g-1,P\cup\{q,r\}}\to \mm_{g,P}$ be the map that
is obtained by identifying the points labelled by $q$ and $r$. We wish to show
that $H^2(\mm_{g,P})$ is generated by tautological classes, assuming the
analogous
statement is known to hold for $\mm_{\gamma,\nu}$ whenever $\gamma<g$ or
$\gamma=g$ and $\nu<\vert P\vert$. We will do this only for
$P\neq\emptyset$, the argument for $P=\emptyset$ being entirely similar. Let
$y$ be any element of $H^2(\mm_{g,P})$. The pullback $\xi^*(y)$ is invariant
under the operation of interchanging $q$ and $r$. Therefore, by the induction
assumptions, it is a linear combination of $\kappa_1$, the $\psi_i$, $i\in P$,
$\psi_q+\psi_r$,
$\delta_{irr}$, and the classes $\delta_{u,U}$, $\delta_{u,U\cup\{q,r\}}$ and
$\delta_{u,U\cup\{q\}}+\delta_{u,U\cup\{r\}}$, where $u$ is any integer between
0 and $g$ and $U$ runs through all subsets of $P$ containing $p$; when $g=3$,
we can even do without $\kappa_1$. Formulas
\Ref{xipullbacks} tell us that there is a linear combination $z$ of
tautological
classes such that the pullback of $x=y-z$ is of the form
$$
\xi^*(x)=f(\psi_q+\psi_r)+\sum_{p\in U\subset P\atop 0\le u\le g-2}
g_{u,U}\delta_{u,U\cup\{q,r\}}+\sum_{p\in U\subset P\atop 0\le u\le g-1}
h_{u,U}(\delta_{u,U\cup\{q\}}+\delta_{u,U\cup\{r\}})\tal{pullbackofx}
$$
for suitable coefficients $f,g_{u,U},h_{u,U}$. In case $g=3$, we may even
assume,
using \Ref{kappag2}, that $f=0$. We will show that, in fact, $\xi^*(x)=0$;
\Ref{injectindirr} and \Ref{injectindirrrem} will then tell us that $x$ itself
vanishes,
proving that $y$ is a linear combination of tautological classes, as desired.

Suppose $s\not\in P\cup\{q,r\}$, and let $\vartheta:\mm_{g-1,P\cup\{s\}}\to
\mm_{g,P}$ be the map that is obtained by attaching a fixed elliptic tail at
the
point labelled by $s$. Look at the diagram
$$
\CD
\mm_{g-1, P\cup \{s\}} @()\L{\vartheta} @(1,0)@()\L{\varphi}@(0,-1)
& \mm_{g, P}\\
\mm_{g-1, P\cup\{q,r\}} @()\l{\xi}@(1,1)
\endCD
$$
where $\varphi$ attaches the point labelled by $t$ of a sphere marked by
$\{t,q,r\}$ to the point labelled by $s$ of a variable curve in $\mm_{g-1,
P\cup
\{s\}}$. This diagram is commutative up to homotopy. The identity
$\varphi^*\xi^*(x)=\vartheta^*(x)$, together with formulas \Ref{pullbackofx}
and
\Ref{thetapullbacks}, applied to $\varphi$, implies that
$$
\vartheta^*(x)=\sum_{p\in U\atop 0\le u\le g-2}
g_{u,U}\delta_{u,U\cup\{s\}}\,.\tal{thetastarx}
$$
Now consider the commutative diagram
$$
\CD
\mm_{g-2,P\cup\{q,r,s\}}@()\L{\beta}@(1,0)
@()\L{\gamma}@(0,-1)&\mm_{g-1,P\cup\{q,r\}}@()\L{\xi}@(0,-1)\\
\mm_{g-1,P\cup\{s\}}@()\L{\vartheta}@(1,0)&\mm_{g,P}
\endCD\tal{compatibility3}
$$
where $\gamma$ and $\beta$ are the analogues of $\xi$ and $\vartheta$,
respectively. If we write down explicitly the identity
$\gamma^*\vartheta^*(x)=\beta^*\xi^*(x)$ using formulas \Ref{pullbackofx},
\Ref{thetastarx}, \Ref{xipullbacks}, and \Ref{thetapullbacks}, we get a
relation
$$
\aligned
\sum_{p\in U\atop 0\le u\le g-2}
g_{u,U}(\delta_{u,U\cup\{s\}}&+\delta_{u-1,U\cup\{s,q,r\}})\\
&\hskip-40pt =f(\psi_q+\psi_r)+\sum_{p\in U\atop 0\le u\le g-2}
g_{u,U}(\delta_{u,U\cup\{q,r\}}+\delta_{u-1,U\cup\{s,q,r\}})\\
&\hskip-40pt\phantom{=} +\sum_{p\in U\atop 0\le u\le g-1}
h_{u,U}(\delta_{u,U\cup\{q\}}+\delta_{u,U\cup\{r\}}
+\delta_{u-1,U\cup\{s,q\}}+\delta_{u-1,U\cup\{s,r\}})
\endaligned\tal{relazione}
$$
in $H^2(\mm_{g-2,P\cup\{q,r,s\}})$. If $g\ge 4$, all the tautological classes
appearing in \Ref{relazione} are independent, so $f=g_{u,U}=h_{u,U}=0$ for all
$u$
and $U$. When $g=3$, we already know that $f=0$; since the boundary classes are
independent in genus 1, we conclude that $g_{u,U}=h_{u,U}=0$ for all $u$
and $U$ in this case as well. This shows that $\xi^*(x)=0$, as desired.

\sez{3}{The initial cases of the induction}

In this section we calculate those cohomology groups which are needed to start
the inductive proofs of \Ref{H1H3} and \Ref{H2}. More exactly, we shall compute
the
$k$-th cohomology group of $\mm_{g,n}$ for all $k$, $g$ and $n$ such that $k\le
3$ or
$k=5$ and $d(g,n)<k$. Our treatment will be elementary and self-contained for
$k\le 3$,
while for $k=5$ we shall use, directly or indirectly, some of the results of
[\Ref{Ezra}],
[\Ref{EzraM22}], and [\Ref{LooijM3}].

We have already settled the case $k=1$ in the body of the proof of
\Ref{H1H3}. For $k=2$, the values of $g$ and $n$ involved are $g=0$ and $n\le
5$, and
$g=1$ and $n\le 2$. For $k=3$ they are $g=0$ and $n\le 6$, $g=1$ and $n\le 3$,
and $g=2$ and $n\le 1$, while for $k=5$ they are $g=0$ and $n\le 8$, $g=1$ and
$n\le
5$, $g=2$ and $n\le 3$, and $g=3$ and $n\le 1$.

As we said, in genus zero we rely on Keel's results [\Ref{Keel}], although the
computations could be easily done directly. What Keel shows, among other
things, is
that $H^k(\mm_{0,n})$ vanishes for all odd $k$, and that $H^2(\mm_{0,n})$ is
generated
by tautological classes, modulo the relations described in \Ref{tautog0}, and
has
dimension $2^{n-1}-{n\choose 2}-1$. In the range we shall have to examine, that
is,
$3\le n\le 6$, the dimensions are $0,1,5,16$, respectively.

We now turn to higher genus. Since $\mm_{1,1}$ is
isomorphic to $\pip^1$, its second cohomology is one-dimensional (and generated
by
$\delta_{irr}$). We next show that there are surjective morphisms
$\alpha:\mm_{0,6}\to \mm_{2,0}$ and $\beta:\mm_{0,7}\to \mm_{2,1}$. This
implies
that the third cohomology groups of $\mm_{2,0}$ and $\mm_{2,1}$ vanish. Let
$(C;p_1,\dots,p_6)$ be a 6-pointed stable genus zero curve. The morphism
$\alpha$
associates to it the stable model of the double admissible covering of $C$
branched at
the $p_i$. As for $\beta$, the image under it of a 7-pointed stable genus zero
curve
$(C;p_1,\dots,p_7)$ is defined to be the stable model of $(D;q)$, where $D$ is
the double
admissible covering of $C$ branched at $p_1,\dots,p_6$, and $q$ is one of the
points
lying above $p_7$. Notice that it is immaterial which of the two possible
choices for $q$
we make, since they yield isomorphic 1-pointed curves.

To complete our analysis for $k=2,3$ it remains to compute the second
cohomology of
$\mm_{1,2}$ and the third cohomology of $\mm_{1,3}$. It will suffice to prove
the
following.

\clail{Lemma}{h3m13}$h^2(\mm_{1,2})=2$, $h^3(\mm_{1,3})=0$.
\endclaim

In fact, there are exactly two boundary classes in $H^2(\mm_{1,2})$, namely
$\delta_{irr}$ and $\delta_{1,\emptyset}$, which are independent by part {\sl
i)} of
\Ref{tautog>0}. Thus the first part of \Ref{h3m13} implies that $\delta_{irr}$
and
$\delta_{1,\emptyset}$ generate $H^2(\mm_{1,2})$.

The proof of the second part of \Ref{h3m13} relies on Theorem \Ref{H2}, and
hence also on the first part, in that we will use the fact that
$H^2(\mm_{1,3})$ is
freely generated by boundary classes. Since the boundary classes are
$\delta_{irr}$,
$\delta_{1,\emptyset}$,
$\delta_{1,\{1\}}$, $\delta_{1,\{2\}}$, and $\delta_{1,\{3\}}$, this shows in
particular
that $h^2(\mm_{1,3})=5$. As we know that $h^1(\mm_{1,2})=h^1(\mm_{1,3})=0$,
Poincar\'e duality implies that $\chi(\mm_{1,2})=2+h^2(\mm_{1,2})$ and
$\chi(\mm_{1,3})=12-h^3(\mm_{1,3})$. Lemma \Ref{h3m13} is
then a consequence of the following result.

\clail{Lemma}{chimm13}$\chi(\mm_{1,2})=4$, $\chi(\mm_{1,3})=12$.
\endclaim

The proof of the lemma rests on the following simple remark. For any space $X$
we
denote by $\chi_c(X)$ the Euler characteristic of $X$ with compact supports,
that is, the
alternating sum of the dimensions of the $\quq$-cohomology groups of $X$ with
compact supports. Now suppose $X$ is a quasi-projective algebraic variety, and
let
$$
X=\overline{X}_d\supset \overline{X}_{d-1}\supset\dots\supset
\overline{X}_1\supset\overline{X}_0
$$
be a filtration of $X$ by closed subvarieties. Suppose that
$X_i=\overline{X}_i\setminus\overline{X}_{i-1}$ is of pure dimension $i$ (or is
empty)
for every $i$. Then
$$
\chi_c(X)=\sum\chi_c(X_i)\,.\tal{chicomp}
$$
This can be proved by induction on $d$. There is nothing to prove when $d=0$.
Now
assume the result known for $\overline{X}_{d-1}$; thus
$$
\chi_c(\overline{X}_{d-1})=\sum_{i<d}\chi_c(X_i)\,.
$$
On the other hand the exact sequence of cohomology with compact supports
$$
\dots\to H_c^j(X_d)\to H_c^j(X)\to H_c^j(\overline{X}_{d-1})\to\cdots
$$
shows that
$$
\chi_c(X)=\chi_c(X_d)+\chi_c(\overline{X}_{d-1})
=\chi_c(X_d)+\sum_{i<d}\chi_c(X_i)
=\sum\chi_c(X_i)\,.
$$
There are two special cases when $\chi_c(X)=\chi(X)$. The first is obviously
the one
when $X$ is compact. The other is when $X$ is an orbifold; in this case
Poincar\'e
duality implies that $h_c^q(X)=h_{2d-q}(X)=h^{2d-q}(X)$, and hence
$$
\chi_c(X)=\sum (-1)^q h_c^q(X)=\sum (-1)^{2d-q} h^{2d-q}(X)=\chi(X)\,.
$$
Thus, when the $X_i$ are orbifolds, \Ref{chicomp} translates into
$$
\chi_c(X)=\sum\chi(X_i)\,.\tal{chisemi}
$$
If $X$ is compact, or an orbifold, this in turn yields
$$
\chi(X)=\sum\chi(X_i)\,.\tal{chi}
$$
We will first use these formulas to compute the Euler characteristics of the
spaces
$\m_{0,n}$ for $n\le 6$. Since $\m_{0,4}$ can be identified with the complex
plane
minus two points,
$$
\chi(\m_{0,4})=-1\,.\tal{chim04}
$$
We next deal with $\m_{0,5}$. This space can be identified with the complement,
inside the $\cic^2$ with coordinates $z_1,z_2$, of the lines $z_1=0$, $z_1=1$,
$z_2=0$,
$z_2=1$, $z_1=z_2$. Put otherwise, if $x_0,x_1,x_2$ are homogeneous coordinates
in
$\pip^2$, then $\m_{0,5}$ can be identified with the complement, inside
$\pip^2$, of the
six projective lines
$$
x_0=0,\ x_1=0,\ x_2=0,\ x_0=x_1,\ x_0=x_2,\ x_1=x_2.
$$
Each of these lines has exactly 3 points which are in common with some of the
others,
and there are 7 points where two or more lines meet. Applying \Ref{chi} to the
filtration $\text{7 points}\subset\text{union of 6 lines}\subset\pip^2$ then
gives
$$
3=\chi(\pip^2)=\chi(\m_{0,5})+6\chi(\m_{0,4})+7=\chi(\m_{0,5})-6+7\,,
$$
that is,
$$
\chi(\m_{0,5})=2\,.\tal{chim05}
$$
The same method can be applied to calculate the Euler characteristic of
$\m_{0,6}$.
This space can be identified with the complement inside $\pip^3$ of the 10
planes
$$
x_i=0,\quad 0\le i\le 3;\qquad
x_i=x_j,\quad 0\le i<j\le 3.
$$
Each nine of these cut on the remaining one a configuration of 6 lines which is
identical
to the one occurring in the analysis of $\m_{0,5}$ we just completed. In
addition, three
of these lines are common to three planes, and three only to two. Hence the
total
number of lines is 25.
There are exactly 15 points where three independent planes meet, namely the
points all of whose homogeneous coordinates are either 0 or 1. Applying formula
\Ref{chi} then gives
$$
\chi(\m_{0,6})=-6\,.\tal{chim06}
$$
We next deal with the genus 1 case. Since $\m_{1,1}$ can be identified with the
complex plane,
$$
\chi(\m_{1,1})=1\,.\tal{chim11}
$$
We will now show that
$$
\chi(\m_{1,2})=1\,.\tal{chim12}
$$
The proof of this fact and the remainder of the computation of the Euler
characteristics
of $\mm_{1,2}$ and $\mm_{1,3}$ require the calculation of the Euler
characteristics of
the quotients of some spaces $\m_{0,n}$ by certain group actions. We will
denote by
$\m'_{0,4}$ and
$\m'_{0,5}$ the quotients of $\m_{0,4}$ and $\m_{0,5}$ modulo the operation of
interchanging the labelling of two of the marked points, and by $\m''_{0,5}$
the
quotient of $\m_{0,5}$ modulo permutations of the labellings of three of the
marked
points. We claim that
$$
\chi(\m'_{0,4})=0\,;\quad \chi(\m'_{0,5})=1\,;\quad
\chi(\m''_{0,5})=1\,.\tal{chiquot}
$$
The map $\m_{0,4}\to\m'_{0,4}$ has degree 2. We claim that there is a unique
fiber
consisting of only one point, so that
$$
-1=\chi(\m_{0,4})=2\chi(\m'_{0,4})-1\,,
$$
proving the first identity in \Ref{chiquot}. Suppose in fact that there is an
isomorphism
between the two 4-pointed curves $(\pip^1;0,\infty,1,x)$ and
$(\pip^1;0,\infty,x,1)$.
This means that there is an automorphism $\alpha$ of $\pip^1$ such that
$\alpha(0)=0$, $\alpha(\infty)=\infty$, $\alpha(1)=x$, and $\alpha(x)=1$. The
first two
conditions imply that $\alpha$ is of the form $\alpha(z)=az$, for some nonzero
complex
number $a$. The last two conditions say that $x=a$ and $a^2=1$. Thus the curve
in
question, up to isomorphism, is $(\pip^1;0,\infty,1,-1)$. The same kind of
argument
proves the last identity in \Ref{chiquot}. In fact, let $x,y$ be distinct and
different
from $0,1,\infty$, and let $\alpha$ be a non-trivial automorphism of $\pip^1$
which
fixes $0$ and $\infty$ and permutes $1,x,y$; observe that none of these last
three
points can be fixed. Clearly, $\alpha$ is of the form $\alpha(z)=az$ for some
nonzero
complex number $a$. On the other hand, suppose for instance that $\alpha(1)=x$,
$\alpha(x)=y$ and $\alpha(y)=1$. Then $x=a$, $y=a^2$ and $a^3=1$. This shows
that the
degree 6 morphism $\m_{0,5}\to\m'_{0,5}$ has exactly one fiber which consists
of
fewer than 6 points, namely the fiber made up of $(\pip^1;0,\infty,1,a,b)$ and
$(\pip^1;0,\infty,1,b,a)$, where $a$ and $b$ are the two non-trivial cubic
roots of unity.
It follows that
$$
2=\chi(\m_{0,5})=6\chi(\m''_{0,5})-4\,,
$$
proving the third identity in \Ref{chiquot}. As for the second identity, the
morphism
$\m_{0,5}\to\m'_{0,5}$ is clearly unramified and of degree 2, so that
$\chi(\m'_{0,5})$
is one half of $\chi(\m_{0,5})$. This completes the proof of \Ref{chiquot}.

We now return to $\m_{1,2}$. Let $U$ be the open subset of $\m_{1,2}$
consisting of
those curves $(C;p_1,p_2)$ such that, if $\tau$ stands for the $-1$ involution
about the
origin $p_2$, then $\tau(p_1)\neq p_1$. We may associate to each 5-pointed
rational
curve $(\pip^1;q_1,\dots,q_5)$ the 2-pointed genus one curve $(C;p_1,p_2)$,
where $C$
is the double covering of $\pip^1$ branched at $q_2,\dots,q_5$, and $p_1$,
$p_2$ map
to $q_1$ and $q_2$, respectively (notice that the two possible choices of $p_2$
give
isomorphic 2-pointed curves). This defines a morphism from $\m_{0,5}$ to
$\m_{1,2}$, which clearly factors through an isomorphism between $\m''_{0,5}$
and
$U$. On the other hand the complement of $U$ in $\m_{1,2}$ can be identified
with
$\m'_{0,4}$, so that applying \Ref{chi} yields
$$
\chi(\m_{1,2})=\chi(\m''_{0,5})+\chi(\m'_{0,4})=1\,,
$$
which is just \Ref{chim12}.
To compute the Euler characteristic of $\mm_{1,2}$ we use the stratification by
graph
type. The open strata other than $\m_{1,2}$ are indexed by the graphs in Figure
1.
Here, and elsewhere, we adopt the convention that a solid dot stands for a
component
of genus zero, and a hollow one for a component of genus one.
\vskip5pt
\epsfxsize=148.9pt
\centerline{\epsfbox{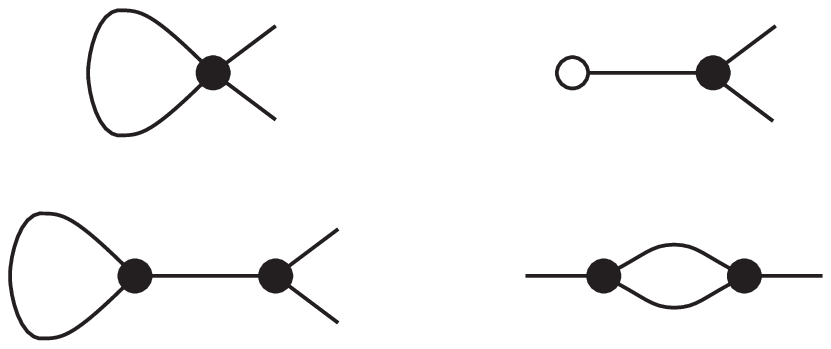}}
\vskip5pt
\centerline{Figure 1.}
\medskip\noindent
There are two one-dimensional strata $V_1$ and $V_2$, corresponding to
the first two graphs, while the zero-dimensional strata are the two points
corrsponding to the last two graphs.
It is evident that $V_1$ is isomorphic to $\m'_{0,4}$, and $V_2$ to $\m_{1,1}$.
Thus
$$
\chi(\mm_{1,2})=\chi(\m_{1,2})+\chi(\m'_{0,4})+\chi(\m_{1,1})+2
=\chi(\m_{1,2})+3=4\,.
$$
This proves the first statement in \Ref{chimm13}, and consequently also the
first one
in \Ref{h3m13}.

The next step is to compute the Euler characteristic of $\m_{1,3}$. We imitate
the
argument used for $\m_{1,2}$. Let
$U$ be the open subset of $\m_{1,3}$ consisting of those curves
$(C;p_1,p_2,p_3)$ such that, if
$\tau$ stands for the $-1$ involution about the origin $p_3$, then
$p_1\neq\tau(p_1)\neq p_2\neq \tau(p_2)$. Also, let $H$ be the subvariety of
$\m_{1,6}$ consisting of all curves $(C;p_1,p_2,p_3,p_4,p_5,p_6)$ such that, if
$\tau$ is as above, then $p_4$, $p_5$ and $p_6$ are fixed by $\tau$ and
$\tau(p_1)\neq p_2$. There are two dominant maps $\alpha:H\to \m_{0,6}$ and
$\beta:H\to U$. The first of them sends $(C;p_1,p_2,p_3,p_4,p_5,p_6)$ to
$(C/\tau;\overline{p}_1,\overline{p}_2,\overline{p}_3,
\overline{p}_4,\overline{p}_5,\overline{p}_6)$, where $\overline{p}_i$ stands
for the
image of $p_i$. The other map just ``forgets'' $p_4$, $p_5$ and $p_6$. The
map $\alpha$ is clearly unramified and of degree 2, while $\beta$ has degree 6.
We
claim that $\beta$ is also unramified. To see this, all we must show is that,
given any
element $(C;p_1,p_2,p_3)$ of $U$, its only automorphism is the identity. Any
automorphism $\sigma$ of $C$ fixing $p_3$ commutes with $\tau$, hence descends
to
an automorphism $\rho$ of $C/\tau=\pip^1$. If, in addition, $\sigma$ fixes
$p_1$ and
$p_2$, then $\rho$ is an automorphism of
$(C/\tau;\overline{p}_1,\overline{p}_2,\overline{p}_3)$.
Since this is a 3-pointed rational curve, $\rho$ is the identity, hence
$\sigma$ equals
$\tau$ or the identity. On the other hand, the definition
of $U$ says that $\tau$ is not an automorphism of $(C;p_1,p_2,p_3)$.

Having proved that $\beta$ is unramified, we can conclude that
$$
\chi(U)=\frac{2}{6}\chi(\m_{0,6})=-2\,.
$$
At this point, to compute the Euler characteristic of $\m_{1,3}$ we observe
that this moduli space is the union of $U$, of three two-dimensional strata
$U_1$,
$U_2$ and $U_3$, and a one-dimensional stratum $U_4$ which are defined as
follows.
The points of $U_1$ are the curves $(C;p_1,p_2,p_3)$ such that $\tau(p_1)=p_1$
and
$\tau(p_2)\neq p_2$, while $U_2$ is like $U_1$, but with the roles of 1 and 2
interchanged. The points of $U_3$ are the curves $(C;p_1,p_2,p_3)$ such that
$\tau(p_1)=p_2$, and those of $U_4$ are the curves $(C;p_1,p_2,p_3)$ such that
$\tau(p_1)=p_1$ and $\tau(p_2)=p_2$. Clearly, $U_1$ and $U_2$ are isomorphic to
$\m'_{0,5}$, $U_3$ is isomorphic to $\m''_{0,5}$ and $U_4$ to $\m_{0,4}$. Thus
$$
\chi(\m_{1,3})=\chi(U)+2\chi(\m'_{0,5})+\chi(\m''_{0,5})+\chi(\m_{0,4})
=-2+2+1-1=0\,.\tal{chim13}
$$
To compute the Euler characteristic of $\mm_{1,3}$ we use the stratification by
graph
type. The open strata other than $\m_{1,3}$ correspond to the graphs in Figure
2
(plus a labelling of the legs by $1,2,3$).
\vskip5pt
\epsfxsize=280.21pt
\centerline{\epsfbox{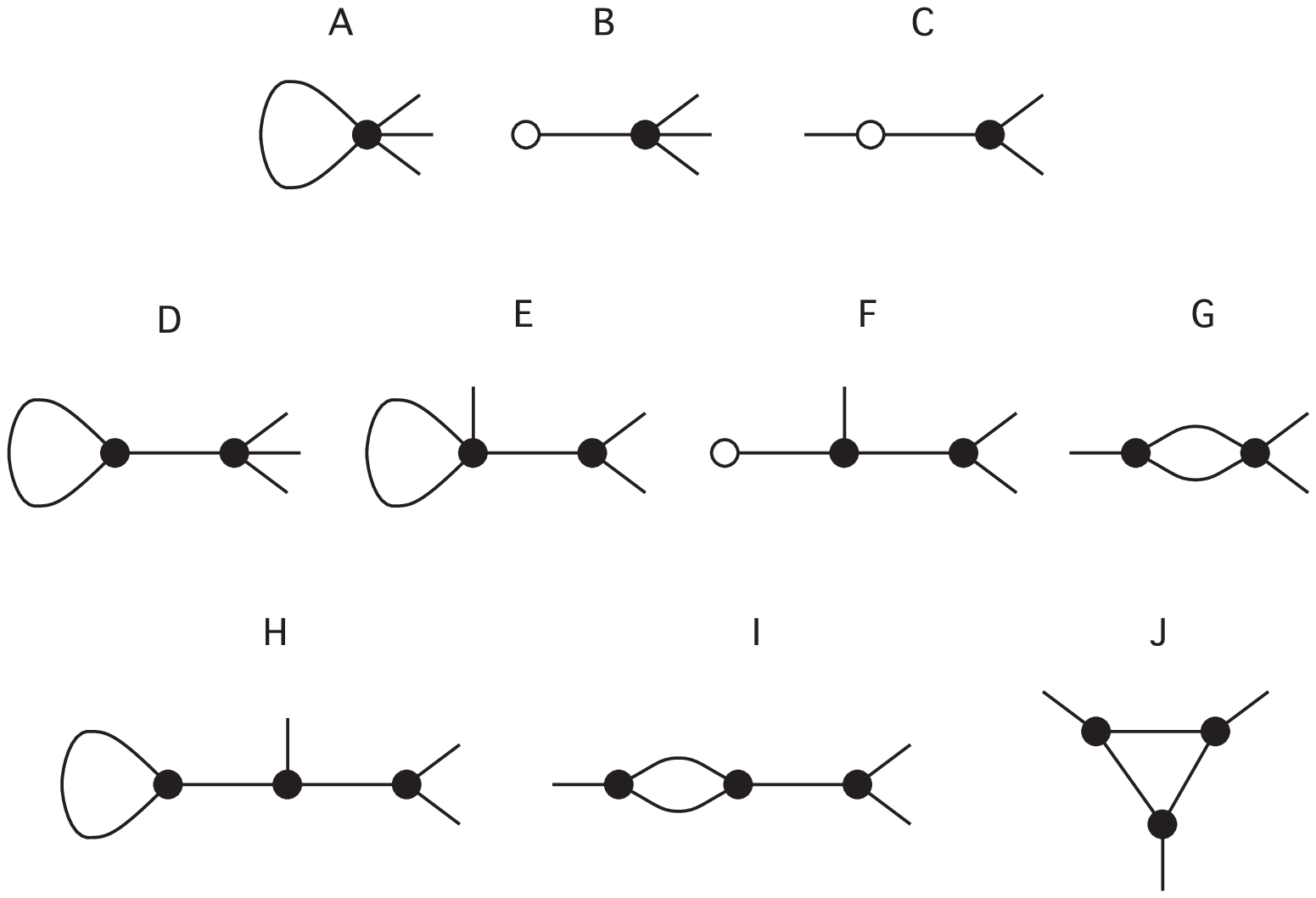}}
\vskip5pt
\centerline{Figure 2.}
\medskip\noindent
There are: one stratum corresponding to graph $A$, isomorphic to $\m'_{0,5}$,
one
corresponding to graph $B$, isomorphic to $\m_{1,1}\times\m_{0,4}$, three
strata
corresponding to graph $C$ and isomorphic to $\m_{1,2}$, one stratum
corresponding
to graph $D$, isomorphic to $\m_{0,4}$, three corresponding to graph $E$ and
isomorphic to $\m'_{0,4}$, three corresponding to graph $F$ and
isomorphic to $\m_{1,1}$, three corresponding to graph $G$ and
isomorphic to $\m'_{0,4}$, and seven zero-dimensional strata, all points, three
corresponding to graph $H$, three to graph $I$, and one to graph $J$. Putting
everything together we conclude that
$$
\aligned
\chi(\mm_{1,3})&=\chi(\m_{1,3})+\chi(\m'_{0,5})+\chi(\m_{1,1})\chi(\m_{0,4})
+3\chi(\m_{1,2})\\
&\quad+\chi(\m_{0,4})+6\chi(\m'_{0,4})+3\chi(\m_{1,1})+7=12\,,
\endaligned
$$
as desired. Theorem \Ref{H1H3} is now completely proved for $k=1,3$.
It remains to examine the initial cases of the induction for $k=5$, in positive
genus.
These are: $g=1$ and $n\le 5$, $g=2$ and $n\le 3$, $g=3$ and $n\le 1$. The
group
$H^5(\mm_{1,n})$ vanishes when $n\le 2$ for dimension reasons, and when $n\le
4$ by
Poincar\'e duality; in fact, $H^5(\mm_{1,3})$ and $H^5(\mm_{1,4})$ are
Poincar\'e dual
to $H^1(\mm_{1,3})$ and $H^3(\mm_{1,4})$. Likewise, $H^5(\mm_{2})$ and
$H^5(\mm_{2,1})$ are Poincar\'e dual to $H^1(\mm_{2})$ and $H^3(\mm_{2,1})$,
which
are both zero. On the other hand, Getzler has shown in [\Ref{Ezra}] that
$H^5(\mm_{1,5})$ vanishes, while in [\Ref{EzraM22}] he has proved that
$H^5(\mm_{2,2})$ is zero and announced that $H^5(\mm_{2,3})$ vanishes as well.
At
this point \Ref{oddindscheme} implies that $H^5(\mm_{g,n})$ vanishes for $g\le
2$ and
all $n$. In genus 3 we may argue as follows. Looijenga [\Ref{LooijM3}] proves
that
$H^7(\m_{3})$ and $H^9(\m_{3,1})$ are zero. By Poincar\'e duality, this is the
same as
saying that $H_c^5(\m_{3})$ and $H_c^5(\m_{3,1})$ vanish, so the exact sequence
of
cohomology with compact supports shows that $H^5(\mm_{3})$ and $H^5(\mm_{3,1})$
inject into $H^5(\del\m_{3})$ and $H^5(\del\m_{3,1})$, respectively. Lemma
\Ref{preinjectinbdry} then says that both $H^5(\mm_{3})$ and $H^5(\mm_{3,1})$
inject into sums $\bigoplus H^5(X_i)$, where the $X_i$ are products of moduli
spaces
$\mm_{g,n}$ with $g<3$. By what has already been proved, $H^5(X_i)=0$ for all
$i$,
hence $H^5(\mm_{3})=H^5(\mm_{3,1})=0$. This concludes the proof of \Ref{H1H3}.

\sez{4}{The induction step in low genus}

In this section we will complete the proof of \Ref{H2}; it remains to deal with
the
genus one and genus two cases.

\medskip
\noindent{\bf Genus 1.} We begin by improving on Lemma \Ref{deltairrtozero}.
Let
$P$ be a finite set, set $n=\vert P\vert $ and let $y$ and $z$ be distinct and
not
belonging to $P$. As usual, we let
$$
\xi:\mm_{0,P\cup\{y,z\}}\to\mm_{1,P}
$$
be the map gotten by identifying the points labelled by $y$ and $z$.

\clail{Lemma}{kerxi}The kernel of $\xi^*:H^2(\mm_{1,P})\to H^2(\mm_{0,P\cup\{y,
z\}})$ is one-dimensional and is generated by $\delta_{irr}$.
\endclaim
Lemma \Ref{deltairrtozero} and part {\sl i)} of Proposition \Ref{tautog>0} say
in
particular that
$\delta_{irr}$ is not zero and belongs to the kernel of $\xi^*$. It remains to
show
that any other element of
$\ker(\xi^*)$ is a multiple of
$\delta_{irr}$. The case $\vert P\vert=1$ is trivial. We have seen in the
previous
section that, when $\vert P\vert=2$, $H^2(\mm_{1,P})$ has dimension two. On the
other hand, the class $\delta_{1,\emptyset}\in H^2(\mm_{1,P})$ maps to
$\delta_{0,\{y,z\}}\in H^2(\mm_{0,P\cup\{y,z\}})$, which is not zero. This
takes care of
the case when $\vert P\vert=2$. We then proceed by induction on $n=\vert
P\vert$.
Let $\alpha$ be an element of the kernel of $\xi^*$. Let $s,t$ be distinct and
not
belonging to $P\cup\{y,z\}$. For any subset $S$ of $P$ with at most $n-2$
elements
consider the diagram
$$
\CD
\mm_{0,S\cup\{y,z,s\}}\times\mm_{0, S^c\cup\{t\}} @(1,0)
@()\L{\eta}@(0,-1) &\mm_{0,P\cup\{y,z\}}@()\L{\xi}@(0,-1)\\
\mm_{1,S\cup\{s\}}\times\mm_{0, S^c\cup\{t\}} @()\L{\nu_S}@(1,0)&\mm_{1,P}
\endCD
$$
where the vertical arrows are obtained by identifying the points labelled by
$y$
and $z$, and the horizontal ones by identifying the points labelled by $s$ e
$t$.
Using the K\"unneth formula and the vanishing of $H^1$ we may write
$\nu_S^*(\alpha)=(\beta,\gamma)$, where $\beta\in H^2(\mm_{1,S\cup\{s\}})$ and
$\gamma\in H^2(\mm_{0, S^c\cup\{t\}})$. From $\xi^*(\alpha)=0$ we deduce that
$\eta^*(\beta,\gamma)=0$. On the other hand
$\eta^*(\beta,\gamma)=(\beta',\gamma)$,
where $\beta'$ is the pullback of $\beta$ to $\mm_{0,S\cup\{y,z,s\}}$. It
follows that
$\gamma=0$ and, by induction hypothesis, that $\beta=a_S\delta_{irr}$ for a
suitable
constant $a_S$. In other words,
$$
\nu_S^*(\alpha)=(a_S\delta_{irr},0)\,.
$$
We now wish to show that, actually, $a=a_S$ does not depend on $S$. This will
conclude
the proof, since then the difference $\alpha-a\delta_{irr}$ will restrict to
zero on all
components of $\partial\mm_{1,P}$ and hence will be zero by \Ref{injectinbdry}.
To
show that $a_S$ is independent of $S$ we proceed as follows. If
$S\neq\emptyset$
write $S=T\cup\{w\}$, where $w\not\in T$, and consider the diagram
$$
\CD
\mm_{1,T\cup\{w\}}\times\mm_{0,\{w,z,s\}}\times\mm_{0, S^c\cup\{t\}}
@()\L{\tau}@(1,0) @()\L{\sigma}@(0,-1) &\mm_{1,T\cup\{y\}}\times\mm_{0,
S^c\cup\{w,r\}}@()\L{\nu_T}@(0,-1)\\ \mm_{1,S\cup\{s\}}\times\mm_{0,
S^c\cup\{t\}}
@()\L{\nu_S}@(1,0)&\mm_{1,P}
\endCD
$$
where the vertical arrows are obtained by identifying the points labelled by
$y$
and $z$, and the horizontal ones by identifying the points labelled by $s$ e
$t$. We find
that
$$
(a_T\delta_{irr},0,0)=\tau^*(a_T\delta_{irr},0)=\tau^*(\nu_T^*(\alpha))
=\sigma^*(\nu_S^*(\alpha))=\sigma^*(a_S\delta_{irr},0)=(a_S\delta_{irr},0,0)\,,
$$
and hence that $a_S=a_T$. Repeated applications of this argument show that
$a_S=a_\emptyset$ for any $S$. This proves the lemma.

\bigskip
Because of the relations between tautological classes
given in Proposition \Ref{tautog>0}, to prove
\Ref{H2} in genus one it suffices to show that $H^2(\mm_{1,P})$ is generated by
boundary classes. More precisely we will show  that the classes $\delta_{irr}$
and
$\delta_{1,S}$, where $S$ runs through all subset of $P$ with at most $n-2$
elements,
generate $H^2 (\mm_{1,P})$.
We prove this claim by induction on $n$. The first case $n=2$ has already been
checked. Denote by $V=V_{1,P}$ the subspace of $H^2 (\mm_{1,P})$ generated by
the
elements $\delta_{1,S}$, where $S$ runs through all subset of $P$ with at most
$n-2$
elements. In view of Lemma \Ref{deltairrtozero}, to prove our claim it suffices
to show
that the morphism $\xi^*$ vanishes modulo $V$. For this we shall use an
explicit basis
of $H^2 (\mm_{0, P\cup\{y,z\}})$ which we will presently describe.

\clail{Lemma}{standardbasis}Let $Q$ be a finite set with at least four
elements, and
let $x,y,z\in Q$ be distinct. Then $\delta_{0,\{y,z\}}$ and all the
classes $\delta_{0,S}$ such that $x\in S$ and $2\le\vert S\vert\le\vert
Q\vert-3$
constitute a basis of $H^2(\mm_{0,Q})$.
\endclaim
As we already recalled, Keel [\Ref{Keel}] shows that the dimension of
$H^2(\mm_{0,n})$
is $2^{n-1}-{n\choose 2}-1$. Now denote by $V$ the subspace spanned by the
classes
listed in the statement of the lemma. Since there are
$2^{\vert Q\vert-1}-{\vert Q\vert\choose 2}-1$ of these, it suffices to show
that all
classes $\delta_{0,T}$ belong to $V$. The only classes that are not already
present
in the list are those of the form $\delta_{0,\{a,b\}}$, where $x\not\in\{a,b\}$
and
$\{a,b\}\neq\{y,z\}$. Let $p,q,r$ be elements of $Q$, all different from $x$
and
such that
$p\neq q\neq r$. We claim that
$$
\delta_{0,\{q,r\}}\equiv \delta_{0,\{p,r\}} \mod V.\tal{cong}
$$
If $p=q$ there is nothing to prove. If $p\neq q$ then, in view of
\Ref{complrel}, one of
Keel's relations \Ref{Keelrel} is
$$
\delta_{0,\{x,p\}}+\sum_{x,p\in S\not\ni q,r\atop 3\le\vert S\vert\le\vert
Q\vert-3}
\delta_{0,S}+\delta_{0,\{q,r\}}=\delta_{0,\{x,q\}}
+\sum_{x,q\in S\not\ni p,r\atop
3\le\vert S\vert\le\vert Q\vert-3} \delta_{0,S}+\delta_{0,\{p,r\}}\,,
$$
which implies \Ref{cong}. One among $y$ and $z$, say $z$, is different from
both $a$
and $b$. Two applications of \Ref{cong} then give
$$
\delta_{0,\{a,b\}}\equiv \delta_{0,\{a,z\}}\equiv \delta_{0,\{y,z\}}\equiv
0\mod V,
$$
proving the lemma.

\bigskip
We will call a basis such as the one constructed in Lemma \Ref{standardbasis} a
{\it
canonical basis} (with respect to $x$, $y$ and $z$). To simplify notation, from
now on in
the genus zero case we shall write $\delta_S$ instead of $\delta_{0, S}$. Let
$Q=P\cup\{y,z\}$. Let $\Cal B$ be the canonical basis of
$H^2 (\mm_{0, Q})$ relative to $x,y,z$.
Let $\alpha\in
H^2 (\mm_{1,P})$. Using the fact that $\xi$ is invariant under the involution
exchanging $y$ and $z$, we can write $\xi^*\alpha$ in terms of $\Cal B$:
$$
\aligned
\xi^*\alpha&=a_{\{y,z\}}\delta_{\{y,z\}}+ \sum_{S\subset X, x\in S, \vert
S\vert \geq2,
\vert X\setminus S\vert \geq1}a_S\delta_S\\ &\quad+\sum_{S\subset X, x\in S,
\vert
X\setminus S\vert \geq3}b_S\delta_{S\cup\{y,z\}}+
\sum_{S\subset X, x\in S, \vert X\setminus S\vert
\geq2}c_S(\delta_{S\cup\{y\}}+\delta_{S\cup\{z\}})\,.
\endaligned\tal{m1n.2}
$$
\noindent
Consider a subset
$R$ of $P$ such that $\vert P\setminus R\vert \geq 2$, and look at the morphism
$$
\vartheta_R: \mm_{1,R\cup\{u\}}\to \mm_{1,P}
$$
defined by taking a varying stable, genus 1, $R\cup\{u\}$-pointed curve, a
fixed stable,
genus zero, $(P\setminus R)\cup\{v\}$-pointed curve $C_0$, and identifying the
points
labelled by $u$ and $v$. As we already noticed before stating Lemma
\Ref{thetapullbacks}, the homomorphism
$\vartheta^*_R$
does not depend on the choice of $C_0$.
By induction hypothesis, we have
$$
\vartheta_R^*(\alpha)=\sum_{S\subset R\,,\,\vert R\setminus S\vert
\geq1}f^R_S\,\delta_{1, S}+\sum_{S\subset R\,,\,\vert R\setminus S\vert
\geq2}g^R_S\,\delta_{1, S\cup\{u\}}+f^R\,\delta_{irr}\,.\tal{m1n.3}
$$
Consider the elements $\delta_{1, S}\in H^2(\mm_{1,P})$ with $\vert P\setminus
S\vert \geq 2$. Recalling the convention about the symbols $\delta_{a,A}$,
Lemma
\Ref{xipullbacks} says that
$$
\aligned
\xi^*\delta_{1, X\setminus \{p, q\}}&=\delta_{\{p, q\}}\notin \Cal
B,\qquad\qquad\ \
\text{if }p\neq x, q\neq x,\\
\xi^*\delta_{1, S}&=\left\{\matrix
\delta_{S\cup \{y, z\}}\in \Cal B,&&&
\text{if }
x\in S, S \neq P\setminus \{p, q\},p\neq x, q\neq x,\\ \delta_{P\setminus S}\in
\Cal
B,\hfill&&& \text{if } x\notin S.\hfill
\endmatrix
\right.
\endaligned\tal{m1n.4}
$$
So all the elements $\delta_{1, S}$ restrict to elements in $\Cal B$ except
when
$S=P\setminus \{p, q\}$. We also have, by lemma \Ref{thetapullbacks}, that
$$
\vartheta^*_R\delta_{1, S}=\left\{
\matrix
\delta_{1, S}\hfill&&&
\text{if }S\subset R,S\neq R,\hfill\\
-\psi_{u}=-\displaystyle{\sum_{S'\subset R, \vert R\setminus S'\vert \geq
1}}\delta_{1,S'}-\frac{1}{12}\delta_{irr}\hfill&&& \text{if }S=R,\hfill\\
\delta_{1, R\setminus (P\setminus
S)\cup\{u\}}\hfill&&&
\text{if }S\supset P\setminus R,\hfill\\
0\hfill&&& \text{otherwise}.\hfill
\endmatrix\right.
\tal{m1n.5}
$$
\clail{Remark}{rem}The only element $\delta_{1, S}$ with $S\subset P$ and
$\vert
P\setminus S\vert \geq 2$ such that, in the expression of
$\vartheta^*_R\delta_{1, S}$,
the element $\delta_{1, (R\setminus \{p, q\})\cup\{u\}}$ appears with non zero
coefficient, is $\delta_{1,P\setminus\{p, q\}}$ and in this case
$\vartheta^*_R\delta_{1,X\setminus\{p, q\}}=\delta_{1, R\setminus \{p,
q\}\cup\{u\}}$
\endclaim
Let us go back to the expression \Ref{m1n.3}. Let $\{r, s\}\subset R'\subset
R$. We
claim that
$$
g_{R\setminus\{r, s\}}^R=g_{R'\setminus\{r, s\}}^{R'}\,.\tal{m1n.6}
$$
Clearly, it suffices to prove the claim in case $R'=R\setminus \{q\}$. Look at
the diagram
$$
\CD
\mm_{1,R\cup\{u\}} @()\L{\xi_R}@(1,0) &
\mm_{1, P}\\
\mm_{1,(R\setminus \{q\}) \cup\{w\}}
@()\l{\xi_{R\setminus \{q\}}}@(1,1)@()\L{\varphi}@(0,1)
\endCD
$$
where the maps are defined in the obious way: we take the smooth rational curve
$C_0$
pointed by $(P\setminus R)\cup\{v\}$, which we used in the definition of
$\vartheta_R$, a smooth rational curve $C''_0$ pointed by $\{u, q, t\}$ and we
define
$C'_0$ by identifying the points labelled by $u$ and $v$ in $C''_0$ and $C_0$
respectively; then
$\varphi$ is defined by identifying the point labelled by $w$ with the point of
$C_0''$
labelled by $t$, while
$\vartheta_{R\setminus
\{q\}}$ is defined by identifying the point labelled by $w$ with the point of
$C_0'$
labelled by $t$. To prove the claim just use the
preceding remark with $P =R\cup\{u\}$, together with the fact that
$\varphi^*\vartheta^*_R\alpha=\vartheta^*_{R\setminus \{q\}}\alpha$.

\bigskip
After this preparation we are ready to modify $\alpha$. Suppose first that $
P=\{x,p,q\}$. The first move consists in adding to $\alpha$ a suitable multiple
of
$\delta_{1,\{\emptyset\}}$ so as to make $a_{\{y,z\}}=0$. The second move
consists in
adding to $\alpha$ a suitable multiple of $\delta_{1,\{x\}}$ so as to make
$f^{\{x\}}_{\emptyset}=0$. The third move consists in adding to $\alpha$ a
suitable
linear combination of $\delta_{1,\{p\}}$ and $\delta_{1,\{q\}}$ so as to make
$a_{\{x,p\}}=a_{\{x,q\}}=0$. The three moves, taken in that order, do not
interfere with
each other. As a result
$$
\aligned
\xi^*\alpha&=c_{\{x\}}(\delta_{\{x,y\}}+\delta_{\{x,z\}}),\\
\vartheta^*_{\{x\}}\alpha&=f^{\{x\}}\delta_{irr}
\endaligned\tal{m1n.7}
$$
Assume now that
$\vert P\vert \geq4$. Let $p\neq x$ and $q\neq y$. Observe that, by
\Ref{m1n.6},
whenever
$\vert P\setminus R\vert \geq 2$, $\vert P\setminus S\vert \geq 2$ and
$\{p,q\}\in
R\cap S$, we have
$$
g^{R}_{R\setminus\{p,q\}}=g^{R\cap
S}_{R\cap S\setminus\{p,q\}}=g^{S}_{S\setminus\{p,q\}}\,,
$$
so that $g^{R}_{R\setminus\{p,q\}}=\gamma_{p,q}$ does not depend on $R$.
Therefore
subtracting from $\alpha$ the class
$\sum_{X\supset\{p,q\}}\gamma_{p,q}\delta_{1,
\,X\setminus\{p,q\}}$, we get that $g^{R}_{R\setminus\{p,q\}}=0$, for all $p$,
$q$ and
$R$ such that $p\neq x$,
$q\neq x$, $\vert P\setminus R\vert \geq 2$ and $\{p,q\}\in R$.

\bigskip
The second move
consists in adding to $\alpha$ a linear combination of elements of type
$\delta_{1, S}$,
with $S\neq P\setminus \{p, q\}$, $p\neq x$, $q\neq x$, in such a way that
$$
\xi^*\alpha=
\sum_{x\in S\subset X, \vert X\setminus
S\vert \geq2}c_S(\delta_{S\cup\{y\}}+\delta_{S\cup\{z\}})\,.\tal{m1n.9}
$$
By the above remark, the second move does not alter what has been accomplished
by
the preceding one. For convenience we shall set $c_T=0$ when $x\not\in T$.

To prove our initial claim
we must prove that all the $c_S$ are equal to 0. Consider the square
$$
\CD
\mm_{0,R\cup\{y,z,u\}} @()\L{\eta_R} @(1,0) @()\L\eta
@(0,-1)&\mm_{0,P\cup\{y,z\}}
@()\L{\xi}@(0,-1)\\\mm_{1,R\cup\{u\}}@()\L{\xi_R}@(1,0) & \mm_{1,P}
\endCD
$$
where $\eta$ is the morphism obtained by identifying the points labelled by $y$
and
$z$, while
$\eta_R$ is obtained by identifying the point labelled by $u$ on the varying
curve in
$\mm_{0,R\cup\{y,z,u\}}$ with the point labelled by $v$ on the fixed curve
$C_0$. We
will examine the consequences of the equalities
$$
\eta_R^*\xi^*\alpha=\eta^*\xi_R^*\alpha\,,\qquad R\subset P,\ \vert P\setminus
R\vert
\geq 2\,.\tal{m1n.10}
$$
We have
$$
\eta^*\xi_R^*\alpha= \sum_{S\subset R,\vert R\setminus S\vert
\geq1}f^R_S\delta_{S\cup\{y,z\}}+\sum_{S\subset R,\vert R\setminus S\vert \geq
2}
g^R_S\delta_{ S\cup\{y,z,u\}}\,.\tal{m1n.11}
$$
Let us look at
$\eta_R^*\xi^*\alpha$. We have, for $S\subset P$, $\vert P\setminus S\vert
\geq2$,
$$
\eta_R^*\delta_{S\cup\{y\}}=
\left\{\matrix
\delta_{S\cup\{y\}}=
\delta_{(R\setminus S)\cup\{z,u\}}\hfill&&& \text{if }S\subset R\,,\hfill\\
\delta_{(P\setminus S)
\cup\{z\}}=\delta_{R\setminus(P\setminus S)\cup\{y,u\}}\hfill&&&
\text{if }P\setminus S\subset R\,,\hfill\\ 0\hfill&&&\text{otherwise,}\hfill
\endmatrix\right.
$$
and a similar relation holds for $\eta_R^*\delta_{S\cup\{z\}}$. Thus
$$
\aligned
\eta_R^*\xi^*\alpha=&
\sum_{S\subset R,\vert P\setminus S\vert \geq2} c_S(\delta_{ S\cup\{y\}}+
\delta_{S\cup\{z\}})\\ &+\sum_{P\setminus S\subset R,\vert P\setminus S\vert
\geq2}
c_S(\delta_{(P\setminus S)\cup\{z\}}+\delta_{(P\setminus S)\cup\{y\}}).
\endaligned\tal{m1n.12}
$$
We are going to prove that $c_T=0$ by descending induction on $\vert T\vert $.
The
first non-trivial case occurs when $\vert T\vert =n-2$.

If $P=\{x, p, q\}$ this amounts to showing that $c_{\{x\}}=0$. Looking at
\Ref{m1n.7},
relation \Ref{m1n.10} for $R=\{x\}$ says that
$$
0=c_{\{x\}}(\delta_{\{x,y\}}+\delta_{\{x,z\}})\,,
$$
and we are done in this case.

Suppose $\vert P\vert \geq 4$, and set $P\setminus \{p, q\}=T$. We can assume
that
$p\neq x$ and $q\neq x$, otherwise $c_T=0$. For $R=\{p, q\}$ the equality
\Ref{m1n.10} now reads
$$
\aligned
&f^{\{p, q\}}_{\emptyset}\delta_{\{y,z \}}+f^{\{p, q\}}_{\{p\}}\delta_{\{p, y,z
\}}+f^{\{p,
q\}}_{\{q\}}\delta_{\{q, y,z \}}
=\\ &c_{P\setminus {\{p, q\}}}
(\delta_{\{p, q, z \}}+\delta_{\{p, q, y \}})\,, \endaligned
$$
where on the left-hand side we used the fact that, by our first move,
$g^{\{p,q\}}_{\emptyset}=0$ while, to simplify the right-hand side, we used
again the
fact that
$c_S=0$ if $x\notin S$. The above is an equality in $H^2(\mm_{0,\{p, q, y,
z,u\}})$,
where as a canonical basis we take $\delta_{\{p, u\}}$, $\delta_{\{q,u\}}$,
$\delta_{\{ y, u\}}$,
$\delta_{\{z, u\}}$, $\delta_{\{ y, z\}}$. It follows that
$$
f^{\{p, q\}}_{\emptyset}=f^{\{p, q\}}_{\{p\}} =f^{\{p, q\}}_{\{q\}}
=c_{P\setminus {\{p, q\}}}=0\,.\tal{m1n.13}
$$
The first step in the induction is completed.

Now let $r\geq 3$, assume
that $c_{S}=0$ if $\vert P\setminus S\vert < r$, let $T$ be such that $\vert
P\setminus
T\vert =r$, and set $R=P\setminus T$. As usual we can assume that $x\notin R$.
Let us
first assume that $\vert T\vert =\vert P\setminus R\vert \geq 2$. Relation
\Ref{m1n.10} reads
$$
\sum_{S\subset
R,\vert R\setminus S\vert \geq1}\hskip-10pt f^R_S\delta_{(R\setminus
S)\cup\{u\}}
+\hskip-10pt\sum_{S\subset R,\vert R\setminus S\vert \geq2}\hskip-10pt
g^R_S\delta_{ S\cup\{y, z, u\}}=
c_T(\delta_{\{y, u\}}+\delta_{\{z, u\}})\,,
$$
where, to simplify the right-hand side, we used the inductive hypothesis
together with
the fact that $c_S=0$ if $x\notin S$. This is an equality in
$H^2(\mm_{0,R\cup\{y, z,
u\}})$. As a basis for $H^2 (\mm_{0,R\cup\{y,z,u\}})$ we take a canonical basis
where
now the role of $Q$ is played by $R\cup\{y, z, u\}$ and the role of $x$ is
played by $u$.
Look at the elements appearing in the above equality; since $x\notin R$, and
$g^R_{R\setminus \{p, q\}}=0$, these elements belong to the canonical basis and
they
are all distinct; so we are done.

Let us finally assume that $\vert T\vert =1$, or what is the same, that
$T=\{x\}$.
We start with a general remark.
Pulling back the class
$$
\xi^*_{\{p, q\}}\alpha=f^{\{p, q\}}_{\emptyset}\delta_{1,\emptyset }+ f^{\{p,
q\}}_{\{p\}}\delta_{1,\{p \}}+ f^{\{p, q\}}_{\{q\}}\delta_{1,\{q \}}+g^{\{p,
q\}}_{\emptyset}\delta_{1,\{u \}}+ f^{\{p,q\}}\delta_{irr}
$$
from $ H^2(\mm_{1,\{p, q, u\}})$ to $H^2(\mm_{1,\{q, w\}})$, comparing it with
$$
\xi^*_{q}\alpha=f^{\{ q\}}_{\emptyset}\delta_{1,\emptyset}+
f^{\{q\}}\delta_{irr}
$$
and looking at the coefficient of
$\delta_{1,\emptyset}$ we get
$$
f^{\{p, q\}}_{\emptyset}- f^{\{p,
q\}}_{\{q\}}=f^{\{q\}}_{\emptyset}\qquad\forall\ p,q. \tal{m1n.14}
$$
But if $\vert T\vert =1$ (and $\vert P\vert \geq4$), then relations
\Ref{m1n.13} have already been proved so that, in particular,
$f^{\{q\}}_{\emptyset}=0$ for $q\neq x$. Using \Ref{m1n.14} again we get
$$
f^{\{x, q\}}_{\emptyset}=f^{\{x, q\}}_{\{q\}}\,.\tal{m1n.15}
$$
Now look at relation \Ref{m1n.10} for $R=\{x, q\}$. Using the induction
hypothesis to
simplify the right-hand side one then gets
$$
\aligned
f^{\{x, q\}}_{\emptyset}\delta_{\{y,z \}}+f^{\{x, q\}}_{\{x\}}\delta_{\{x,
y,z\}}+
f^{\{x, q\}}_{\{q\}}\delta_{\{q, y,z\}}&+ g^{\{x,q\}}_{\emptyset}\delta_{\{x, q
\}}\\
&=c_x(\delta_{\{x, z\}}+\delta_{\{x, y \}})\,,
\endaligned
$$
which, by Keel's relations, can be written as
$$
\aligned
&f^{\{x, q\}}_{\emptyset}\delta_{\{y,z \}}+ f^{\{x, q\}}_{\{x\}}
\left(\delta_{\{x, q\}}+\delta_{\{x, u\}}-\delta_{\{x,y\}}- \delta_{\{x,
z\}}+\delta_{\{y,z\}}\right)\\&\hskip40pt +f^{\{x, q\}}_{\{q\}}\delta_{\{x,
u\}}+g^{\{x,
q\}}_{\emptyset}\delta_{\{x, q \}} =c_{\{x\}}(\delta_{\{x, z\}}+\delta_{\{x, y
\}})\,.
\endaligned
$$
Take
$$
\delta_{\{x, q\}}, \delta_{\{x, u\}}, \delta_{\{x, y\}}, \delta_{\{x, z\}},
\delta_{\{y, z\}}
$$
as a canonical basis for $H^2(\mm_{0,\{x, q, y, z, u \}})$. Looking at
the coefficient of $\delta_{\{y, z\}}$ we get $
f^{\{x,q\}}_{\emptyset}=-f^{\{x,q\}}_{\{q\}}$ which, together with
\Ref{m1n.15}, implies
that all coefficients in the above identity must vanish, concluding the proof
of the
lemma.

\medskip
\noindent{\bf Genus 2.}
In order to prove \Ref{H2} in genus two we must show that, for any finite set
$P$,
the space $H^2(\mm_{2, P})$ is generated by the classes $\psi_q$, with $q\in
P$, and by the boundary classes $\delta_{irr}$, $\delta_{1,A}$,
$\delta_{2,B}$, where $A$ and $B$ run through all subset of $P$ such that
$\vert
B^c\vert \geq 2$, and such that, if $P\neq\emptyset$, then $A$ contains a
preassigned
point $p\in P$. We set $n=\vert P\vert$.

We first consider the cases $n=0$ and $n=1$. By Theorem \Ref{injectindirr},
$H^2(\mm_{2})$ injects in $H^2(\mm_{1,2})$ via $\xi^*$. On the other hand,
$H^2(\mm_{1,2})$ is two-dimensional and the classes $\delta_{irr}$ and
$\delta_{1,\emptyset}$ are independent in $H^2(\mm_{2})$, so that $\xi^*$ is an
isomorphism, and the result follows in this case.

We next consider the case $n=1$. Again $H^2(\mm_{2,\{p\}})$ injects in
$H^2(\mm_{1,\{p,x,y\}})$ via $\xi^*$. Look at the diagrams
$$
\CD
\mm_{0,\{x,y,p,z\}} @()\L{\nu} @(1,0) @()\L\mu @(0,-1)&\mm_{1,\{p,z\}}
@()\L{\vartheta}@(0,-1)\\\mm_{1,\{p,x,y\}}@()\L{\xi}@(1,0) & \mm_{2, \{p\}}
\endCD
\tal{r6}
$$
$$
\CD
\mm_{1, \{p,z\}} @()\L{\vartheta} @(1,0)@()\L{\varphi}@(0,-1) & \mm_{2,
\{p\}}\\
\mm_{1, \{p,x,y\}} @()\l{\xi}@(1,1)
\endCD
\tal{r7}
$$
The map $\vartheta$ consists in attaching a fixed
one-pointed smooth elliptic curve $C_0$ to the point labelled by $z$ of a
variable
2-pointed elliptic curve. The map $\varphi$ consists in attaching
the point labelled by $w$ of a fixed smooth rational curve $C'_0$, marked by
the set
$\{x,y,w\}$, to the point labelled by
$z$ of a variable 2-pointed elliptic curve. The two diagrams are commutative up
to
homotopy.

Using Lemma \Ref{xipullbacks}, Lemma \Ref{thetapullbacks} and Proposition
\Ref{tautog>0}, we find that
$$
\aligned
\xi^*(\psi_p)&=\psi_p=\frac{1}{12}\delta_{irr}
+\delta_{1,\{x\}}+\delta_{1,\{y\}}+\delta_{1,\emptyset}\,,\\
\xi^*(\delta_{1,\emptyset})&=\delta_{1,\emptyset}+\delta_{1,\{p\}}\,,\\
\xi^*(\delta_{irr})&=\delta_{irr}-\psi_x-\psi_y+\delta_{1,\{x\}}
+\delta_{1,\{y\} }=
\frac{5}{6}\delta_{irr}-2\delta_{1,\{p\}}-2\delta_{1,\emptyset}\,,\\
\vartheta^*(\psi_p)&=\psi_p=\frac{1}{12}\delta_{irr}+\delta_{1,\emptyset}\,,\\
\vartheta^*(\delta_{1,\emptyset})&=\delta_{1,\emptyset}-
\psi_z=-\frac{1}{12}\delta_{irr}\,,\\
\vartheta^*(\delta_{irr})&=\delta_{irr}\,.\\
\endaligned
$$
A priori, given a class $\alpha$ in $H^2(\mm_{2,\{p\}})$, we have
$$
\aligned
\xi^*(\alpha)&= a\delta_{irr}+b(\delta_{1,\{x\}}+\delta_{1,\{y\}})
+c\delta_{1,\{p\}}+d\delta_{1,\emptyset},\\
\vartheta^*(\alpha)&=t\delta_{1,\emptyset}+s\delta_{irr},\\
\endaligned
$$
so that, by adding to $\alpha$ a suitable linear combination of $\psi_ p$,
$\delta_{irr}$,
and $\delta_{1,\{p\}}$, one can assume that
$s=b=d=0$. Then the equality $\nu^*\vartheta^*(\alpha)=\mu^*\xi^*(\alpha)$
gives
$c=t$, while the equality $\varphi^*\xi^*(\alpha)=\vartheta^*(\alpha)$ gives
$$
a\delta_{irr}-t\psi_z=a\delta_{irr}-
t(\frac{1}{12}\delta_{irr}+\delta_{1,\emptyset}) =t\delta_{1,\emptyset}\,,
$$
so that $a=t=c=0$, and the result is proved also in this case.

We now turn to the
general case of $\mm_{2, P}$ with $\vert P\vert =n\geq 2$, and proceed by
induction
on $n$. Fix, once and for all, a point $p$ in $P$. Consider subsets $R\subset
P$ such
that $p\in R^c$ and such that $R^c=P\setminus R$ contains two or more points.
We
look at the maps
$$
\aligned
&\xi:\mm_{1, P\cup\{x,y\}}\to \mm_{2, P}\,,\\
&\vartheta_R:\mm_{2, R\cup\{ z\}}\to\mm_{2, P}\,,\\
&\vartheta:\mm_{1, P\cup\{ s\}}\to \mm_{2, P}\,. \endaligned
$$
The map
$\vartheta_R$ is defined by identifying the point labelled by $w$ of a fixed
irreducible
smooth rational curve marked by the set $R^c\cup\{w\}$ with the point labelled
by
$z$ of a variable curve in $\mm_{2, R\cup\{ z\}}$. The map $\vartheta$ is
defined by
identifying the point labelled by
$t$ of a fixed 1-pointed elliptic curve with the point labelled by $s$ of a
variable
curve in $\mm_{1, P\cup\{ s\}}$. Given a class $\alpha\in H^2(\mm_{2, P})$, a
priori
one has
$$
\aligned
\xi^*(\alpha)&=a\delta_{irr}+\sum_{S\subset P}a_S\delta_{1,S} +\sum_{S\subset
P,
\vert S^c\vert \geq2}b_S\delta_{1,S\cup\{x,y\}} +\sum_{S\subset P, \vert
S^c\vert
\geq1}c_S(\delta_{1,S\cup\{x\}}+\delta_{1,S\cup\{y\}})\,,\\
\vartheta_R^*(\alpha)&=a^R\delta_{irr}+\sum_{r\in R}b_r^R\psi_r+b^R_z\psi_z+
\hskip-4pt\sum_{S\subset R, \vert S\vert \geq 2}\hskip-4pt c^R_S
\delta_{0,S}+\hskip-4pt\sum_{S\subset R, \vert S\vert \geq 1}\hskip-4pt
d^R_S\delta_{0,S\cup\{z\}}+ \hskip-4pt\sum_{S\subset R}\hskip-4pt
h^R_S\delta_{1,S}\,.
\endaligned
$$
Let us show that, by adding to $\alpha$ a suitable linear combination of
$\delta_{irr}$
and of the $\psi_i$ with $i\neq p$, one can assume that $a^R=b_r^R=0$ for all
$r\in R$
and for all $R\subset P$ such that $\vert R^c\vert \geq 2$. For each proper
subset
$R'$ of $R$ there is an obvious diagram
$$
\CD
\mm_{2, R'\cup\{z'\}} @()\L{\vartheta_{R'}} @(1,0)@()\L{\zeta}@(0,-1) &\mm_{2,
P}\\
\mm_{2,R\cup\{z\}} @()\l{\vartheta_R}@(1,1)
\endCD\tal{r8}
$$
which is commutative up to homotopy,
and it is evident that $ a^R=a^{R'}$ and $b_r^R=b_r^{R'}$ whenever $r\in R'$.
But
then it suffices to annihilate $a^{\emptyset}$ and $b_r^{\{r\}}$, for every
$r\in
P\setminus\{p\}$, which can be achieved by adding to $\alpha$ a suitable linear
combination of $\delta_{irr}$ and $\psi_r$ for $r\in P\setminus\{p\}$. Next,
recall that
$$
\xi^*(\psi_p)=\frac{1}{12}\delta_{irr}+
\sum_{T\subset P\cup\{x,y\},p\in P}\delta_{0,T}\,,
$$
so that, by adding to $\alpha$
a suitable multiple of $\psi_p$, we can assume that $a=0$. Finally, since
$$
\aligned
\xi^*(\delta_{2,S})&=\delta_{1,S\cup\{x,y\}}\,,\\
\xi^*(\delta_{1,S})&=\delta_{1,S}+\delta_{0,S\cup\{x,y\}}
=\delta_{1,S}+\delta_{1,S^c}\,,\\
\endaligned
$$
we can assume that $b_S=0$ and that
$a_S=0$ if $p\notin S$. In conclusion
$$
\xi^*(\alpha)=\sum_{S\subset P,p\in S}a_S\delta_{1,S} +\sum_{S\subset P, \vert
S^c\vert \geq1}c_S(\delta_{1,S\cup\{x\}}+\delta_{1,S\cup\{y\}})\,, \tal{r9}
$$
$$
\vartheta_R^*(\alpha)=b^R_z\psi_z+\hskip-8pt
\sum_{S\subset R, \vert S\vert \geq 2}\hskip-8pt
c^R_S\delta_{0,S}+ \hskip-8pt\sum_{S\subset R, \vert
S\vert \geq 1}\hskip-8pt
d^R_S\delta_{0,S\cup\{z\}}+ \sum_{S\subset
R}h^R_S\delta_{1,S}\,.\tal{r10}
$$
Now set $R=P\setminus\{p,q\}$, for some
$q\in P$ with $q\neq p$, and look at the diagram
$$
\CD
\mm_{1, R\cup\{x,y,z\}} @()\L{\nu} @(1,0) @()\L\mu @(0,-1)&\mm_{1,
P\cup\{x,y\}}
@()\L{\xi}@(0,-1)\\\mm_{2, R\cup\{z\}}@()\L{\vartheta_R}@(1,0) & \mm_{2, P}
\endCD
$$
Using Lemma (1.4), we get
$$
\aligned
\nu^*\xi^*(\alpha)&=\sum_{ \{p,q\}\subset S\subset
P}a_S\delta_{1,(S\setminus\{p,q\})\cup\{z\}} +\sum_{S\subset R}
c_S(\delta_{1,S\cup\{x\}}+\delta_{1,S\cup\{y\}})\\
&\phantom{=}+\sum_{\{p,q\}\subset S\subset
P,\vert S^c\vert \geq1}c_S(\delta_{1,(S\setminus\{p,q\})\cup\{x,z\}}+
\delta_{1,(S\setminus\{p,q\})\cup\{y,z\}})\,.
\endaligned\tal{r12}
$$
$$
\aligned
\mu^*\vartheta_R^*(\alpha)&=b^R_z\psi_z+
\sum_{S\subset R, \vert S\vert \geq 2}c^R_S\delta_{1,(R\setminus
S)\cup\{x,y,z\}}\\
&\phantom{=}+
\sum_{S\subset R, \vert S\vert \geq 1}d^R_S\delta_{1,(R\setminus S)\cup\{x,y\}}
+\sum_{S\subset R}h^R_S(\delta_{1,S}+\delta_{1,(R\setminus
S)\cup\{z\}})\,,
\endaligned\tal{r13}
$$
The equality $\nu^*\xi^*(\alpha)=\mu^*\vartheta_R^*(\alpha)$ provides a
relation of
linear dependence between elements of the following types
$$
\aligned
&\delta_{1,T\cup\{z\}},\ \delta_{1,T\cup\{x\}},\ \delta_{1,T\cup\{y\}},\
\delta_{1,T\cup\{x,z\}},\
\delta_{1,T\cup\{y,z\}},\ \\
&\quad\psi_z,\
\delta_{1,T\cup\{x,y,z\}},\
\delta_{1,T\cup\{x,y\}},\
\delta_{1,T}+\delta_{1,R\setminus T\cup\{z\}},\
\endaligned
$$
where $T\subset R$. By the results on $H^2(\mm_{1,\nu})$, and since
$\psi_z=(1/12)\delta_{irr}+\cdots$, the above elements are linearly
independent, so
that all the coefficient in \Ref{r12} and \Ref{r13} are zero. As $q$ is any
element in $P$
different from $p$, this means that \Ref{r9} becomes
$$
\aligned
\xi^*(\alpha)&=a_{\{p\}}\delta_{1,\{p\}}
+c_{\{p\}}(\delta_{1,\{p,x\}}+\delta_{1,\{p,y\}})\\
&\phantom{=}+c_{(P\setminus\{p\})}
(\delta_{1,(P\setminus\{p\})\cup\{x\}}
+\delta_{1,(P\setminus\{p\})\cup\{y\}})\,,
\endaligned\tal{r14}
$$
while $\vartheta_R^*(\alpha)=0$.
We now look at the following diagram, which is the analogue of \Ref{r7}.
$$
\CD
\mm_{1, P\cup\{s\}} @()\L{\vartheta} @(1,0)@()\L{\varphi}@(0,-1) & \mm_{2, P}
\\\mm_{1, P\cup\{x,y\}} @()\l{\xi}@(1,1) \endCD
$$
The identity $\varphi^*\xi^*(\alpha)=\vartheta^*(\alpha)$, together with
\Ref{thetapullbacks}, applied to $\varphi$, gives
$$
\vartheta^*(\alpha)= a_{\{p\}}\delta_{1,\{p\}}\,.\tal{r11}
$$
Finally, we consider the diagram
$$
\CD
\mm_{0, P\cup\{x,y,s\})} @()\L{A} @(1,0) @()\L{B} @(0,-1)&\mm_{1, P\cup\{x,y\}}
@()\L{\xi}@(0,-1)\\\mm_{1, P\cup\{s\}}@()\L{\vartheta}@(1,0) & \mm_{2, P}
\endCD
$$
The identity $B^*\vartheta^*(\alpha)=A^*\xi^*(\alpha)$ gives
$$
a_{\{p\}}\delta_{\{p,x,y\}}=a_{\{p\}}\delta_{\{p,s\}}
+c_{\{p\}}(\delta_{\{p,x,s\}}+
\delta_{\{p,y,s\}})+c_{(P\setminus \{p\})}(\delta_{(P\setminus
\{p\})\cup\{x,s\}}
+\delta_{(P\setminus \{p\})\cup\{y,s\}})\,.
$$
As long as $\vert P\vert \geq 3$, the boundary classes appearing in the above
relation
belong to the canonical basis of $H^2(\mm_{0, P\cup\{x, y, s\}})$, relative to
$p,x,y$, so
that all coefficients must vanish and we are done. If
$P=\{p,q\}$, we further simplify notation and rewrite the above relation as
$$
a\delta_{\{q,s\}}=a\delta_{\{p,s\}}+c(\delta_{\{q,y\}}+
\delta_{\{q,x\}})+d(\delta_{\{p,y\}}+
\delta_{\{p,x\}})\,.
$$
We choose $\delta_{\{q,s\}}$, $\delta_{\{q,x\}}$, $\delta_{\{q,y\}}$,
$\delta_{\{q,p\}}$
and $\delta_{\{p,s\}}$ as a basis for $\mm_{0,5}$ and, using the relations
$$
\aligned
\delta_{\{p,y\}}&=\delta_{\{p,s\}}+\delta_{\{q,y\}}-\delta_{\{q,s\}}\,,\\
\delta_{\{p,x\}}&=\delta_{\{p,s\}}+\delta_{\{q,x\}}-\delta_{\{q,s\}}\,,
\endaligned
$$
we get $a=-2d$ and $c=-d$. Thus
$$
\xi^*(\alpha)
=c(2\delta_{1,\{p\}}+\delta_{1,\{p,x\}}+\delta_{1,\{p,y\}}
-\delta_{1,\{q,x\}}-\delta_{1,\{q,y\}})\,.
$$
On the other hand, using \Ref{xipullbacks} and \Ref{psig1}, one finds that
$$
2\delta_{1,\{p\}}+\delta_{1,\{p,x\}}+\delta_{1,\{p,y\}}
-\delta_{1,\{q,x\}}-\delta_{1,\{q,y\}}=\xi^*(\delta_{1,\{p\}}
+\psi_q-\psi_p)\,,
$$
so \Ref{injectindirr} implies that $\alpha=c(\delta_{1,\{p\}}+\psi_q-\psi_p)$.
The
proof of \Ref{H2} is now complete.

\newcounter\bibno
\makebib
\ininbook

\bib
\no\bibno\label{ArbCo}
\by Enrico Arbarello, Maurizio Cornalba
\paper Combinatorial and algebro-geometric
cohomology classes on the moduli spaces of curves
\jour J.~Alg.~Geom.
\vol 5
\yr 1996
\pp 705--749
\endbib

\bib
\no\bibno\label{BoPik}
\paper Galois covers of moduli of curves
\by Marco Boggi, Martin Pikaart
\paperinfo preprint 1997
\endbib

\bib
\no\bibno\label{MDT}
\paper On the projectivity of the moduli spaces of
curves
\by Maurizio Cornalba
\jour J.~reine angew. Math.
\pp 11--20
\vol 443
\yr 1993
\endbib

\bib
\no\bibno\label{HodgeIII}
\paper Th\'eorie de Hodge III
\by Pierre Deligne
\jour I.H.E.S.~Publ.~Math.
\pp 5--77
\vol 44
\yr 1974
\endbib

\bib
\no\bibno\label{Faber}
\paper Chow rings of moduli spaces of curves
\by Carel Faber
\paperinfo thesis, Universiteit van Amsterdam, 1988
\endbib

\bib
\no\bibno\label{Ezra}
\paper Intersection theory on $\mm_{1,4}$ and elliptic Gromov-Witten invariants
\by Ezra Getzler
\jour J.~Amer.~Math.~Soc.
\pp 973--998
\vol 10
\yr 1997
\endbib

\bib
\no\bibno\label{EzraM22}
\paper Topological recursion relations in genus 2
\by Ezra Getzler
\paperinfo {\tt math.AG/9801003}
\endbib

\bib
\no\bibno\label{HarerH2}
\paper The second homology group of the mapping
class group of an orientable surface
\by John Harer
\jour Inv.~Math.
\pp 221--239
\vol 72
\yr 1982
\endbib

\bib
\no\bibno\label{HarerVirtDim}
\paper The virtual cohomological dimension of the mapping
class group of an orientable surface
\by John Harer
\jour Inv.~Math.
\pp 157--176
\vol 84
\yr 1986
\endbib

\bib
\no\bibno\label{HarerH3}
\paper The third homology group of the moduli space of curves
\by John Harer
\jour Duke Math.~J.
\pp 25--55
\vol 65
\yr 1991
\endbib

\bib
\no\bibno\label{HarerH4}
\paper The fourth homology group of the moduli space of curves
\by John Harer
\paperinfo to appear
\endbib

\bib
\no\bibno\label{HarrisMumf}
\by Joseph Harris, David Mumford
\paper On the Kodaira dimension of the moduli space of curves
\jour Inv.~Math.
\vol 67
\yr 1982
\pp 23--86
\endbib

\bib
\no\bibno\label{Keel}
\paper Intersection theory of moduli space of stable N-pointed curves of genus
zero
\by Se\`an Keel
\jour Trans.~AMS
\pp 545--574
\vol 330
\yr 1992
\endbib

\bib
\no\bibno\label{LooijM3}
\paper Cohomology of $\m_3$ and $\m^1_3$
\by Eduard Looijenga
\inbook ``Mapping class groups and moduli spaces of Riemann surfaces''
\eds C.-F.~B\"odigheimer and R.~M.~Hain
\bookinfo Contemp.~Math.~150
\pp 205--228
\publ Amer.~Math.~Soc.
\publaddr Providence, RI
\yr 1993
\endbib

\bib
\no\bibno\label{Looij2}
\paper Smooth Deligne-Mumford compactifications by means of Prym level
structures
\by Eduard Looijenga
\jour J.~Alg.~Geom.
\pp 283--293
\vol 3
\yr 1994
\endbib

\bib
\no\bibno\label{Looij}
\paper Cellular decompositions of compactified moduli spaces of pointed curves
\by Eduard Looijenga
\inbook ``The Moduli Space of Curves''
\eds R.~ Dijkgraaf, C.~Faber, G.~van der Geer
\bookinfo Progress in Mathematics 129
\pp 369--399
\publ Birkh\"auser
\publaddr Boston
\yr 1995
\endbib

\bib
\no\bibno\label{Mumfenseign}
\by David Mumford
\paper Stability of projective varieties
\jour L'Ens.~Math.
\vol 23
\yr 1977
\pp 39--110
\endbib

\bib
\no\bibno\label{Mumfenum}
\paper Towards an enumerative geometry of the
moduli space of curves
\inbook ``Arithmetic and Geometry''
\eds M.~Artin, J.~Tate
\bookinfo vol.~2
\publ Birkh\"auser
\publaddr Boston
\by David Mumford
\pp 271--328
\yr 1983
\endbib

\endmakebib

\bigskip\bigskip
Enrico Arbarello:

Scuola Normale Superiore,
Piazza dei Cavalieri 7,
56126 Pisa, Italia

e-mail: ea\@sabsns.sns.it
\medskip

Maurizio Cornalba:

Dipartimento di Matematica,
Universit\`a di Pavia,
via Ferrata 1,
27100 Pavia, Italia

e-mail: cornalba\@unipv.it

\enddocument